\documentclass[10pt]{article}
\usepackage{amssymb,latexsym}
\usepackage{epsfig}
\usepackage{eufrak}
\usepackage{amsmath}
\usepackage{mathrsfs}
\usepackage{color}

\newtheorem{theorem}{Theorem}
\newtheorem{lemma}{Lemma}

\newcommand{\be}{\begin{equation}}
\newcommand{\ee}{\end{equation}}
\newcommand{\bee}{\begin{eqnarray*}}
\newcommand{\eee}{\end{eqnarray*}}
\newcommand{\bel}{\begin{eqnarray}}
\newcommand{\eel}{\end{eqnarray}}
\newcommand{\bec}{\begin{cases}}
\newcommand{\eec}{\end{cases}}
\newcommand{\bem}{\begin{bmatrix}}
\newcommand{\eem}{\end{bmatrix}}

\newcommand{\la}{\label}
\newcommand{\li}{\left}
\newcommand{\ri}{\right}

\newcommand{\ovl}{\overline}
\newcommand{\udl}{\underline}

\newcommand{\lc}{\lceil}
\newcommand{\rc}{\rceil}

\newcommand{\ep}{\epsilon}
\newcommand{\vep}{\varepsilon}

\newcommand{\Lm}{\Lambda}
\newcommand{\Up}{\Upsilon}

\newcommand{\vsi}{\varsigma}
\newcommand{\de}{\delta}

\newcommand{\ga}{\gamma}

\newcommand{\vse}{\vartheta}
\newcommand{\se}{\theta}
\newcommand{\Se}{\Theta}

\newcommand{\al}{\alpha}
\newcommand{\ba}{\beta}

\newcommand{\vro}{\varrho}
\newcommand{\ro}{\rho}
\newcommand{\ka}{\kappa}
\newcommand{\om}{\omega}

\newcommand{\f}{\frac}
\newcommand{\sq}{\sqrt}
\newcommand{\cd}{\cdots}

\newcommand{\qu}{\quad}
\newcommand{\qqu}{\qquad}
\newcommand{\fa}{\forall}

\newcommand{\mscr}{\mathscr}
\newcommand{\mcal}{\mathcal}
\newcommand{\mbf}{\mathbf}
\newcommand{\bb}{\mathbb}

\newcommand{\wh}{\widehat}

\newcommand{\bs}{\boldsymbol}

\newcommand{\ap}{\approx}

\newcommand{\sh}{\slash}

\newcommand{\tx}{\text}

\newcommand{\iy}{\infty}

\newcommand{\jm}{\jmath}
\newcommand{\pa}{\partial}

\newcommand{\bed}{\begin{description}}
\newcommand{\eed}{\end{description}}
\newcommand{\bei}{\begin{itemize}}
\newcommand{\eei}{\end{itemize}}
\newcommand{\ben}{\begin{enumerate}}
\newcommand{\een}{\end{enumerate}}
\newcommand{\bib}{\bibitem}
\newcommand{\beL}{\begin{lemma}}
\newcommand{\eeL}{\end{lemma}}
\newcommand{\beT}{\begin{theorem}}
\newcommand{\eeT}{\end{theorem}}
\newcommand{\sect}{\section}

\newcommand{\bpf}{\begin{pf}}
\newcommand{\epf}{\end{pf}}
\newcommand{\bsk}{\bigskip}

\newcommand{\pfbox}{\hfill\mbox{$\Box$}}

\newenvironment{pf}{\paragraph*{Proof{\rm.}}}{\pfbox\bigskip}

\begin{document}

\title{{\bf Asymptotically  Optimal Sequential Estimation of the Mean Based on Inclusion Principle}
\thanks{The author have been  working with Department of Electrical Engineering at Louisiana
State University at Baton Rouge, LA 70803, USA, and  Department of Electrical Engineering at Southern University and A\&M College, Baton Rouge,
LA 70813, USA; Email: chenxinjia@gmail.com. The main results of this paper appeared in Proceedings of SPIE Conferences.} }

\author{Xinjia Chen}

\date{June 2013}

\maketitle

\begin{abstract}

A large class of problems in sciences and engineering can be formulated as the general problem of constructing random intervals with
pre-specified coverage probabilities for the mean.  Wee propose a general approach for statistical inference of mean values based on accumulated
observational data.  We  show that the construction of such random intervals can be accomplished by comparing the endpoints of random intervals
with confidence sequences for the mean. Asymptotic results are obtained for such sequential methods.

\end{abstract}

\tableofcontents

\section{Introduction}

As information technology is pervasive to our society,  statistical inferential methods are becoming increasingly important for the purpose of
harnessing the power of information.  In particular, a wide variety of machine learning problems can be cast into the general framework of
parameter estimation.  A frequent problem is the estimation of the mean value of a random variable. Familiar examples include, PAC learning and
algorithmic learning \cite{Mitchell, Vapnik, vind}, statistical pattern recognition \cite{Fu, Kostiantis, Vapnik}, data mining based on huge
data sets \cite{John, Kiv}, discovery of association rules \cite{Agra, Mitchell}, quiry size estimation \cite{Has, Liptona, Liptonb}, to name
but a few. An important issue of parameter estimation is the determination of sample sizes. However, the appropriate sample size usually depends
on the parameters to be estimated from the sampling process. To overcome this difficulty, sequential estimation methods have been proposed in
the literature (see, e.g. \cite{Chow1965, Lai, Gosh, Wald} and the references therein).  In a sequential method, the sample size is not fixed in
advance. Instead, data is evaluated as it is collected and further sampling is stopped in accordance with a pre-defined stopping rule as
significant results are observed.  Despite the richness of sequential methods, it seems that there exists no unified approach for sequential
estimation which deals with various precision requirements and make full use of the parametric information available. In this respect, we
propose to apply the inclusion principle proposed in \cite{Chenestimation, Chen_rule} for constructing multistage procedures for parameter
estimation.  In our paper \cite{Chenestimation, Chen_rule}, we have demonstrated that a wide variety of sequential estimation problems can be
cast into the general framework of constructing a sequential random interval of a prescribed level of coverage probability. To ensure the
requirement of coverage probability, we propose to use a sequence of confidence intervals, referred to as controlling confidence sequence, to
define a stopping rule such that  the sampling process is continued until the controlling confidence sequence is included by the sequential
random interval at some stage.  Due to the inclusion relationship, the associated coverage probability can be controlled by the confidence
coefficients of the sequence of confidence intervals. Moreover, as will be seen in the sequel, stopping rules derived from the inclusion
principle usually achieve best possible efficiency.

The remainder of this paper is organized as follows.  In Section 2, we demonstrate that a lot of estimation problems can be formulated as the
construction of an asymptotically symmetrical random interval for the mean of a random variable. In Section 3, we apply the inclusion principle
to construct fully sequential sampling schemes for estimation of mean values.  In Section 4, we apply the inclusion principle to construct
multistage sequential sampling schemes for estimation of mean values.  Section 5 is the conclusion.  The statistical methodology proposed in
this paper has been applied to electrical engineering, computer science and medical sciences, see recent literature \cite{Chen_SPIE2013,
Chen_SPIE8387, Chen_SPIE, Chen_SPIE11a, Chen_SPIE11b, Chen_SPIE11c, Chen_J, ChenJPS} and the references therein.

Throughout this paper, we shall use the following notations. The set of real numbers is denoted by $\bb{R}$. The set of integers is denoted by
$\bb{Z}$. The set of positive integers is denoted by $\bb{N}$. The notation ``$\Pr$'' denotes the probability measure.   The expectation of a
random variable is denoted by $\bb{E}[.]$.   We use ``$X$'' to denote the random variable of which its mean value is the subject of estimation.
Since we are only concerned with events which are generated by $X$ or its i.i.d. samples, the probability measure $\Pr$ is always associated
with the distribution of $X$.  Since all random variables involved in the sequel can be expressed in terms of $X$ or its i.i.d. samples,  the
expectation $\bb{E} [.]$ is always taken under the probability measure associated with the distribution of $X$.  Let $\Phi(.)$ denote the
standard normal distribution function, that is, $\Phi(x) = \f{1}{\sq{2 \pi}} \int_{-\iy}^x e^{- \f{x^2}{2} } dx$ for any $x \in \bb{R}$. Let
$\Phi^{-1}(.)$ denote the inverse function of $\Phi(.)$. The other notations will be made clear as we proceed.

\sect{Problem Formulation}

As mentioned in the introduction, many problems in sciences and engineering boils down to the estimation of the expectation $\mu$ of a random
variable $X$ based on its i.i.d. samples $X_1, X_2, \cd$.  Assume that the domain of the expectation $\mu = \bb{E} [X]$ is an interval denoted
by $\Se$. Let $\nu$ be the variance of $X$, i.e., $\nu = \bb{E} [ |X - \mu|^2]$.  Naturally,  the estimator $\wh{\bs{\mu}}$ for $\mu$ is taken
as the sample mean, where the sample number can be deterministic or random.  Since the estimator $\wh{\bs{\mu}}$ is of random nature, it is
important to provide a measure of its reliability. Let $\vep > 0, \; \vro \in (0, 1)$ and $\de \in (0, 1)$.  Based on different error criteria,
the estimation problems are typically posed in the following ways: \bsk

(I) Construct an estimator $\wh{\bs{\mu}}$  such that $\Pr \{ | \wh{\bs{\mu}} - \mu | < \vep  \} \geq 1 - \de$ for $\mu \in \Se$.

(II) Construct an estimator $\wh{\bs{\mu}}$ such that $\Pr \{ | \wh{\bs{\mu}} - \mu | < \vep |\mu| \} \geq 1 - \de$ for $\mu \in \Se$.

(III) Construct an estimator $\wh{\bs{\mu}}$  such that $\Pr \{ | \wh{\bs{\mu}} - \mu | <  \vro \vep  \; \tx{or} \;  | \wh{\bs{\mu}} - \mu | <
\vep |\mu|   \} \geq 1 - \de$ for $\mu \in \Se$.

(IV) Construct an estimator $\wh{\bs{\mu}}$  such that $\Pr \{ (1 - \vep) \wh{\bs{\mu}} < \mu < (1 + \vep) \wh{\bs{\mu}}  \} \geq 1 - \de$ for
$\mu \in \Se$.

\bsk

It can be seen that the above four problems can be accommodated as special cases of the following more general problem:

Construct an estimator $\wh{\bs{\mu}}$  such that $\Pr \{ \mscr{L} (\wh{\bs{\mu}}, \vep) < \mu < \mscr{U} (\wh{\bs{\mu}}, \vep) \} \geq 1 - \de$
for $\mu \in \Se$, where the functions $\mscr{L}(.,.)$ and $\mscr{U}(.,.)$ have the following properties: \bsk

(i) $\mscr{L}(z, \vep)$ and $\mscr{U}(z, \vep)$ are continuous functions of $z$ and $\vep > 0$.

(ii) For any $\vep > 0$,
\[ \mscr{L} (z, \vep) \leq z  \leq \mscr{U} (z, \vep) \qqu \fa z \in \bb{R},
\]
\[
\mscr{L} (\mu, \vep) < \mu < \mscr{U} (\mu, \vep) \qqu \fa \mu \in \Se.
\]

(iii) For any $\mu \in \Se$, \be \la{uniformasp} \lim_{ z \to \mu \atop{\vep  \downarrow 0 }  }  \f{ z - \mscr{L} (z, \vep) } {\vep} = \lim_{ z
\to \mu \atop{\vep  \downarrow 0 }  }  \f{ \mscr{U} (z, \vep) - z } {\vep} = \ka(\mu) > 0, \ee where $\ka(\mu)$ is a continuous function of
$\mu$.

\bsk

In view of (\ref{uniformasp}), the random interval $( \mscr{L} (\wh{\bs{\mu}}, \vep), \;  \mscr{U} (\wh{\bs{\mu}}, \vep) )$ is referred as to
{\it asymptotically symmetrical random interval}.  In the sequel, our objective is to construct such random intervals for $\mu$.  For
generality, we also consider random interval $( \mcal{L} (\wh{\bs{\mu}}, \vep), \;  \mcal{U} (\wh{\bs{\mu}}, \vep) )$, where the  functions
$\mcal{L}(., .)$ and $\mcal{U}(. , .)$  are ``in the same class''  as  $\mscr{L}(. , .)$ and $\mscr{U}(. , .)$  in the sense that for any $\mu
\in \Se$,
\[
\lim_{ z \to \mu \atop{\vep  \downarrow 0 }  }  \f{ z - \mcal{L} (z, \vep) } {\vep} = \lim_{ z \to \mu \atop{\vep  \downarrow 0 }  }  \f{
\mcal{U} (z, \vep) - z } {\vep} = \lim_{ z \to \mu \atop{\vep  \downarrow 0 }  }  \f{ z - \mscr{L} (z, \vep) } {\vep} = \lim_{ z \to \mu
\atop{\vep \downarrow 0 }  }  \f{ \mscr{U} (z, \vep) - z } {\vep}.
\]
For simplicity of notations, let $\mbf{I}$ and $\bs{\mcal{I}}$ denote the random intervals $( \mscr{L} (\wh{\bs{\mu}}, \vep), \;  \mscr{U}
(\wh{\bs{\mu}}, \vep) )$ and $( \mcal{L} (\wh{\bs{\mu}}, \vep), \;  \mcal{U} (\wh{\bs{\mu}}, \vep) )$, respectively.

Consider a fixed-size sampling scheme of sample size $n = \mcal{N} (\vep, \de, \mu, \nu)$, where \be \la{minimuss}
 \mcal{N} (\vep, \de, \mu,
\nu) = \f{\nu \mcal{Z}^2}{[ \ka (\mu) \vep ]^2} \qqu \tx{with $\mcal{Z} = \Phi^{-1} \li ( 1 - \f{\de}{2} \ri )$}. \ee If the estimator
$\wh{\bs{\mu}}$ for $\mu$ is taken as the sample mean $\f{ \sum_{i=1}^n X_i }{n}$, then it can be shown that \be \la{aspss889} \lim_{\vep
\downarrow 0} \Pr \{ \mu \in  \mbf{I} \} = \lim_{\vep \downarrow 0} \Pr \{  \mu \in  \bs{\mcal{I}} \} = 1 - \de \ee for all $\mu \in \Se$.
Hence, for small $\vep > 0$, the minimum sample size for a fixed-sample-size procedure can be approximated as $\mcal{N} (\vep, \de, \mu, \nu)$.
In the special case that the distribution of $X$ can be determined by its expectation $\mu$, then the minimum sample size can be approximated as
\[
\mcal{N} (\vep, \de, \mu)  = \f{\mcal{V} (\mu) \mcal{Z}^2}{[ \ka (\mu) \vep ]^2},  \qqu \tx{where} \qqu \mcal{V} (\mu)  = \bec \nu  & \tx{if $\mu  \in \Se$},\\
0 & \tx{otherwise} \eec \]

From (\ref{minimuss}), it can be seen that the minimum sample size satisfying (\ref{aspss889}) depends on the unknown mean $\mu$ and variance
$\nu$.  This is an inherent issue of a fixed-sample-size procedure.  Hence, it is desirable to develop sequential or multistage sampling schemes
such that (\ref{aspss889}) can be satisfied without the knowledge of $\mu$ and $\nu$. This goal will be accomplished in the remainder of the
paper.

\section{Fully Sequential Estimation}

In this section, we first propose general sequential estimation methods for $\mu = \bb{E} [ X ]$ by virtue of the  inclusion principle.
Afterward, we state their general properties and describe more concrete sampling schemes constructed with various confidence sequences.

\subsection{Sequential Sampling Schemes from Inclusion Principle} \la{sec31}

Let $X_1, X_2, \cd$ be i.i.d. sample of $X$.  Define
\[ Y_n = \f{ \sum_{i=1}^{n} X_i  }{n}, \qqu n \in \bb{N}.
\]
Let $\{ \al_n \}_{n  \in \bb{N}}$ be a sequence of numbers such that $\al_n \in (0, 1)$.  Let $m_\vep$ be a positive integer dependent on
$\vep$.  Let $(L_n, U_n), \; n \geq m_\vep$ be a sequence of confidence intervals such that for any $n \geq m_\vep$, $(L_n, U_n)$ is defined in
terms of $X_1, \cd, X_n$ and that the coverage probability $\Pr \{ L_n < \mu < U_n \}$ is no less than or close to $1 -\al_n$. By the inclusion
principle of \cite{Chenestimation}, we propose the following general stopping rule:

Continue sampling until the confidence interval $(L_n, U_n)$ is included by interval $(\mscr{L} (Y_n, \vep), \; \mscr{U} (Y_n, \vep) )$ for some
$n \geq m_\vep$.

Let  $\mbf{n}$ denote the sample size at the termination of the sampling process.  Clearly, $\mbf{n}$ is a random number.  The estimator for
$\mu$ is taken as $\wh{\bs{\mu}} = \f{ \sum_{i=1}^\mbf{n} X_i}{\mbf{n}} = Y_{\mbf{n}}$.  As before, let $\mbf{I}$ and $\bs{\mcal{I}}$ denote the
random intervals $( \mscr{L} (\wh{\bs{\mu}}, \vep), \; \mscr{U} (\wh{\bs{\mu}}, \vep ) )$ and  $( \mcal{L} (\wh{\bs{\mu}}, \vep), \; \mcal{U}
(\wh{\bs{\mu}}, \vep) )$, respectively.   It can be shown that \[
 \Pr \{ \mu \notin \mbf{I} \} \leq \sum_{ n = m_\vep  }^\iy  \min ( \al_n, \Pr \{ \mbf{n} = n \} )
\] provided that $\Pr \{ \mbf{n} < \iy \} = 1$ and that $\Pr \{  L_n < \mu < U_n \} \geq 1 - \al_n$ for $n \in \bb{N}$. For simplicity of the
stopping rule, we shall make effort to eliminate the computation of confidence limits.

\subsection{Characteristics of Sampling Schemes}

We expect that by appropriate choice of confidence sequence,  the sampling scheme proposed in Section \ref{sec31} have properties as follows.

\bed

\item [Property (I)]: \be \la{prop1} \Pr \{ \mbf{n} < \iy \} = 1, \qqu \bb{E} [ \mbf{n} ] < \iy.  \ee

\item [Property (II)]: \be \la{prop3} \Pr \li \{  \lim_{\vep \downarrow 0} \f{ \mbf{n} }{ \mcal{N} (\vep, \de, \mu, \nu) }  = 1 \ri \} = 1. \ee

\item [Property (III)]: \be \la{prop4} \lim_{\vep \downarrow 0} \Pr \{ \mu \in \bs{\mcal{I}} \} = 1 - \de.
 \ee

\item [Property (IV)]: \be \la{prop5} \lim_{\vep \downarrow 0} \f{ \bb{E} [ \mbf{n} ] }{ \mcal{N} (\vep, \de, \mu, \nu) }  = 1. \ee

\eed

\subsection{Using Confidence Sequences from Distribution Functions}

In this section, we consider the estimation of $\mu = \bb{E} [X]$ under the assumption that the distribution of $X$ is determined by $\mu$ such
that $\mcal{V} (\mu)$ is a positive continuous function of $\mu \in \Se$  and that $\bb{E} [ | X - \mu |^3 ]$ is a continuous function of $\mu
\in \Se$.  For $n \in \bb{N}$, define \[
 F_{Y_n} (z, \mu) = \bec \Pr \{ Y_n \leq z \} & \tx{if $\mu \in \Se$},\\
0 & \tx{otherwise}
 \eec  \qqu \qqu \qqu
 G_{Y_n} (z, \mu) = \bec \Pr \{ Y_n \geq z \} & \tx{if $\mu \in \Se$},\\
0 & \tx{otherwise} \eec \] Let $\{ \de_n \}_{n \in \bb{N}}$ be a sequence of positive numbers less than $1$ such that $\de_n \to \de$ as $n \to
\iy$. Let $m_\vep$ be a positive integer dependent on $\vep$ such that \[ \tx{$m_\vep \to \iy$ and $\vep^2 m_\vep \to 0$ as $\vep \downarrow
0$}. \] By virtue of the inclusion principle, we propose the following stopping rule:

Continue sampling until \be \la{ST893683}
 F_{Y_n} (Y_n, \mscr{U} (Y_n, \vep) ) \leq \f{\de_n}{2}, \qqu  G_{Y_n} (Y_n, \mscr{L} (Y_n, \vep)) \leq
\f{\de_n}{2} \ee for some integer $n$ no less than $m_\vep$.

For such stopping rule, we have established the following results.

\beT \la{ThmCDF}   The sampling scheme possesses the properties described by (\ref{prop1})--(\ref{prop5}). Moreover, if for arbitrary $n \in
\bb{N}, \; z \in \bb{R}$, the probability $\Pr \{ Y_n \leq z \}$ is non-increasing with respect to $\mu \in \Se$, then, \be \la{VIP886}
 \Pr \li
\{ \mu \notin \mbf{I} \ri \} \leq \sum_{n = m_\vep}^\iy \min \li ( \de_n, \Pr \{ \mbf{n} = n   \} \ri ). \ee

\eeT

See Appendix \ref{ThmCDF_app} for a proof.  It can be shown that the assumption that $\Pr \{ Y_n \leq z \}$ is non-increasing with respect to
$\mu \in \Se$ can be satisfied if either one of the following hold.

(i) The likelihood ratio
\[
\f{ f (X_1, \cd, X_n; \mu_2) }{ f (X_1, \cd, X_n; \mu_1)  }, \qqu \mu_1 < \mu_2, \qqu \mu_1, \; \mu_2 \in \Se
\]
is a non-decreasing function of $Y_n$, where $f$ denotes the probability density or mass functions.

(ii) $Y_n$ is a maximum likelihood estimator of $\mu$ in the sense that for any realization $(x_1, \cd, x_n)$, of $(X_1, \cd, X_n)$, the
likelihood function $f(x_1, \cd, x_n; \mu)$ is non-decreasing with respect to $\mu \leq \f{\sum_{i=1}^n x_i}{n}$ and non-increasing with respect
to $\mu \geq \f{\sum_{i=1}^n x_i}{n}$.

In order to derive the above stopping rule from the inclusion principle, we have used confidence intervals $(L_n, U_n), \; n \in \bb{N}$ with
$F_{Y_n} (Y_n, U_n) = \f{\de_n}{2}$ and $G_{Y_n} (Y_n, L_n ) = \f{\de_n}{2}$ for $n \in \bb{N}$.  Under the assumption that for arbitrary $n \in
\bb{N}, \; z \in \bb{R}$, the probability $\Pr \{ Y_n \leq z \}$ is non-increasing with respect to $\mu \in \Se$, it can be shown that $\Pr \{
L_n < \mu < U_n \} \geq 1 - \de_n$ for $n \in \bb{N}$ and that the stopping rule based on the inclusion principle is equivalent to the above
stopping rule defined by (\ref{ST893683}).

\subsection{Using Confidence Sequences from Large Deviation Theory} \la{sec3489}

In this section, we consider the estimation of $\mu = \bb{E} [X]$ in the context that the distribution of $X$ is determined by $\mu$ such that
the following assumptions are satisfied.

\bed

\item (i) The cumulant function \be \psi (s, \mu) = \ln \bb{E} [ \exp ( s (X - \mu ) ) ] \ee is well defined on $\{ (s, \mu): a(\mu) < s < b(\mu),
\; \mu \in \Se \}$, where $a(\mu)$ is $-\iy$ or a negative continuous function of $\mu \in \Se$; and similarly, $b(\mu)$ is $\iy$ or a positive
continuous function of $\mu \in \Se$.

\item (ii) $\f{\pa^3 \psi (s, \mu) }{\pa s^3}$ is a continuous function of $\mu \in \Se$ and $s \in \mscr{S}_\mu$, where
\[
\mscr{S}_\mu = \{s: a(\mu) < s < b(\mu) \}.
\]

\item (iii) $\bb{E} [ | X |^3 ] < \iy$  for $\mu \in \Se$.

\eed

To describe our sampling scheme, define \be \la{defM896}
\mscr{M} (z, \mu) = \bec \inf_{s \in \mscr{S}_\mu } [ s ( \mu - z )  + \psi (s, \mu) ]  & \tx{if $\mu \in \Se$},\\
- \iy & \tx{otherwise}. \eec \ee Let $\{ \de_n \}_{n \in \bb{N}} $ be a sequence of positive numbers less than $1$ such that $\Phi \li (  \sq{ 2
\ln \f{2}{\de_n} } \ri ) \to 1 - \f{\de}{2}$ as $n \to \iy$.  Let $m_\vep$ be a positive integer dependent on $\vep$ such that $m_\vep \to \iy$
and $\vep^2 m_\vep \to 0$ as $\vep \downarrow 0$.  By virtue of the inclusion principle, we propose the following stopping rule:

Continue sampling until \be \la{ST893683CH}  \mscr{M} (Y_n, \mscr{U} ( Y_n, \vep) ) \leq \f{ 1 }{n} \ln \f{\de_n}{2}, \qqu   \mscr{M} (Y_n,
\mscr{L} ( Y_n, \vep ) ) \leq \f{ 1 }{n} \ln \f{\de_n}{2} \ee for some integer $n$ no less than $m_\vep$.

With regard to such stopping rule, we have obtained the following results.

\beT \la{ThmCHB}   The sampling scheme possesses the properties described by
 (\ref{prop1})--(\ref{prop5}).  Moreover, if for arbitrary $n \in
\bb{N}, \; z \in \bb{R}$, the probability $\Pr \{ Y_n \leq z \}$ is non-increasing with respect to $\mu \in \Se$, then \be \la{VIP886b} \Pr \li
\{ \mu \notin \mbf{I} \ri \} \leq \sum_{n = m_\vep}^\iy \min \li ( \de_n, \Pr \{ \mbf{n} = n   \} \ri ). \ee

\eeT

See Appendix \ref{ThmCHB_app} for a proof.

For purpose of deriving the above stopping rule from the inclusion principle, we can define confidence interval $(L_n, U_n)$ with $\mscr{M}
(Y_n, U_n) = \f{ 1 }{n} \ln \f{\de_n}{2}, \; \mscr{M} (Y_n, L_n ) = \f{ 1 }{n} \ln \f{\de_n}{2}$ and $L_n \leq Y_n \leq U_n$ for $n \in \bb{N}$.
Under the assumption that for arbitrary $x \in \bb{R}$, the probability $\Pr \{ X \leq x \}$ is non-increasing with respect to $\mu \in \Se$, it
can be shown by large deviation theory that $\Pr \{ L_n < \mu < U_n \} \geq 1 - \de_n$ for $n \in \bb{N}$ and that the stopping rule based on
the inclusion principle is equivalent to the above stopping rule defined by (\ref{ST893683CH}).

\subsection{Using Confidence Sequences from Normal Approximation with Linear Interpolation}

In this section, we consider the estimation of $\mu = \bb{E} [X]$ under the assumption that the distribution of $X$ is determined by $\mu$ such
that $\mcal{V} (\mu)$ is a positive continuous function of $\mu \in \Se$  and that $\bb{E} [ | X|^3 ] < \iy$  for $\mu \in \Se$.

Let $\{ \de_n \}_{n \in \bb{N}}$ be a sequence of positive numbers less than $1$ such that {\small $\Phi \li (  \sq{ \ln \f{1}{\de_n}  } \ri )
\to 1 - \f{\de}{2}$} as $n \to \iy$. Let $\ro \in [0, 1]$.  Let $m_\vep$ be a positive integer dependent on $\vep$ such that $m_\vep \to \iy$
and $\vep^2 m_\vep \to 0$ as $\vep \downarrow 0$. By virtue of the inclusion principle, we propose the following stopping rule:

Continue sampling until {\small \bee  &  & \mcal{V} \li ( Y_n + \ro [\mscr{U} (Y_n, \vep) - Y_n ] \ri ) \leq \f{ n [ \mscr{U} (Y_n, \vep) - Y_n
]^2 }{ \ln \f{1}{\de_n} },  \qqu  \mcal{V} \li ( Y_n + \ro [\mscr{L} (Y_n, \vep) - Y_n ] \ri ) \leq \f{ n [  \mscr{L} (Y_n, \vep)  - Y_n ]^2 }{
\ln \f{1}{\de_n} }  \eee} for some integer $n$ no less than $m_\vep$.

For such stopping rule, we have proved the following result.

\beT \la{ThmNPL}  The sampling scheme possesses the properties described by (\ref{prop1})--(\ref{prop5}).

\eeT

See Appendix \ref{ThmNPL_app} for a proof.

In the course of deriving the above stopping rule from the inclusion principle, we have used the confidence sequences constructed as follows.
Let $\{ \al_n \}_{n \in \bb{N}}$ be a sequence of numbers such that $\li [ \Phi^{-1} \li ( 1 - \f{\al_n}{2} \ri ) \ri ]^2 = \ln \f{1}{\de_n}$.
For $n \in \bb{N}$, we can  construct a confidence interval $(L_n, U_n)$ such that $\Pr \{ L_n < \mu < U_n \} \ap 1 - \al_n$,  where the lower
confidence  limit $L_n$ and upper confidence limit $U_n$  satisfy $L_n \leq Y_n \leq U_n$ and
\[
\Phi \li ( \f{ \sq{{n}} [  L_n  - Y_n ]
 }{  \sq{ \mcal{V} \li ( Y_n + \ro [L_n - Y_n ] \ri ) } } \ri ) = \f{\al_n}{2}, \qqu
\Phi \li ( \f{ \sq{{n}} [  Y_n  - U_n ]
 }{  \sq{ \mcal{V} \li ( Y_n + \ro [U_n - Y_n ] \ri ) } } \ri ) = \f{\al_n}{2}.
 \]
These confidence intervals are constructed based on normal approximation with a linear approximation, which had been proposed in
\cite{Chenestimation}.

\subsection{Using Distribution-Free Confidence Sequences}

In this section, we consider the estimation of $\mu = \bb{E} [X]$ under the assumption that $\bb{E} [ | X |^6 ] < \iy$ for $\mu \in \Se$. Note
that the distribution of $X$ is not necessarily determined by $\mu$.  To describe our sampling scheme, define
\[
V_n = \f{ \sum_{i=1}^n (X_i - Y_n)^2 }{n}
\]
for $n \in \bb{N}$.  Let $\{ \de_n \}_{n \in \bb{N}}$ be a sequence of positive numbers less than $1$ such that {\small $\Phi \li (  \sq{ \ln
\f{1}{\de_n} } \ri ) \to 1 - \f{\de}{2}$} as $n \to \iy$.  Let $\ro \in (0, 1]$.  According to the inclusion principle, we propose the following
stopping rule:

 Continue sampling until \[ V_{n} + \f{\ro}{n} \leq  \f{ n}{  \ln \f{1}{\de_n} }   \times \min \li \{ | \mscr{U}( Y_n, \vep  ) - Y_n |^2, \; |
\mscr{L}(  Y_n, \vep  ) - Y_n |^2  \ri \}  \]
 for some $n \in \bb{N}$.

With respect to such stopping rule, we have shown the following result.

\beT \la{Seqgenfull}  The sampling scheme possesses the properties described by (\ref{prop1})--(\ref{prop5}). \eeT

See Appendix \ref{Seqgenfull_app} for a proof.

To derive the above stopping rule from the inclusion principle, we can use confidence interval $(L_n, U_n)$ with \bee L_n = Y_n -  \Phi^{-1} \li
( 1 - \f{\al_n}{2} \ri )  \sq{ \f{V_n}{n}  + \f{\ro}{n^2}  }, \qqu U_n = Y_n + \Phi^{-1} \li ( 1 - \f{\al_n}{2} \ri ) \sq{ \f{V_n}{n} +
\f{\ro}{n^2} } \eee for $n \in \bb{N}$, where $\{ \al_n \}_{n \in \bb{N}}$ is a sequence of positive numbers less than $1$ such that $\li [
\Phi^{-1} \li ( 1 - \f{\al_n}{2} \ri ) \ri ]^2 = \ln \f{1}{\de_n}$.  By the central limit theorem, it can be shown that $\Pr \{ L_n < \mu < U_n
\} \ap 1 - \al_n$ for large $n \in \bb{N}$.  The intervals $(L_n, U_n), \; n \in \bb{N}$ are called distribution-free confidence intervals
because their construction requires no information of the underlying distribution of $X$.

 \sect{Multistage Stopping Rules}

In this section, we first propose general multistage estimation methods for $\mu = \bb{E} [ X ]$ based on the  inclusion principle. Afterward,
we state their general properties and describe more concrete sampling schemes constructed with various confidence sequences.

\subsection{Multistage Sampling Schemes from Inclusion Principle} \la{sec42}

Let $\{ C_\ell \}_{\ell \in \bb{Z}}$ be an increasing  sequence of positive numbers such that \[ 1 < \inf_{\ell \in \bb{N}} \f{C_{\ell +
1}}{C_{\ell}} \leq \sup_{\ell \in \bb{N}} \f{C_{\ell + 1}}{C_{\ell}} < \iy. \]  Let $\{ \Up_\ell \}_{\ell \in \bb{Z}}$ be a non-decreasing
sequence of positive numbers such that
\[ 1 \leq \inf_{\ell \in \bb{N}} \f{\Up_{\ell + 1}}{\Up_{\ell}} \leq \sup_{\ell \in \bb{N}} \f{\Up_{\ell + 1}}{\Up_{\ell}} < \iy.
\]
Let $\tau_\vep$ be an integer-valued function of $\vep$ such that $\lim_{\vep \downarrow 0} \tau_\vep = 1$.  Such integer $\tau_\vep$ is
referred to as the {\it starting index}.   Choose sample sizes as an increasing sequence of positive integers $\{ n_\ell \}_{\ell \geq
\tau_\vep}$ dependent on $\vep$  such that \be \la{gensamplesize}
 \f{ \vep^2 n_\ell }{ C_\ell \Up_\ell } \to 1 \qu \tx{uniformly
for $\ell \geq 1$ as $\vep \downarrow 0$. } \ee  Define \[ Z_\ell = \f{ \sum_{i=1}^{n_\ell} X_i  }{n_\ell}, \qqu \ell \geq \tau_\vep,
\]
where $X_1, X_2, \cd$ are i.i.d. samples of $X$.   Let $(L_\ell, U_\ell), \; \ell \geq \tau_\vep$ be a sequence of confidence intervals such
that for any $\ell \geq \tau_\vep$, $(L_\ell, U_\ell)$ is defined in terms of $X_i, \; i = 1, \cd, n_\ell$  and that the coverage probability
$\Pr \{ L_\ell < \mu < U_\ell  \}$ is no less than or close to $1 -\de_\ell$. By the inclusion principle of \cite{Chenestimation}, we propose
the following general stopping rule:

Continue sampling until the confidence interval $(L_\ell, U_\ell)$ is included by interval $( \mscr{L} (Z_\ell, \vep), \; \mscr{U} (Z_\ell,
\vep) )$ for some $\ell \geq \tau_\vep$.

For $\ell \geq \tau_\vep$, let $\bs{D}_\ell$ be a decision variable such that $\bs{D}_\ell = 1$ if the confidence interval $( L_\ell, U_\ell )$
is included by interval $( \mscr{L} (Z_\ell, \vep), \; \mscr{U} (Z_\ell, \vep) )$ and $\bs{D}_\ell = 0$ otherwise.  Define $\bs{D}_\ell = 0$ for
$\ell < \tau_\vep$. Then, the stopping rule can be restated  as follows:
 \be
 \la{STgeneral}
\tx{Continue sampling until $\bs{D}_\ell = 1$ for some $\ell \geq \tau_\vep$}. \ee  The index $\bs{l}$ at the termination of sampling process is
referred to as the {\it stopping index}.  Clearly, $\bs{l}$ is a random number.  The sample number $\mbf{n}$ at the termination of sampling
process is $n_{\bs{l}}$.  The estimator for $\mu$ is taken as $\wh{\bs{\mu}} = \f{ \sum_{i=1}^\mbf{n} X_i}{\mbf{n}} = Z_{\bs{l}}$.  As before,
let $\mbf{I}$ and $\bs{\mcal{I}}$ denote the random intervals $( \mscr{L} (\wh{\bs{\mu}}, \vep), \; \mscr{U} (\wh{\bs{\mu}}, \vep ) )$ and  $(
\mcal{L} (\wh{\bs{\mu}}, \vep), \; \mcal{U} (\wh{\bs{\mu}}, \vep) )$, respectively.  It can be shown that \[
 \Pr \{ \mu \notin \mbf{I} \} \leq \sum_{ \ell = \tau_\vep  }^\iy  \min ( \de_\ell, \Pr \{ \bs{l} = \ell \} )
\] provided that $\Pr \{ \bs{l} < \iy \} = 1$ and that $\Pr \{  L_\ell < \mu < U_\ell \} \geq 1 - \de_\ell$ for $\ell \geq
\tau_\vep$.

For simplicity of the stopping rule, we propose to eliminate the computation of confidence limits.

\subsection{General Characteristics  of Sampling Schemes}

To describe the characteristic of the sampling schemes proposed in Section \ref{sec42}, we need to introduce some quantities.  Define \[
\xi_\ell = \ka (\mu) \sq{ \f{ C_\ell \Up_\ell}{ \nu } }, \qqu  \Lm_\ell (\mu, \nu) =  \f{ \ka^2 (\mu) C_\ell }{\nu} \] for $\ell \in \bb{Z}$.
Specially, if $\nu$ is a function of $\mu$, that is $\nu = \mcal{V} (\mu)$, we write \[ \Lm_\ell (\mu) =  \f{ \ka^2 (\mu) C_\ell }{\mcal{V}
(\mu)} \] for $\ell \in \bb{Z}$.  For $\ell \geq \tau_\vep$, let $\bs{\mcal{I}}_\ell$ denote the random interval $( \mcal{L} (Z_\ell, \vep), \;
\mcal{U} (Z_\ell, \vep ) )$.  Define {\small \bee & & \ovl{P} = \sum_{\ell = \tau_\vep}^\iy \Pr \{ \mu \notin \bs{\mcal{I}}_\ell, \;
\bs{D}_{\ell - 1} = 0, \;
\bs{D}_\ell = 1 \}, \\
&  &  \udl{P} = 1 - \sum_{\ell = \tau_\vep}^\iy \Pr \{ \mu \in \bs{\mcal{I}}_\ell, \; \bs{D}_{\ell - 1} = 0, \; \bs{D}_\ell = 1 \}, \\
&  & Q = \sum_{\ell = \tau_\vep}^\iy \min \li [ \Pr \{ \mu \notin  \bs{\mcal{I}}_\ell \}, \; \Pr \{ \bs{D}_{\ell - 1} = 0 \}, \; \Pr \{
\bs{D}_\ell = 1 \}  \ri ] \eee} and \be \la{defim89} \jm = \min \{ \ell \geq 1: \Lm_\ell (\mu, \nu) > 1 \}. \ee Actually, $\jm$ is a function of
$\mu$ and $\nu$. We expect that by appropriate choice of sample sizes and confidence sequences, the sampling scheme proposed in Section
\ref{sec42} have the following general properties.

\bed

\item [Property (I)]:  \be \la{perty891}
 \Pr \{ \mbf{n} < \iy \} = 1, \qqu \bb{E} [ \mbf{n} ] < \iy. \ee

\item [Property (II)]: \be
 \la{propertyii1}
 \Pr
\li \{ \jm - 1 \leq \liminf_{\vep \downarrow 0} \bs{l} \leq \limsup_{\vep \downarrow 0} \bs{l}  \leq \jm \ri \} = 1. \ee Moreover,  \be
\la{propertyii2} \Pr \li \{ \lim_{\vep \downarrow 0} \bs{l} = \jm \ri \} = 1 \qu \tx{provided that $\jm = 1$ or both $\jm
> 1$ and $\Lm_{\jm-1} (\mu, \nu) < 1$ hold}. \ee

\item [Property (III)]: \be \la{propertyiii1} \liminf_{\vep \downarrow 0 } \Pr \{  \mu \in \bs{\mcal{I}} \} \geq 2 [ \Phi \li ( \xi_{\jm -1} \ri )  + \Phi \li ( \xi_\jm
\ri ) ] - 3.  \ee Moreover, \be \la{propertyiii2} \lim_{\vep \downarrow 0 } \Pr \{  \mu \in \bs{\mcal{I}} \} = 2 \Phi \li ( \xi_\jm \ri ) - 1
\qu \tx{provided that $\jm = 1$ or both $\jm
> 1$ and $\Lm_{\jm-1} (\mu, \nu) < 1$ hold}. \ee

\item [Property (IV)]: \bel  \li ( \f{ \xi_{\jm-1}  }{ \mcal{Z}} \ri )^2 \leq \liminf_{\vep \downarrow 0} \f{ \bb{E} [ \mbf{n} ] } { \mcal{N} (\vep, \de, \mu, \nu) }
\leq \limsup_{\vep \downarrow 0} \f{ \bb{E} [ \mbf{n} ] } { \mcal{N} (\vep, \de, \mu, \nu)  }  \leq \li ( \f{ \xi_\jm  }{ \mcal{Z}} \ri )^2.
\la{propertyiv1} \eel Moreover, \be \lim_{\vep \downarrow 0 } \f{ \bb{E} [ \mbf{n} ] } { \mcal{N} (\vep, \de, \mu, \nu) }  = \li ( \f{ \xi_\jm
}{ \mcal{Z}} \ri )^2 \qu \tx{provided that $\jm = 1$ or both $\jm
> 1$ and $\Lm_{\jm-1} (\mu, \nu) < 1$ hold}. \la{propertyiv3} \ee

\item [Property (V)]:  For all $\ell \in
\bb{Z}$, \bel  &  & \Pr \{ \bs{l} > \ell \} \leq \Pr \{ \bs{D}_\ell = 0 \}, \qqu \qqu  \lim_{\vep \downarrow 0} [  \Pr \{ \bs{D}_\ell = 0 \} -
\Pr \{ \bs{l} > \ell \}   ] = 0.  \la{propertyv2} \eel

\item [Property (VI)]: \bel  &  & \udl{P} \leq \Pr \{ \mu \notin \bs{\mcal{I}} \} \leq \ovl{P} \leq Q, \la{propertyvi1} \\
&  & \lim_{\vep \downarrow 0} | \Pr \{ \mu \notin \bs{\mcal{I}} \} - \ovl{P} | = \lim_{\vep \downarrow 0}  | \Pr \{ \mu \notin \bs{\mcal{I}} \}
- \udl{P} | = 0. \la{propertyvi2}
 \eel
Moreover, \be \la{propertyvi3} \lim_{\vep \downarrow 0} [ Q -  \Pr \{ \mu \notin \bs{\mcal{I}} \} ] = 0 \qu \tx{provided that $\jm = 1$ or both
$\jm > 1$ and $\Lm_{\jm-1} (\mu, \nu) < 1$ hold}. \ee

\item [Property (VII)]:  \bel
&  & \bb{E} [ \mbf{n} ] \leq n_{\tau_\vep} + \sum_{\ell = \tau_\vep}^\iy (n_{\ell + 1} - n_\ell)  \Pr \{ \bs{D}_\ell = 0  \}, \la{propertyvii1}\\
&  &  \lim_{\vep \downarrow 0} \f{ 1 }{\bb{E} [ \mbf{n} ] } \li [ n_{\tau_\vep} + \sum_{\ell = \tau_\vep}^\iy (n_{\ell + 1} - n_\ell)  \Pr \{
\bs{D}_\ell = 0 \} \ri ] = 1. \la{propertyvii2} \eel

\eed

\subsection{Using Confidence Sequences from Distribution Functions}

In this section, we consider the multistage estimation of $\mu = \bb{E} [X]$ under the assumption that the distribution of $X$ is determined by
$\mu$ such that $\mcal{V} (\mu)$ is a positive continuous function of $\mu \in \Se$  and that $\bb{E} [ | X |^3 ]$ is a continuous function of
$\mu \in \Se$.  For $\ell \geq \tau_\vep$, define \[
 F_{Z_\ell} (z, \mu) = \bec \Pr \{ Z_\ell \leq z  \} & \tx{if $\mu \in \Se$},\\
0 & \tx{otherwise}
 \eec  \qqu \qqu \qqu \qqu
 G_{Z_\ell} (z, \mu) = \bec \Pr \{ Z_\ell \geq z  \} & \tx{if $\mu \in \Se$},\\
0 & \tx{otherwise} \eec \] Let $\{ \de_\ell \}_{\ell \in \bb{Z}}$ be a non-increasing sequence of positive numbers less than $1$. Assume that
the sample sizes at all stages are given by (\ref{gensamplesize}) with
\[
\Up_\ell =  \li [ \Phi^{-1} \li ( 1 - \f{\de_\ell}{2} \ri ) \ri ]^2, \qqu \ell \in \bb{Z}
\]
such that $\de_\ell ( C_\ell \Up_\ell)^2  \to \iy$ as $\ell \to \iy$.    According to the inclusion principle, we propose the following stopping
rule:

Continue sampling until \be  F_{Z_\ell} \li ( Z_\ell, \mscr{U} (Z_\ell, \vep) \ri ) \leq \f{\de_\ell}{2}, \qqu   G_{Z_\ell} \li ( Z_\ell,
\mscr{L} (Z_\ell, \vep) \ri ) \leq \f{\de_\ell}{2} \la{stCDF1}\ee for some $\ell \geq \tau_\vep$.

To express this stopping rule in the general form (\ref{STgeneral}), the decision variables $\bs{D}_\ell, \; \ell \geq \tau_\vep$ can be defined
such that $\bs{D}_\ell = 1$ if (\ref{stCDF1}) is  true and $\bs{D}_\ell = 0$ otherwise.  With regard to such sampling scheme, we have
established the following results.

\beT \la{CHCDF}  The sampling scheme possesses properties (I) -- (VII) described by (\ref{perty891}) -- (\ref{propertyvii2}).  Moreover, if for
arbitrary $\ell \geq \tau_\vep$ and $z \in \bb{R}$, the probability $\Pr \{ Z_\ell \leq z \}$ is non-increasing with respect to $\mu \in \Se$,
then \be \la{VIP886c} \Pr \li \{ \mu \notin \mbf{I} \ri \} \leq \sum_{\ell = \tau_\vep}^\iy \min \li ( \de_\ell, \Pr \{ \bs{l} = \ell   \} \ri
). \ee \eeT

See Appendix \ref{CHCDF_app} for a proof.

\subsection{Using Confidence Sequences from Large Deviation Theory}

In this section, we consider the estimation of $\mu = \bb{E} [X]$ in the context that the distribution of $X$ is determined by $\mu$ such that
the assumptions (i), (ii), (iii) given in Section \ref{sec3489} are satisfied.

Let $\{ \de_\ell \}_{\ell \in \bb{Z}}$ be a non-increasing sequence of positive numbers less than $1$.  Assume that the sample sizes at all
stages are given by (\ref{gensamplesize}) with
\[
\Up_\ell = 2 \ln \f{2}{\de_\ell}, \qqu \ell \in \bb{Z}
\]
and $C_\ell \to \iy$ as $\ell \to \iy$.   Recall the definition of the bivariate function $\mscr{M}(.,.)$ given by (\ref{defM896}).  By virtue
of the inclusion principle, we propose the following stopping rule:

Continue sampling until \be  \mscr{M} (Z_\ell, \mscr{U} ( Z_\ell, \vep) ) \leq \f{ 1 }{{n_\ell}} \ln \f{\de_\ell}{2}, \qqu \qqu   \mscr{M}
(Z_\ell, \mscr{L} ( Z_\ell, \vep ) ) \leq \f{ 1 }{{n_\ell}} \ln \f{\de_\ell}{2}  \la{stCH1} \ee for some $\ell \in \bb{Z}$ no less than
$\tau_\vep$.

To express this stopping rule in the general form (\ref{STgeneral}), the decision variables $\bs{D}_\ell, \; \ell \geq \tau_\vep$ can be defined
such that $\bs{D}_\ell = 1$ if (\ref{stCH1}) is true and $\bs{D}_\ell = 0$ otherwise.

For such stopping rule, we have the following results.

\beT \la{CHST}  The sampling scheme possesses properties (I) -- (VII) described by (\ref{perty891}) -- (\ref{propertyvii2}).  Moreover, if for
arbitrary $\ell \geq \tau_\vep$ and $z \in \bb{R}$, the probability $\Pr \{ Z_\ell \leq z \}$ is non-increasing with respect to $\mu \in \Se$,
then \be \la{VIP886d} \Pr \li \{ \mu \notin \mbf{I} \ri \} \leq \sum_{\ell = \tau_\vep}^\iy \min \li ( \de_\ell, \Pr \{ \bs{l} = \ell   \} \ri
). \ee \eeT

See Appendix \ref{CHST_app} for a proof.

\subsection{Using Confidence Sequences from Normal Approximation with Linear Interpolation}

In this section, we consider the multistage  estimation of $\mu = \bb{E} [X]$ under the assumption that the distribution of $X$ is determined by
$\mu$ such that $\mcal{V} (\mu)$ is a positive continuous function of $\mu \in \Se$  and that $\bb{E} [ | X |^3 ] < \iy$ for $\mu \in \Se$.

Let $\{ \de_\ell \}_{\ell \in \bb{Z}}$ be a non-increasing sequence of positive numbers less than $1$ such that $\de_\ell \in (0, 1)$ for all
$\ell \in \bb{Z}$. Assume that the sample sizes at all stages are given by (\ref{gensamplesize}) with $\Up_\ell = \ln \f{1}{\de_\ell},  \; \;
\ell \in \bb{Z}$ and $C_\ell \to \iy$ as $\ell \to \iy$.  Let $\ro \in [0, 1]$.   According the inclusion principle, we propose the following
stopping rule:

Continue sampling until {\small \be  \mcal{V} \li ( Z_\ell + \ro [\mscr{U} (Z_\ell, \vep) - Z_\ell ] \ri )  \leq \f{ n_\ell [  \mscr{U} (Z_\ell,
\vep) - Z_\ell ]^2 }{ \ln \f{1}{\de_\ell} }, \qqu  \mcal{V} \li ( Z_\ell + \ro [\mscr{L} (Z_\ell, \vep) - Z_\ell ] \ri ) \leq \f{ n_\ell [
\mscr{L} (Z_\ell, \vep)  -  Z_\ell ]^2}{\ln \f{1}{\de_\ell}} \la{stnal1}  \ee} for some $\ell \geq \tau_\vep$.

To express this stopping rule in the general form (\ref{STgeneral}), the decision variables $\bs{D}_\ell, \; \ell \geq \tau_\vep$ can be defined
such that $\bs{D}_\ell = 1$ if (\ref{stnal1}) is  true and $\bs{D}_\ell = 0$ otherwise.

For such stopping rule, we have the following result.

\beT \la{NLST}

The sampling scheme possesses properties (I) -- (VII) described by (\ref{perty891}) -- (\ref{propertyvii2}).

\eeT

See Appendix \ref{NLSTapp} for a proof.

\subsection{Using Distribution-Free Confidence Sequences}

In this section, we consider the multistage estimation of $\mu = \bb{E} [X]$ under the assumption that $\bb{E} [ | X  |^6 ] < \iy$ for $\mu \in
\Se$.  Note that the distribution of $X$ is not necessarily determined by $\mu$.  To describe our sampling scheme, define
\[
W_\ell = \f{ \sum_{i=1}^{n_\ell} (X_i - Z_\ell)^2 }{n_\ell}
\]
for $n \in \bb{N}$. Let $\{ \de_\ell \}_{\ell \in \bb{Z}}$ be a non-increasing sequence of positive numbers less than $1$.   Let $\ro \geq 0$.
Assume that the sample sizes at all stages are given by (\ref{gensamplesize}) with
\[
\Up_\ell =  \ln \f{1}{\de_\ell},  \qqu \ell \in \bb{Z}
\]
and $C_\ell \to \iy$ as $\ell \to \iy$.  By virtue of the inclusion principle, we propose the following stopping rule:

Continue sampling until \be W_\ell + \f{\ro}{n_\ell} \leq \f{ n_\ell}{ \ln \f{1}{\de_\ell} } \times \min \li \{ [ \mscr{U}( Z_\ell, \vep ) -
Z_\ell ]^2, \;  [ \mscr{L}( Z_\ell, \vep  ) - Z_\ell ]^2 \ri \}  \la{stfree} \ee
 for some $\ell \geq \tau_\vep$ in $\bb{Z}$.

To express this stopping rule in the general form (\ref{STgeneral}), the decision variables $\bs{D}_\ell, \; \ell \geq \tau_\vep$ can be defined
such that $\bs{D}_\ell = 1$ if (\ref{stfree}) is true and $\bs{D}_\ell = 0$ otherwise.

With regard to such stopping rule, we have established the following result.

\beT \la{Seqgen}  The sampling scheme possesses properties (I) -- (VII) described by (\ref{perty891}) -- (\ref{propertyvii2}).
 \eeT

See Appendix \ref{Seqgen_app} for a proof.

\section{Concluding Remarks}

In this paper, we have proposed to apply the inclusion principle to develop sequential schemes for estimating mean values.  The central idea of
the inclusion principle is to use a sequence of confidence intervals to construct stopping rules so that the sampling process is continued until
a confidence interval is included by an interval defined in terms of the samples.  Our asymptotic results demonstrate that the stopping rules
derived from the inclusion principle usually possess the best possible efficiency.

We would like to point out that for all sampling schemes proposed in this paper, the prescribed confidence level $1 - \de$ is asymptotically
satisfied as $\vep \to 0$.  However, $\vep \to 0$ implies that the average sample number tends to infinity.  To overcome this issue,  one can
apply the exact computational techniques developed in \cite{Chen_SPIE, Chenestimation} to adjust the confidence parameter $\de$ so that the
pre-specified confidence level for a given $\vep$ can be guaranteed with little conservatism.

\appendix

\sect{Concentration Inequalities}

To study the sampling schemes proposed in this paper, we need some concentration inequalities.

\beL

\la{nonuniformBE}

Define $\mscr{W} = \bb{E} [ | X - \mu |^3 ]$. The following statements hold true.
 \bel
&  & \li | \Pr \{ Y_n \leq z  \} - \Phi \li ( \f{ \sq{n} ( z - \mu  )   }
{  \sq{ \nu }  } \ri ) \ri | < \f{ \mscr{W}}{ \sq{n \nu^3}} \qu \tx{for any $z \in \bb{R}$}, \la{berry55}\\
&  & \li | \Pr \{ Y_n \geq z  \} - \Phi \li ( \f{ \sq{n} ( \mu - z  )   }
{  \sq{ \nu }  } \ri ) \ri | < \f{ \mscr{W}}{ \sq{n \nu^3}} \qu \tx{for any $z \in \bb{R}$}, \la{berry56} \\
&  & \Pr \{ Y_n \leq z  \}  < \f{1}{2} \exp \li ( - \f{n}{2} \f{ | z - \mu |^2 }{  \nu } \ri ) +  \f{ C }{n^2}
\f{ \mscr{W} }{ |z - \mu|^3 } \qu \tx{for $z$ less than $\mu \in \Se$}, \la{beryb1}\\
&  & \Pr \{ Y_n \geq z  \}  < \f{1}{2} \exp \li ( - \f{n}{2} \f{ | z - \mu |^2 }{  \nu } \ri ) +  \f{ C }{n^2} \f{ \mscr{W} }{
|z - \mu|^3 } \qu \tx{for $z$ greater than $\mu \in \Se$},  \la{beryb2} \\
&  & \Pr \{ | Y_n - \mu | \geq \ga   \}  < \exp \li ( - \f{n}{2} \f{ \ga^2 }{  \nu} \ri ) + \f{2 C }{n^2} \f{ \mscr{W} }{ \ga^3 } \qu \tx{for
any $\ga > 0$},  \la{mostvip} \eel where $C$ is an absolute constant.

\eeL

\bpf
 Note that   \[ \Pr \{ Y_n \leq z  \} = \Pr \li \{ \f{ \sq{n} ( Y_n
- \mu )   }{ \sq{ \nu} } \leq \f{ \sq{n} ( z - \mu  ) }{ \sq{ \nu }  }   \ri \}  \] for $\mu \in \Se$. Hence, (\ref{berry55}) follows from
Berry-Essen inequality.  In a similar manner,  making use of Berry-Essen inequality and the observation that
\[ \Pr \{ Y_n \geq z  \} = \Pr \li \{ \f{ \sq{n} ( \mu
- Y_n  )   }{ \sq{ \nu} } \leq \f{ \sq{n} ( \mu - z ) }{ \sq{ \nu }  }  \ri \}
\] for $\mu \in \Se$, we have that (\ref{berry56}) is true.

By the non-uniform version of Berry-Essen's inequality, \be \la{byv} \li | \Pr \{ Y_n \leq z  \} - \Phi \li ( \f{ \sq{n} ( z - \mu ) }{ \sq{
\nu} } \ri ) \ri | <  \f{ C \mscr{W}}{\sq{ n \nu^3  }  + n^2 | z - \mu |^3 } < \f{C \mscr{W} }{ n^2 |z - \mu|^3 } \ee for $z \in \bb{R}$ and
$\mu \in \Se$.  The inequality (\ref{beryb1}) immediately follows from (\ref{byv}) and the fact that $\Phi(x) < \f{1}{2} \exp (- \f{x^2}{2} ),
\; \fa x < 0$. In a similar manner, we can show inequality (\ref{beryb2}).   Finally, combining (\ref{beryb1}) and (\ref{beryb2}) yields
(\ref{mostvip}). This completes the proof of the lemma.

\epf

\beL

\la{vipnow}

Assume that $\bb{E} [ |X|^6 ] < \iy$.  Then, for $\eta > 0$,  there exists a constant $C > 0$ such that \bee &  & \Pr \{ | V_n - \nu | \geq \eta
\} \leq \exp \li ( - \f{n}{4} \f{ \eta}{ \nu } \ri ) +  \exp \li ( - \f{n}{8} \f{ \eta^2 }{ \varpi } \ri ) + \f{4 C }{n^2 \eta^3}  \li ( \sq{2}
\mscr{W} \eta^{3\sh 2} +  4 \mscr{V} \ri ) \eee where
\[ \nu = \bb{E} [ | X - \mu |^2 ] < \iy, \qqu \mscr{W} =
\bb{E} [ | X - \mu |^3 ] < \iy, \]
\[ \varpi = \bb{E} [ | ( X - \mu
)^2 - \nu |^2 ] < \iy, \qqu \mscr{V} = \bb{E} [ | ( X - \mu )^2 - \nu |^3 ] < \iy.
\]
\eeL

\bpf

By the assumption that $\bb{E} [ |X|^6 ] < \iy$, it must be true that
\[
\nu < \iy, \qqu \mscr{W} < \iy, \qqu \varpi < \iy, \qqu \mscr{V} < \iy.
\]
From (\ref{mostvip}) of Lemma \ref{nonuniformBE}, we have that there exists a constant $C > 0$ such that \bee &  & \Pr \li \{  | Y_n - \mu |
\geq \sq{ \f{\eta}{2} } \ri \} \leq \exp \li ( - \f{n}{4} \f{ \eta }{  \nu } \ri ) + \f{4
\sq{2} C }{n^2} \f{ \mscr{W} }{ \eta^{3\sh 2} },\\
&  & \Pr \li \{  | U_n - \nu | \geq \f{\eta}{2} \ri \} \leq \exp \li ( - \f{n}{8} \f{ \eta^2 }{ \varpi } \ri ) + \f{16 C }{n^2} \f{ \mscr{V} }{
\eta^3 }, \eee where \[ U_n = \f{ \sum_{i=1}^n (X_i - \mu)^2 }{n}. \] Since $U_n  = V_n + (Y_n - \mu)^2$, we have $\Pr \{  V_n \geq \nu + \eta
\} \leq \Pr \{  U_n \geq \nu + \eta  \}$. Moreover, \bee \Pr \{  V_n \leq \nu - \eta  \} & \leq & \Pr \li \{  U_n - (Y_n - \mu)^2 \leq \nu -
\eta, \;  (Y_n - \mu)^2 < \f{\eta}{2} \ri \}
+ \Pr \li \{ (Y_n - \mu)^2 \geq \f{\eta}{2}  \ri \}\\
& \leq & \Pr \li \{  U_n  \leq \nu - \f{\eta}{2} \ri \} + \Pr \li \{ | Y_n - \mu | \geq \sq{ \f{\eta}{2} } \ri \}. \eee  It follows that
\bee \Pr \{ | V_n - \nu | \geq \eta  \} & \leq  & \Pr \li \{ | Y_n - \mu | \geq \sq{ \f{\eta}{2} } \ri \} + \Pr \li \{ | U_n  - \nu | \geq \f{\eta}{2} \ri \}\\
& \leq  & \exp \li ( - \f{n}{4} \f{ \eta}{ \nu } \ri ) + \f{4 \sq{2} C }{n^2} \f{ \mscr{W} }{ \eta^{3\sh 2} }  + \exp \li ( - \f{n}{8} \f{
\eta^2 }{ \varpi } \ri ) + \f{16 C }{n^2} \f{
\mscr{V} }{ \eta^3 }\\
&  =  & \exp \li ( - \f{n}{4} \f{ \eta}{ \nu } \ri ) +  \exp \li ( - \f{n}{8} \f{ \eta^2 }{ \varpi } \ri ) + \f{4 C }{n^2 \eta^3}  \li (  \sq{2}
\mscr{W} \eta^{3\sh 2} +  4 \mscr{V} \ri ).  \eee This completes the proof of the lemma.

\epf

\section{Analysis of Function $\mscr{M}$}

This section is devoted to the analysis of the function $\mscr{M}$.  The results are useful for investigating sampling schemes associated with
the function $\mscr{M}$.  To simplify notations, define a multivariate function \be \la{defg3}
 g(s, \se, z) = s (\se - z) + \psi (s, \se)
\ee for $z \in \bb{R}, \; \se \in \Se$ and $a(\se) < s < b (\se)$.

\beL \la{lemgooday} Let $[u_1, u_2] \subseteq \Se$ and $d > 0$. Assume that  $\f{\pa \psi (s, \se) }{\pa s}$ is a continuous function of $s \in
(-d, d)$ and $\se \in (u_1, u_2)$. Define
\[
\mcal{M} (z, \se) = \inf_{s \in [-c, c]}  [s (\se - z) + \psi (s, \se) ]
\]
with $c \in (0, d)$.  Then, $\mcal{M} (z, \se)$ is a continuous function of $\se \in (u_1, u_2)$ and $z \in \bb{R}$. Moreover,
\[
\mcal{M} (z, \se) < 0 \qu \tx{if} \qu z \neq \se.
\]

\eeL

\bpf

Since the variance of $X$ is greater than $0$ for all $\se \in \Se$, the cumulant function $\psi (s, \se)$ is a convex function of $s$ for all
$\se \in \Se$. As a consequence of the assumption that $\f{\pa \psi (s, \se) }{\pa s}$ is a continuous function of $\se \in (u_1, u_2)$ and $s
\in (-d, d)$, we have that $g(s, \se, z)$ is a continuous function of $\se \in (u_1, u_2), \; z \in \bb{R}$ and $s \in [-c, c]$.  It follows
that, for given $\se \in (u_1, u_2), \; z \in \bb{R}$, the minimum of $g(s, \se, z)$ with respect to $s \in [-c, c]$ is attained at a number
$s_0 \in [-c, c]$, that is, $\mcal{M} (z, \se) = g(s_0, \se, z)$. Clearly,  $s_0 \in [-c, c]$ is a function of $z, \; \se$.  To establish the
continuity of $\mcal{M} (z, \se)$, it suffices to consider three cases as follows.

\bsk

Case (i): $\f{\pa g (s, \se, z) }{\pa s} = 0$ for $s = s_0 \in [-c, c]$.

Case (ii): $\f{\pa g (s, \se, z) }{\pa s} > 0$ for $s \in [-c, c]$.

Case (iii) $\f{\pa g (s, \se, z) }{\pa s} < 0$ for $s \in [-c, c]$.

\bsk

To show the continuity of $\mcal{M} (z, \se)$ in Case (i), note that for  $\ep > 0$,  there exists a positive number $\eta$ less than $\min \{
u_2 - \se, \se - u_1, d - s_0, s_0 + d \}$ such that \be \la{graet89}
 | g(s, \vse, y) - g(s_0, \se, z) | < \ep \qu \tx{for $|s - s_0| < \eta, \;
| \vse - \se | < \eta, \; |y - z| < \eta$}, \ee as a consequence of the continuity of  $g(s, \vse, y)$ with respect to $s, \vse$ and $y$. Note
that $\f{\pa g (s, \se, z) }{\pa s} > 0$ for $s = s_0 + \eta$ and that $\f{\pa g (s, \se, z) }{\pa s} < 0$ for $s = s_0 - \eta$ because $g(s,
\se, z)$ is convex on $s$. Since $\f{\pa g (s, \se, z) }{\pa s}$ is a continuous function of $s, \se$ and $z$, there exists $\vsi \in (0, \eta)$
such that
\[
\f{\pa g (s, \vse, y) }{\pa s} |_{ s_0 + \eta } > 0, \qqu \f{\pa g (s, \vse, y) }{\pa s} |_{ s_0 - \eta } < 0
\]
for $| \vse - \se | \leq \vsi$ and $|y - z| \leq \vsi$.  Hence, for any $(\vse, y)$ with $| \vse - \se | \leq \vsi, \; |y - z| \leq \vsi$, there
exists $s_1 \in (s_0 - \eta, s_0 + \eta)$ such that $\f{\pa g (s, \vse, y) }{\pa s} = 0$ for $s = s_1$.  Let $s_2$ denote the minimizer of $g(s,
\vse, y)$ with respect to $s \in [-c, c]$ for $| \vse - \se | \leq \vsi, \; |y - z| \leq \vsi$.  If $s_1 \in [-c, c]$, then $s_2 = s_1$ and thus
$|s_2 - s_0| \leq | s_2 - s_1 | + | s_1 - s_0| = | s_1 - s_0| < \eta$ because $g(s, \vse, y)$ is convex on $s \in [-c, c]$.  If $s_1 < -c$, then
$s_0 - \eta < s_1 < s_2 = - c \leq s_0$, which implies that $| s_2 - s_0| < \eta$. If $s_1
> c$, then $s_0 + \eta
> s_1 > s_2 =  c \geq s_0$, which implies that $| s_2 - s_0| < \eta$.
Therefore, we have $| s_2 - s_0| < \eta$ for all $(\vse, y)$ with $| \vse - \se | \leq \vsi, \; |y - z| \leq \vsi$.   It follows from
(\ref{graet89}) that
\[
| g(s_2, \vse, y) - g(s_0, \se, z) | < \ep \qu \tx{for $| \vse - \se | < \vsi, \; |y - z| < \vsi$}.
\]
Thus, we have $\vsi > 0$ such that $u_1 < \se - \vsi < \se + \vsi < u_2$ and that
\[
| \mcal{M} (y, \vse) - \mcal{M} (z, \se)  | = | g(s_2, \vse, y) - g(s_0, \se, z) | < \ep
\]
for $| \vse - \se | < \vsi, \; |y - z| < \vsi$.  This proves the continuity of $\mcal{M} (z, \se)$  on $\se \in (u_1, u_2)$ and $z \in \bb{R}$
for Case (i).

To show the continuity of $\mcal{M} (z, \se)$ in Case (ii), note that for  $\ep > 0$,  there exists a positive number $\eta$ less than $\min \{
u_2 - \se, \se - u_1\}$ such that
\[
| g(-c, \vse, y) - g(-c, \se, z) | < \ep \qu \tx{for $| \vse - \se | < \eta, \; |y - z| < \eta$},
\]
as a consequence of the continuity of $g(-c, \vse, y)$ with respect to $\vse$ and $y$.  Since $\f{\pa g (s, \se, z) }{\pa s} > 0$ for $s \in
[-c, c]$ and $\f{\pa g (s, \vse, y) }{\pa s}$ is a continuous function of $s, \vse$ and $y$, there exists $\vsi \in (0, \eta)$ such that
\[
\f{\pa g (s, \vse, y) }{\pa s} > 0 \qu \tx{for $| \vse - \se | \leq \vsi, \; |y - z| \leq \vsi$ and $s \in [-c, c]$}.
\]
This implies that the minimum of $g(s, \vse, y)$ with respect to $s \in [-c, c]$ is attained at $s = - c$, that is, $\mcal{M} (y, \vse) = g(-c,
\vse, y)$. Thus, we have $\vsi \in (0, \eta)$ such that
\[
| \mcal{M} (y, \vse) - \mcal{M} (z, \se)  | < \ep \qu \tx{for $| \vse - \se | < \vsi, \; |y - z| < \vsi$}.
\]
This proves the continuity of $\mcal{M} (z, \se)$  on $\se \in (u_1, u_2)$ and $z \in \bb{R}$ for Case (ii).

To show the continuity of $\mcal{M} (z, \se)$ in Case (iii), note that for  $\ep > 0$,  there exists a positive number $\eta$ less than $\min \{
u_2 - \se, \se - u_1\}$ such that
\[
| g(c, \vse, y) - g(c, \se, z) | < \ep \qu \tx{for $| \vse - \se | < \eta, \; |y - z| < \eta$},
\]
as a consequence of the continuity of $g(c, \vse, y)$ with respect to $\vse$ and $y$.   Since $\f{\pa g (s, \se, z) }{\pa s} < 0$ for $s \in
[-c, c]$ and $\f{\pa g (s, \vse, y) }{\pa s}$ is a continuous function of $s, \vse$ and $y$, there exists $\vsi \in (0, \eta)$ such that
\[
\f{\pa g (s, \vse, y) }{\pa s} < 0
\]
for $| \vse - \se | \leq \vsi, \; |y - z| \leq \vsi$ and $s \in [-c, c]$.  This implies that the minimum of $g(s, \vse, z)$ with respect to $s
\in [-c, c]$ is attained at $s = c$, that is, $\mcal{M} (y, \vse) = g(c, \vse, y)$.  Thus, we have $\vsi \in (0, \eta)$ such that
\[
| \mcal{M} (y, \vse) - \mcal{M} (z, \se)  | < \ep \qu \tx{for $| \vse - \se | < \vsi, \; |y - z| < \vsi$}.
\]
This establishes the continuity of $\mcal{M} (z, \se)$  for Case (iii).  So, we have proved that $\mcal{M} (z, \se)$ is a continuous function of
$\se \in (u_1, u_2)$ and $z \in \bb{R}$.

Clearly, $\f{\pa g (s, \se, z) }{\pa s} = z - \se + \f{\pa \psi (s, \se) }{\pa s}$ is continuous on $s$. Since $z \neq \se$ and $\f{\pa \psi (s,
\se) }{\pa s} = 0$ for $s = 0$, we have that $\f{\pa g (s, \se, z) }{\pa s} = z - \se \neq 0$ for $s = 0$. By continuity of $\f{\pa g (s, \se,
z) }{\pa s}$ on $s$, we have that $\f{\pa g (s, \se, z) }{\pa s} \neq 0$ for $s$ at a neighborhood of $0$.  Making use of this observation and
the fact $g (0, \se, z) = 0$, we have that $\mcal{M} (z, \se) < 0$.  This completes the proof of the lemma.

\epf

\beL

\la{remove asp}

 Assume that $\f{\pa \psi (s, \se) }{\pa s}$ is a continuous function of $\se \in \Se$ and $s \in \mscr{S}_\se$.
 Then, for any $\mu \in \Se$,  there exists a number $\ga > 0$ such that $(\mu - \ga, \mu + \ga) \subseteq \Se$ and that \[ \sup_{\se
\in (\mu - \ga, \mu + \ga) } \mscr{M} (\se, \mscr{U} ( \se, \vep) ) < 0, \qqu \sup_{\se \in (\mu - \ga, \mu + \ga) } \mscr{M} (\se, \mscr{L} (
\se, \vep) ) < 0.
\]

\eeL

\bpf

Clearly, for any $\mu \in \Se$,  there exists a number $\ga > 0$ such that $(\mu - \ga, \mu + \ga) \subseteq \Se$.  For such $\ga$, we claim
that $\sup_{\se \in (\mu - \ga, \mu + \ga) } \mscr{M} (\se, \mscr{U} ( \se, \vep) ) < 0$.  To prove this claim, define
\[
S = \{ \se \in (\mu - \ga, \mu + \ga):  \mscr{U} ( \se, \vep) \in \Se \}.
\]
If $S = \emptyset$, then $\mscr{M} (\se, \mscr{U} ( \se, \vep) ) = - \iy$ for any $\se \in (\mu - \ga, \mu + \ga)$.  It remains to consider the
case that $S \neq \emptyset$.  For simplicity of notations, define
\[
v_1 = \inf_{\se \in \Se} \mscr{U} ( \se, \vep), \qqu v_2 = \sup_{\se \in \Se} \mscr{U} ( \se, \vep).
\] Since $\mscr{U} ( \se, \vep)$ is continuous on $\se$, there exist $u_1$ and $u_2$ such that
\[
\inf \Se <  u_1 < v_1 \leq v_2 < u_2  < \sup \Se.
\]
Define \[ q = \inf_{ \se \in (\mu - \ga, \mu + \ga) } [ \mscr{U} ( \se, \vep) - \se ]
\]
and
\[
\ovl{a} = \max_{ \se \in [u_1, u_2] } a (\se), \qqu \udl{b} = \min_{ \se \in [u_1, u_2] } b (\se).
\]
Note that $q > 0$ and
\[
\sup_{\se \in (\mu - \ga, \mu + \ga) } \mscr{M} (\se, \mscr{U} ( \se, \vep) ) \leq \sup_{z \in [\mu - \ga, \mu + \ga] \atop{ \se \in [v_1, v_2]
 \atop{ \se - z \geq \f{q}{2}} } } \mscr{M} (z, \se ).
\]
By the continuity of $a(\se)$ and $b(\se)$ with respect to $\se \in [u_1, u_2]$, we have $\ovl{a}  < 0$ and $\udl{b} > 0$.  Hence, $\psi (s,
\se)$ is well defined for $\ovl{a} < s < \udl{b}$ for all $\se \in (u_1, u_2)$. Let $d = \min \{ - \ovl{a}, \udl{b} \}$. Then, $\f{\pa \psi (s,
\se) }{\pa s}$ is a continuous function of $\se \in (u_1, u_2)$ and $s \in (-d, d)$.   Define $\mcal{M} (z, \se) = \inf_{s \in [-c, c]}  [s (\se
- z) + \psi (s, \se) ]$ with $c \in (0, d)$.  According to Lemma \ref{lemgooday}, we have that $\mcal{M} (z, \se)$ is a continuous function of
$z \in [\mu - \ga, \mu + \ga] \subseteq \Se$ and $\se \in (u_1, u_2) \subseteq \Se$. By such continuity of $\mcal{M} (z, \se)$ and the fact that
$\mcal{M} (z, \se) < 0$ if $z \neq \se$, we have that
\[
\sup_{z \in [\mu - \ga, \mu + \ga] \atop{ \se \in [v_1, v_2]
 \atop{ \se - z \geq \f{q}{2}} } } \mscr{M} (z, \se ) \leq \sup_{z \in [\mu - \ga, \mu + \ga] \atop{ \se \in [v_1, v_2]
 \atop{ \se - z \geq \f{q}{2}} } } \mcal{M} (z, \se ) < 0,
\]
which implies that $\sup_{\se \in (\mu - \ga, \mu + \ga) } \mscr{M} (\se, \mscr{U} ( \se, \vep) ) < 0$. This proves our claim.   By a similar
argument, we can show that $\sup_{\se \in (\mu - \ga, \mu + \ga) } \mscr{M} (\se, \mscr{L} ( \se, \vep) ) < 0$.  The proof of the lemma is thus
completed.

\epf

\beL \la{uniformcon}  Let $\se_\ep$ be a continuous function of $\ep$ and $z \in (a, b) \subseteq \Se$ such that  $\se_\ep - z \to 0$ uniformly
on $z \in (a, b)$ as $\ep \to 0$. Assume that $\f{\pa^3 \psi (s, z) }{\pa s^3}$ is continuous with respect to $s$ and $z$ at a neighborhood of
$s = 0, \; z = \mu$, where $\mu \in (a, b)$. Then, for any $\vsi \in (0, 1)$, there exist $\ga
> 0$ and $\vep^*
> 0$ such that {\small $1 - \vsi < - \f{ 2 \mcal{V} (\se_\ep)  \mscr{M} (z, \se_\ep) } {  (\se_\ep - z)^2  }  < 1 + \vsi$} for all $z \in (\mu -
\ga, \mu + \ga)$ and $\vep \in (0, \vep^*)$.  \eeL

\bpf Recalling (\ref{defg3}), we have $g (s, \se_\ep, z) = s ( \se_\ep - z ) + \psi (s, \se_\ep)$.  For simplicity of notations, define \[ R(s,
\se_\ep) = \psi (s, \se_\ep) - \mcal{V} (\se_\ep) \f{s^2}{2}, \qqu r(s, \se_\ep) = \f{\pa [R(s, \se_\ep)]  } { \pa s }.
\]
Then, $\mscr{M} (z, \se_\ep) = \inf_{s \in \mscr{S}_{\se_\ep} } g (s, \se_\ep, z)$ and \bee g(s, \se_\ep, z) = (\se_\ep - z) s + \mcal{V}
(\se_\ep) \f{s^2}{2} + R(s, \se_\ep), \qqu \qu \f{\pa }{\pa s} g (s, \se_\ep, z ) = \se_\ep - z + s \mcal{V} (\se_\ep) + r(s, \se_\ep). \eee
Since $\psi (s, \se_\ep)$ is a convex function of $s$, it follows that $g (s, \se_\ep, z)$ is also a convex function of $s$.  By the assumption
of the lemma, there exist $u \in (a, \mu), \; v \in (\mu, b)$ and $c > 0$ such that $\f{\pa^3 \psi (s, z) }{\pa s^3}$ is continuous with respect
to $z \in [u, v]$ and $s \in [-c, c]$.   Clearly,  $\f{\pa }{\pa s} g (s, \se_\ep, z) = 0$ is equivalent to  \be \la{devide}
 s = \f{ z - \se_\ep }{ \mcal{V} (\se_\ep) } - \f{r(s, \se_\ep)} { \mcal{V} (\se_\ep)}.
\ee   We claim that there exist $\eta > 0$ and $\vep^* > 0$ satisfying the following requirements:

(i) $\se_\ep \in [u, v]$ for all $z \in [\mu - \eta, \mu + \eta]$ and $\ep \in (0, \vep^*)$;

(ii) There exists a unique root $s_\ep$, which is a function of $\ep \in (0, \vep^*)$ and $z \in [\mu - \eta, \mu + \eta]$ for equation
(\ref{devide}) with respect to $s$;

(iii) $s_\ep \to 0$ uniformly with respect to $z \in [\mu - \eta, \mu + \eta]$ as $\ep \downarrow 0$.

To establish the above claim, note that as a consequence of the continuity of $\se_\ep$ with respect to $z$ and $\ep$, there exist $\eta
> 0$ and $\vep_1 > 0$ such that $\se_\ep \in [u, v]$ for $z \in [\mu - \eta, \mu + \eta]$ and $\ep \in (0, \vep_1)$.
Note that $\f{\pa \psi(s, \se) }{\pa s} > 0$ for $s = c$ and $\se \in [u, v] \subseteq \Se$ because $\f{\pa \psi(s, \se) }{\pa s} = 0$ for $s =
0$ and $\psi(s, \se)$ is convex with respect to $s$.  Since $\f{\pa \psi(s, \se) }{\pa s} \mid_{s = c}$  is continuous on $\se$ and $\se_\ep \in
[u, v]$ for all $\ep \in (0, \vep_1)$, it follows that \be \la{use896a}
 \inf_{\ep \in (0, \vep_1)} \li. \f{\pa \psi(s, \se_\ep) }{\pa s} \ri
|_{s = c}
> 0.
\ee Similarly, we can show that \be \la{use896b} \sup_{\ep \in (0, \vep_1)} \li. \f{\pa \psi(s, \se_\ep) }{\pa s} \ri |_{s = - c} < 0. \ee By
virtue of (\ref{use896a}),  (\ref{use896b}) and  the assumption that $\se_\ep \to z$ uniformly with respect to $z \in (a, b)$ as $\ep \downarrow
0$, we have that there exists $\vep_2 \in (0, \vep_1)$ such that
\[
\f{\pa \psi(s, \se_\ep) }{\pa s} \mid_{s = - c} <  z - \se_\ep <  \f{\pa \psi(s, \se_\ep) }{\pa s} \mid_{s = c}
\]
for  $z \in [\mu - \eta, \mu + \eta]$ and $\ep \in (0, \vep_2)$.  As a consequence of this observation and the fact that $\f{\pa \psi(s,
\se_\ep) }{\pa s}$ is increasing on $s$, we have that for $z \in [\mu - \eta, \mu + \eta]$ and $\ep \in (0, \vep_2)$, there exists a unique
number $s_\ep \in [-c, c]$ such that $\f{\pa }{\pa s} g (s, \se_\ep, z_\ep) = 0$ for $s = s_\ep$.  To complete the proof of the claim, it
remains to show that $s_\ep \downarrow 0$ uniformly with respect to $z \in [\mu - \eta, \mu + \eta]$ as $\ep \downarrow 0$.  For this purpose,
define \[ h(s) = \bec \min_{\se \in [u, v]} \f{\pa \psi(s, \se)}{\pa s} & \tx{for} \; s \geq 0,\\
\max_{\se \in [u, v]} \f{\pa \psi(s, \se)}{\pa s} & \tx{for} \; s < 0
 \eec
\]
Clearly,  $h(s)$ is equal to $0$ for $s = 0$ and is continuous, increasing with respect to $s \in [ - c, c]$.  Moreover, there exists $\vep^*
\in (0, \vep_2)$ such that $z - \se_\ep < h(c)$ for all $\ep \in (0, \vep^*)$ and $z \in [\mu - \eta, \mu + \eta]$.  Hence, for every $\ep \in
(0, \vep^*)$ such that $z - \se_\ep > 0$, there exists $t_\ep \geq s_\ep
> 0$ such that $z - \se_\ep = h(t_\ep)$.  Since $\se_\ep \to z$ uniformly with respect to $z \in [\mu - \eta, \mu + \eta]$, it follows that
$t_\ep \downarrow 0$ uniformly with respect to $z \in [\mu - \eta, \mu + \eta]$ as $\ep \downarrow 0$. Therefore,  for $\ep \in (0, \vep^*)$
such that $z - \se_\ep > 0$, we have that $s_\ep \downarrow 0$ uniformly with respect to $z \in [\mu - \eta, \mu + \eta]$ as $\ep \downarrow 0$.
Similarly, it can be argued that there exist $\eta > 0$ and $\vep^* > 0$ such that for $\ep \in (0, \vep^*)$ with $z - \se_\ep < 0$,  $s_\ep
\uparrow 0$ uniformly with respect to $z \in [\mu - \eta, \mu + \eta]$ as $\ep \downarrow 0$. Therefore, we have shown the above claim.

Now consider $R(s, \se_\ep)$ and $r(s, \se_\ep)$. Since the third partial derivative of $\psi(t, \se_\ep)$ with respect to  $t$ is continuous
for $t \in [-c, c]$, we have \bee \psi (s, \se_\ep) & = & \psi (t, \se_\ep) |_{t = 0} + s \times \f{\pa [\psi (t, \se_\ep)]  } { \pa t } |_{t =
0} + \f{s^2}{2} \times \f{\pa^2 [\psi (t, \se_\ep)]  } { \pa t^2 } |_{t = 0}
+ \f{s^3}{6} \times  \f{\pa^3 [\psi (t, \se_\ep)]  } { \pa t^3 } |_{t = \vsi s}\\
&  = & \mcal{V} (\se_\ep) \f{s^2}{2} + \f{s^3}{6} \times  \f{\pa^3 [\psi (t, \se_\ep)]  } { \pa t^3 } |_{t = \vsi s},   \eee where $\vsi \in (0,
1)$ depends on $s$ and $\se_\ep$.  Thus, $R(s, \se_\ep) = \f{s^3}{6} \times  \f{\pa^3 [\psi (t, \se_\ep)]  } { \pa t^3 } |_{t = \vsi s}$. Since
$\f{\pa^3 [\psi (t, \se)]  } { \pa t^3 }$ is continuous with respect to $t \in [-c, c]$ and $\se \in [u, v]$,  the minimum and maximum of
$\f{\pa^3 [\psi (t, \se)]  } { \pa t^3 }$ are attained at some $t \in [-c, c]$ and $\se \in [u, v]$, which are finite. That is, there exist two
constants $A$ and $B$ independent of $\ep$ and $z$ such that $A < \f{\pa^3 [\psi (t, \se_\ep)] } { \pa t^3 } |_{t = \vsi s} < B$. It follows
that \be \la{imp89a}
 \f{A}{6} < \f{ R(s, \se_\ep) }{s^3} < \f{B}{6}
\ee for  $\ep \in (0, \vep^*), \; z \in [\mu - \eta, \mu + \eta]$ and $0 < |s| \leq c$.   Since the  third partial derivative of $\psi(t,
\se_\ep)$ with respect to  $t$ is continuous for $t \in [ - c, c]$, we have  {\small \bee \f{\pa [\psi (s, \se_\ep)] } { \pa s } = \f{\pa [\psi
(t, \se_\ep)] } { \pa t } |_{t = 0} + s \times \f{\pa^2 [\psi (t, \se_\ep)]  } { \pa t^2 } |_{t = 0}  + \f{s^2}{2} \times  \f{\pa^3 [\psi (t,
\se_\ep)]  } { \pa t^3 } |_{t = \vsi s} = \mcal{V} (\se_\ep) s + \f{s^2}{2} \times  \f{\pa^3 [\psi (t, \se_\ep)]  } { \pa t^3 } |_{t = \vsi s},
\eee}  where $\vsi \in (0, 1)$ depends on $s$ and $\se_\ep$. Therefore, \[ \f{r(s, \se_\ep)}{s^2} = \f{1}{s^2} \li [ \f{\pa [\psi (s, \se_\ep)]
} { \pa s } - \mcal{V} (\se_\ep) s \ri ]  = \f{1}{2}   \f{\pa^3 [\psi (t, \se_\ep)]  } { \pa t^3 } |_{t = \vsi s}.
  \]
It follows that \be \la{imp89b} \f{A}{2} < \f{r(s, \se_\ep)}{s^2}  < \f{B}{2} \ee for $\ep \in (0, \vep^*), \; z \in [\mu - \eta, \mu + \eta]$
and $0 < |s| \leq c$.  Making use of (\ref{devide}), we have  \bel  - \f{ 2 \mcal{V} (\se_\ep)  \mscr{M} (z, \se_\ep) } { (\se_\ep - z)^2 } & =
& - \f{ 2 \mcal{V} (\se_\ep) } { (\se_\ep - z)^2 } \li [ (\se_\ep - z) s_\ep + \mcal{V}
(\se_\ep) \f{s_\ep^2}{2} +  R(s_\ep, \se_\ep) \ri ] \nonumber\\
&  = &  \li [ 1 + \f{r(s_\ep, \se_\ep)}{s_\ep \mcal{V} (\se_\ep)} \ri ]^{-2} \li [ 1 + 2 \f{r(s_\ep, \se_\ep)}{s_\ep \mcal{V} (\se_\ep)} - \f{2
s_\ep}{ \mcal{V} (\se_\ep) } \f{ R(s_\ep, \se_\ep) }{ s_\ep^3} \ri ] \la{onestep}  \eel Clearly, $0 < \min_{\se \in [u, v]} \mcal{V} (\se) \leq
\mcal{V} (\se_\ep) \leq \max_{\se \in [u, v]} \mcal{V} (\se)$ for  $z \in [\mu - \eta, \mu + \eta]$ and $\ep \in (0, \vep^*)$.  Recall that
$s_\ep \to 0$ uniformly on $z \in [\mu - \eta, \mu + \eta]$ as $\ep \downarrow 0$.  Using these facts and (\ref{imp89a}), (\ref{imp89b}) and
(\ref{onestep}) concludes the proof of the lemma.

\epf

\section{Proof of Theorem \ref{ThmCDF} } \la{ThmCDF_app}

For simplicity of notations, define
\[
\mscr{W} (\mu) = \bb{E} [ | X - \mu |^3 ], \qqu  \mcal{B} (\mu) = \f{ \bb{E} [ | X - \mu |^3 ] }{ \bb{E}^{ \f{3}{2} } [ | X - \mu |^2 ] }, \]
which are continuous functions of $\mu \in \Se$.

\subsection{Proof of Property (I)}

To prove (\ref{prop1}), we need the following preliminary result.

\beL \la{asn889} For $\mu \in \Se$, there exist a number $\ga > 0$ and a positive integer $m > 0$  such that $\{ \mbf{n} > n \} \subseteq \{ |
Y_n - \mu | \geq \ga   \}$ for all $n > m$.

\eeL

\bpf

To show the lemma, we need to prove two claims as follows.

(i): For $\mu \in \Se$, there exist a number $\ga > 0$ and  a
positive integer $m > 0$  such that
\[
\li \{  F_{Y_n} ( Y_n, \mscr{U} (Y_n, \vep) ) > \f{\de_n}{2} \ri \} \subseteq \{  | Y_n - \mu | \geq \ga   \} \qu \tx{for all $n > m$}.
\]

(ii): For $\mu \in \Se$, there exist a number $\ga > 0$ and  a
positive integer $m > 0$  such that
\[
\li \{  G_{Y_n} ( Y_n, \mscr{L} (Y_n, \vep) ) > \f{\de_n}{2} \ri \} \subseteq \{  | Y_n - \mu | \geq \ga   \} \qu \tx{for all $n > m$}.
\]

To prove the first claim, note that there exists $\ga > 0$ such that  $\mu - \ga, \; \mu + \ga \in \Se$ and that $\mscr{U} (\se, \vep) - \se
\geq \f{ \mscr{U} (\mu, \vep) - \mu }{2} > 0$ for all $\se \in (\mu - \ga, \mu + \ga) \subseteq \Se$.  Define $S = \{ \se \in (\mu - \ga, \mu +
\ga): \mscr{U} (\se, \vep) \in \Se \}$. If $S = \emptyset$, then $\li \{ F_{Y_n} ( Y_n, \mscr{U} (Y_n, \vep) ) > \f{\de_n}{2}, \; | Y_n - \mu |
< \ga \ri \} = \emptyset$ for $n \in \bb{N}$.  In the other case that $S \neq \emptyset$, we have that there exists a constant $D < \iy$ such
that
\[
\sup_{ z \in S } \f{\mscr{W}(\mscr{U} (z, \vep) )}{| z - \mscr{U} ( z, \vep)  |^3} < D
\]
It follows from (\ref{beryb1}) that
\bee &  &  \li \{  F_{Y_n} ( Y_n, \mscr{U} (Y_n, \vep) ) > \f{\de_n}{2}, \; | Y_n - \mu | < \ga  \ri \}\\
&  & = \li \{ F_{Y_n} ( Y_n, \mscr{U} (Y_n, \vep) ) > \f{\de_n}{2}, \; \mscr{U} (Y_n, \vep) \in \Se, \; |
Y_n - \mu | < \ga  \ri \}\\
&  & \subseteq \li \{ \f{1}{2} \exp \li ( - \f{ n | Y_n - \mscr{U} (Y_n, \vep) |^2 }{  2 \mcal{V} (\mscr{U} (Y_n, \vep))   } \ri ) + \f{C
\mscr{W}(\mscr{U} (Y_n, \vep) )}{n^2 | Y_n - \mscr{U} (Y_n, \vep)  |^3} > \f{\de_n}{2}, \;
 \mscr{U} (Y_n, \vep) \in \Se, \; | Y_n - \mu | < \ga \ri \} \\
&  & \subseteq \li \{  \f{1}{2} \exp \li ( - \f{ n | Y_n - \mscr{U} (Y_n, \vep) |^2 }{  2 \mcal{V} (\mscr{U} (Y_n, \vep))   } \ri )
+ \f{C D}{n^2} > \f{\de_n}{2}, \; \mscr{U} (Y_n, \vep) \in \Se, \; | Y_n - \mu | < \ga \ri \} \\
&  & \subseteq \li \{  \exp \li ( - \f{ n | Y_n - \mscr{U} (Y_n, \vep) |^2 }{  2 \mcal{V} (\mscr{U} (Y_n, \vep))   } \ri )  > \f{\de}{4}, \;
 \mscr{U} (Y_n, \vep) \in \Se, \; | Y_n - \mu | < \ga \ri \} \\
&  & \subseteq \li \{  \exp \li ( - \f{ n | Y_n - \mscr{U} (Y_n, \vep) |^2 }{ 2 \mcal{V}^* } \ri )
> \f{\de}{4}, \; | Y_n - \mu | < \ga \ri \} \\
&  & \subseteq \li \{  \mscr{U} (Y_n, \vep) - Y_n < \f{ \mscr{U} (\mu, \vep) - \mu  }{2}, \; | Y_n - \mu | < \ga \ri \} = \emptyset \eee for
large enough $n$, where $\mcal{V}^* = \sup_{\se \in S} \mcal{V} (\se)$.  Hence, in both cases that $S = \emptyset$ and $S \neq \emptyset$, \bee
\li \{ F_{Y_n} ( Y_n, \mscr{U} (Y_n, \vep) ) > \f{\de_n}{2} \ri \} & \subseteq & \li \{ F_{Y_n} ( Y_n, \mscr{U} (Y_n, \vep) )
> \f{\de_n}{2}, \; | Y_n - \mu | < \ga  \ri \}  \cup \{ |
Y_n - \mu | \geq \ga \}\\
& = & \{ | Y_n - \mu | \geq \ga \}  \eee for large enough $n$. This establishes the first claim.   Similarly, we can show the second claim.

By the definition of the stopping rule, we have that {\small \be
\la{observa8899}
 \{ \mbf{n} > n \} \subseteq \li \{  F_{Y_n} (
Y_n, \mscr{U} (Y_n, \vep) ) > \f{\de_n}{2} \ri \} \cup \li \{ G_{Y_n} ( Y_n, \mscr{L} (Y_n, \vep) ) > \f{\de_n}{2} \ri \}. \ee} Finally, the
proof of the lemma is completed by making use of (\ref{observa8899}) and our two established claims.

\epf

We are now in a position to prove (\ref{prop1}). According to Lemma \ref{asn889}, there exist an integer $m$ and a positive number $\ga > 0$
such that \be \la{makeuse}
 \Pr \{ \mbf{n} > n \}
\leq \Pr \{ | Y_n - \mu | \geq \ga \} \qqu \tx{for all $n > m$}. \ee Making use of (\ref{makeuse}) and the weak law of large numbers, we have
\[ 0 \leq  \limsup_{n \to \iy} \Pr \{ \mbf{n}
> n \} \leq \lim_{n \to \iy} \Pr \{ | Y_n - \mu | \geq \ga \} = 0,
\]
which implies that $\lim_{n \to \iy} \Pr \{ \mbf{n} > n \} = 0$.  It follows that
\[
\Pr \{ \mbf{n} < \iy \} = 1 - \lim_{n \to \iy} \Pr \{ \mbf{n} > n \} = 1.
\]
By virtue of (\ref{mostvip}) of Lemma \ref{nonuniformBE} and (\ref{makeuse}), we have \bee  \bb{E} [  \mbf{n}  ]  \leq  m + \sum_{n > m} \Pr \{
| Y_n - \mu | \geq \ga \} \leq  m + \sum_{n > m}  \li [ \exp \li ( - \f{n}{2} \f{ \ga^2 }{  \mcal{V} (\mu) } \ri ) + \f{2 C }{n^2} \f{
\mscr{W}(\mu) }{ \ga^3 } \ri ] < \iy.  \eee This completes the proof of (\ref{prop1}).

\subsection{Proof of Property (II)}

To prove (\ref{prop3}),  define
\[
\bs{L}_{\mbf{n}} = \mscr{L} ( Y_\mbf{n}, \vep ), \qqu \bs{L}_{\mbf{n} - 1} = \mscr{L} ( Y_{\mbf{n} - 1}, \vep ), \qqu \bs{U}_{\mbf{n}} =
\mscr{U} ( Y_\mbf{n}, \vep ), \qqu \bs{U}_{\mbf{n} - 1} = \mscr{U} ( Y_{\mbf{n} - 1}, \vep )
\]
and {\small \bee &  & \mscr{E}_1 = \li \{ F_{Y_\mbf{n}} (Y_\mbf{n}, \bs{U}_{\mbf{n}} ) \leq \f{\de_{\mbf{n}}}{2}, \;  G_{Y_\mbf{n}} (Y_\mbf{n},
\bs{L}_{\mbf{n}}
) \leq \f{\de_{\mbf{n}}}{2} \ri \}, \\
&  & \mscr{E}_2 = \{  \mbf{n} = m_\vep  \} \bigcup \li \{ F_{Y_{\mbf{n}-1}} (Y_{\mbf{n}-1}, \bs{U}_{\mbf{n}-1} ) > \f{\de_{\mbf{n}-1}}{2}, \;  \mbf{n} > m_\vep  \ri \}
 \bigcup \li \{ G_{Y_{\mbf{n}-1}} (Y_{\mbf{n}-1}, \bs{L}_{\mbf{n}-1} ) > \f{\de_{\mbf{n}-1}}{2}, \;  \mbf{n} > m_\vep  \ri \}, \\
&  & \mscr{E}_3 = \li \{ \lim_{\vep \downarrow 0} \mbf{n} = \iy, \; \lim_{\vep \downarrow 0} Y_{\mbf{n}} = \mu, \; \lim_{\vep \downarrow 0}
Y_{\mbf{n} - 1} = \mu \ri \}, \\
&  & \mscr{E}_4 = \mscr{E}_1 \cap \mscr{E}_2 \cap  \mscr{E}_3, \\
&  &  \mscr{E}_5 = \li \{ \lim_{\vep \downarrow 0} \f{ \mbf{n} }{ \mcal{N} (\vep, \de, \mu) } = 1 \ri \}.  \eee} By the assumption that $m_\vep
\to \iy$ as $\vep \downarrow 0$, we have that $\{ \lim_{\vep \downarrow 0} \mbf{n} = \iy \}$ is a sure event.  It follows from the strong law of
large numbers that $\mscr{E}_3$ is an almost sure event. By the definition of the stopping rule, we have that
 $\mscr{E}_1 \cap \mscr{E}_2$ is an almost sure event. Hence, $ \mscr{E}_4$ is an almost sure event. To show (\ref{prop3}), it suffices to show that
$\mscr{E}_4 \subseteq \mscr{E}_5$. For this purpose, we let $\om \in \mscr{E}_4$ and expect to show that $\om \in \mscr{E}_5$. For simplicity of
notations, let \[ n = \mbf{n} (\om), \qqu y_n = Y_\mbf{n} (\om), \qqu L_n = \mscr{L} ( y_n, \vep), \qqu U_n = \mscr{U} ( y_n, \vep), \]
\[ y_{n-1} = Y_{\mbf{n} - 1} (\om), \qqu L_{n-1} = \mscr{L} ( y_{n-1}, \vep), \qqu U_{n-1} = \mscr{U} ( y_{n-1}, \vep).  \]
Let $\eta > 0$ be small enough so that $\mu + \eta, \; \mu - \eta \in \Se$.  Then, there exists $\vep^\star > 0$ such that
\[ | L_n - \mu | < \eta, \qu | U_n- \mu | < \eta, \qu | L_{n-1} - \mu | < \eta, \qu
| U_{n-1} - \mu | < \eta \] for all $\vep \in (0, \vep^\star)$.  It follows from the continuity of $\mcal{B}(.)$ that there exists a constant $K
> 0$ such that
\[
\mcal{B} (L_n) < K, \qu \mcal{B} (U_n) < K, \qu \mcal{B} (L_{n-1}) < K, \qu \mcal{B} (U_{n-1}) < K
\]
for all $\vep \in (0, \vep^\star)$.  Since $\om \in \mscr{E}_1$, we have
\[
F_{Y_n} (y_n, U_n) \leq \f{\de_n}{2},  \qqu G_{Y_n} (y_n, L_n ) \leq \f{\de_n}{2}.
\]
Applying (\ref{berry55}), (\ref{berry56}) of Lemma \ref{nonuniformBE} and these inequalities, we have
\[
\Phi \li ( \f{ \sq{n} [  y_n - U_n ]   }{  \sq{ \mcal{V} ( U_n) } } \ri ) - \f{ \mcal{B} (  U_n) }{\sq{n}} \leq \f{\de_n}{2}, \qqu \Phi \li (
\f{ \sq{n} [ L_n - y_n  ]   }{  \sq{ \mcal{V} ( L_n ) } } \ri ) - \f{ \mcal{B} ( L_n ) }{\sq{n}} \leq \f{\de_n}{2}
\]
which can be written as
\[
\f{n} { \mcal{N} (\vep, \de, \mu)  }  \geq \f{1}{\mcal{Z}^2} \li [  \Phi^{-1}   \li ( 1 -  \f{\de_n}{2} - \f{ \mcal{B} (  U_n) }{\sq{n}} \ri )
\ri ]^2  \li [ \f{ \ka (\mu) \vep }{  U_n- y_n  } \ri ]^2 \f{ \mcal{V} ( U_n) } { \mcal{V} (\mu) },
\]
\[
\f{n} { \mcal{N} (\vep, \de, \mu)  }  \geq \f{1}{\mcal{Z}^2} \li [ \Phi^{-1}   \li ( 1 -  \f{\de_n}{2} - \f{ \mcal{B} (  L_n ) }{\sq{n}} \ri )
\ri ]^2  \li [ \f{ \ka (\mu) \vep }{ y_n  - L_n  } \ri ]^2 \f{ \mcal{V} ( L_n ) } { \mcal{V} (\mu) }.
\]
These two inequalities imply that \be \la{infeq} \liminf_{\vep \downarrow 0} \f{n} { \mcal{N} (\vep, \de, \mu)  }  \geq 1. \ee As a consequence
of (\ref{infeq}) and the assumption that $\lim_{\vep \downarrow 0} \vep^2 m_\vep = 0$, it must be true that $n \neq m_\vep$ for small enough
$\vep \in (0, \vep^\star)$.  By virtue of this observation and $\om \in \mscr{E}_2$, we have that either
\[
F_{Y_{n-1}} (y_{n-1}, U_{n-1} ) > \f{\de_{n-1}}{2} \qqu \tx{or} \qqu  G_{Y_{n-1}} (y_{n-1}, L_{n-1} ) > \f{\de_{n-1}}{2}
\]
must be true for small enough $\vep \in (0, \vep^\star)$. It follows from (\ref{berry55}) and (\ref{berry56}) that either
\[
\Phi \li ( \f{ \sq{{n - 1}} [  y_{n - 1} - U_{n-1}  ]   }{  \sq{ \mcal{V} ( U_{n-1} ) } } \ri ) + \f{\mcal{B} (  U_{n-1} ) }{\sq{{n - 1}}} >
\f{\de_{n-1}}{2}
\]
or
\[
\Phi \li ( \f{ \sq{{n - 1}} [ L_{n - 1} - y_{n - 1}  ]   }{  \sq{ \mcal{V} ( L_{n - 1} ) } } \ri ) + \f{ \mcal{B} (  L_{n - 1} ) }{\sq{{n - 1}}}
> \f{\de_{n-1}}{2}
\]
must be true for small enough $\vep \in (0, \vep^\star)$. Hence, either
\[
\f{n} { \mcal{N} (\vep, \de, \mu)  }    < \f{1}{\mcal{Z}^2} \li [ \Phi^{-1}   \li ( 1 -  \f{\de_{n-1}}{2} - \f{ \mcal{B} (  U_{n-1} ) }{\sq{n}}
\ri ) \ri ]^2  \li [ \f{ \ka (\mu) \vep }{ U_{n-1} - y_{n - 1} } \ri ]^2 \f{ n \mcal{V} ( U_{n-1} ) } { (n - 1) \mcal{V} (\mu) }
\]
or
\[
\f{n} { \mcal{N} (\vep, \de, \mu)  }   < \f{1}{\mcal{Z}^2} \li [ \Phi^{-1}   \li ( 1 -  \f{\de_{n-1}}{2} - \f{ \mcal{B} (  L_{n - 1} ) }{\sq{n}}
\ri ) \ri ]^2  \li [ \f{ \ka (\mu) \vep }{ y_{n - 1}  - L_{n - 1}  }  \ri ]^2 \f{ n \mcal{V} ( L_{n - 1} ) } { (n - 1) \mcal{V} (\mu) }
\]
must be true for small enough $\vep \in (0, \vep^\star)$. Making use of these two inequalities, we have that \be \la{supeq} \limsup_{\vep
\downarrow 0} \f{n} { \mcal{N} (\vep, \de, \mu) } \leq 1. \ee Combining (\ref{infeq}) and (\ref{supeq}) yields
\[ \lim_{\vep \downarrow 0} \f{ \mbf{n} }{ \mcal{N} (\vep, \de, \mu) } = 1.
\]
Thus, we have shown $\mscr{E}_4 \subseteq \mscr{E}_5$, which implies (\ref{prop3}).

\subsection{Proof of Property (III)}

For purpose of simplifying notations, define
\[
\bs{\mcal{L}}_{\mbf{n}} = \mcal{L} ( Y_\mbf{n}, \vep), \qu \bs{\mcal{U}}_{\mbf{n}} = \mcal{U} ( Y_\mbf{n}, \vep).
\]
To prove (\ref{prop4}), we need to show the following preliminary results.

\beL

\la{use9668}
 \bel &  & \Pr \li \{
 \lim_{\vep \downarrow 0} \f{ \sq{ \mbf{n} } [ Y_\mbf{n} - \bs{\mcal{U}}_{\mbf{n}} ]
} { \sq{ \mcal{V} (\mu) } } = - \mcal{Z}   \ri \}  = 1, \la{consistency1}\\
&  & \Pr \li \{ \lim_{\vep \downarrow 0} \f{ \sq{ \mbf{n} } [ \bs{\mcal{L}}_{\mbf{n}} - Y_\mbf{n} ] } { \sq{ \mcal{V} (\mu) } } = - \mcal{Z}
\ri \}  = 1.  \la{consistency2} \eel

\eeL

\bpf  Let $\mscr{E}_i, \; i = 1, \cd, 5$ be defined as before in the proof of property (II).  Define
\[
\mscr{E}_6 = \li \{
 \lim_{\vep \downarrow 0} \f{ \sq{ \mbf{n} } [ Y_\mbf{n} - \bs{\mcal{U}}_{\mbf{n}} ]
} { \sq{ \mcal{V} (\mu) } } = - \mcal{Z}   \ri \}.
\]
To show equation (\ref{consistency1}), it suffices to show $\mscr{E}_4 \subseteq \mscr{E}_6$.  For this purpose, we let $\om \in \mscr{E}_4$ and
expect to show that $\om \in \mscr{E}_6$.

For simplicity of notations, let $n = \mbf{n} (\om), \;  y_n = Y_\mbf{n} (\om), \; y_{n-1} = Y_{\mbf{n}-1} (\om)$, \bee &  &  L_n = \mscr{L} (
y_n, \vep), \qu U_n = \mscr{U} (
y_n, \vep), \qu \mcal{L}_n = \mcal{L} ( y_n, \vep), \qu \mcal{U}_n = \mcal{U} ( y_n, \vep), \\
&  & L_{n-1} = \mscr{L} ( y_{n-1}, \vep), \qu U_{n-1} = \mscr{U} ( y_{n-1}, \vep), \qu \mcal{L}_{n-1} = \mcal{L} ( y_{n-1}, \vep), \qu
\mcal{U}_{n-1} = \mcal{U} ( y_{n-1}, \vep).  \eee Let $\eta > 0$ be small enough so that $\mu + \eta, \; \mu - \eta \in \Se$.  Then, there
exists $\vep^\star
> 0$ such that
\[ | L_n - \mu | < \eta, \qu | U_n- \mu | < \eta, \qu | L_{n-1} - \mu | < \eta, \qu
| U_{n-1} - \mu | < \eta \] for all $\vep \in (0, \vep^\star)$.  It follows from the continuity of $\mcal{B}(.)$ that there exists a constant $K
> 0$ such that
\[
\mcal{B} (L_n) < K, \qu \mcal{B} (U_n) < K, \qu \mcal{B} (L_{n-1}) < K, \qu \mcal{B} (U_{n-1}) < K
\]
for all $\vep \in (0, \vep^\star)$.  Since $\om \in \mscr{E}_1$, we have $F_{Y_n} (y_n, U_n) \leq \f{\de_n}{2}$. Making use of this inequality
and (\ref{berry55}), we have that
\[
\f{ \sq{n} [ \mcal{U}_n - y_n ] }{  \sq{ \mcal{V} (\mu) } } \geq  \Phi^{-1}   \li ( 1 -  \f{\de_n}{2} - \f{ \mcal{B} ( U_n) }{\sq{n}} \ri ) \sq{
\f{ \mcal{V} \li ( U_n\ri ) } { \mcal{V} (\mu) } } \f{ \mcal{U}_n - y_n }{ U_n- y_n},
\]
which implies that
\[
\liminf_{\vep \downarrow 0} \f{ \sq{n} [ \mcal{U}_n - y_n ] }{  \sq{ \mcal{V} (\mu) } } \geq \mcal{Z}.
\]
On the other hand, as a consequence of $\om \in \mscr{E}_2$ and the fact that  $n \neq m_\vep$ for small enough $\vep \in (0, \vep^\star)$, we
have  that either $F_{Y_{n-1}} (y_{n-1}, U_{n-1} ) > \f{\de_{n-1}}{2}$ or $G_{Y_{n-1}} (y_{n-1}, L_{n-1} )
> \f{\de_{n-1}}{2}$ must be true for small enough $\vep \in (0, \vep^\star)$.  Making use of this observation and (\ref{berry55}), (\ref{berry56}),
we have that either {\small \bee &  &  \f{ \sq{n} [ \mcal{U}_n - y_{ n } ] }{  \sq{ \mcal{V} (\mu)  }  }  < \Phi^{-1} \li ( 1 - \f{\de_{n-1}}{2}
- \f{ \mcal{B} (  U_{n-1} ) }{\sq{n - 1}} \ri )  \f{ U_n- y_{ n } }{ U_{n-1}- y_{ n - 1 } } \sq{ \f{ n \mcal{V} \li ( U_{n-1}\ri )  } { (n - 1)
\mcal{V} (\mu) } } \f{ \mcal{U}_n - y_n }{
U_n- y_n} \qu \tx{or}\\
&  & \f{ \sq{n} [ \mcal{U}_n - y_{ n  } ] }{  \sq{ \mcal{V} (\mu)  }  }  < \Phi^{-1} \li ( 1 -  \f{\de_{n-1} }{2} - \f{ \mcal{B} ( L_{n - 1} )
}{\sq{n - 1}} \ri )   \f{ U_n- y_{ n  }  }{ y_{ n - 1 } - L_{n - 1} } \sq{ \f{ n \mcal{V} \li ( L_{n - 1}\ri ) } { (n - 1) \mcal{V} (\mu) } }
\f{ \mcal{U}_n - y_n }{ U_n- y_n} \eee} must be true for small enough $\vep \in (0, \vep^\star)$. This implies that
\[
\limsup_{\vep \downarrow 0} \f{ \sq{n} [ \mcal{U}_n - y_n ] }{  \sq{ \mcal{V} (\mu) } } \leq \mcal{Z}.
\]
Hence,
\[
\lim_{\vep \downarrow 0} \f{ \sq{n} [ \mcal{U}_n - y_n ] }{  \sq{ \mcal{V} (\mu) } } = \mcal{Z}.
\]
This shows $\om \in \mscr{E}_6$ and thus $\mscr{E}_4 \subseteq \mscr{E}_6$. It follows that $\Pr \{  \mscr{E}_6  \} = 1$, which implies
(\ref{consistency1}).  Similarly, we can show (\ref{consistency2}).  This completes the proof of the lemma.

\epf

We are now in a position to prove (\ref{prop4}). Note that \bee \Pr \li \{ \bs{\mcal{U}}_{\mbf{n}} \leq \mu \ri \} & = & \Pr \li \{  \f{ \sq{
\mbf{n} } ( Y_\mbf{n} - \mu ) } { \sq{ \mcal{V} (\mu) } }  \leq \f{ \sq{ \mbf{n} } [ Y_\mbf{n} - \bs{\mcal{U}}_{\mbf{n}} ]
} { \sq{ \mcal{V} (\mu) } } \ri \}\\
& \leq & \Pr \li \{  \f{ \sq{ \mbf{n} } ( Y_\mbf{n} - \mu ) } { \sq{ \mcal{V} (\mu) } }  \leq - \mcal{Z} ( 1 - \eta) \ri \}  + \Pr \li \{
 \f{ \sq{ \mbf{n} } [ Y_\mbf{n} - \bs{\mcal{U}}_{\mbf{n}} ]
} { \sq{ \mcal{V} (\mu) } } \notin [ - \mcal{Z} ( 1 + \eta), - \mcal{Z} ( 1 - \eta) ]  \ri \} \eee and that \bee \Pr \li \{
\bs{\mcal{U}}_{\mbf{n}} \leq \mu \ri \}  \geq  \Pr \li \{  \f{ \sq{ \mbf{n} } ( Y_\mbf{n} - \mu ) } { \sq{ \mcal{V} (\mu) } }  \leq - \mcal{Z} (
1 + \eta) \ri \}  - \Pr \li \{ \f{ \sq{ \mbf{n} } [ Y_\mbf{n} - \bs{\mcal{U}}_{\mbf{n}} ] } { \sq{ \mcal{V} (\mu) } } \notin [ - \mcal{Z} ( 1 +
\eta), - \mcal{Z} ( 1 - \eta) ]  \ri \}. \eee From Anscombe's random central limit theorem, we have \bee &  & \lim_{\vep \downarrow 0} \Pr \li
\{  \f{ \sq{ \mbf{n} } ( Y_\mbf{n} - \mu
) } { \sq{ \mcal{V} (\mu) } }  \leq - \mcal{Z} ( 1 - \eta) \ri \}  = \Phi ( - \mcal{Z} ( 1 - \eta)  ),\\
&  & \lim_{\vep \downarrow 0} \Pr \li \{  \f{ \sq{ \mbf{n} } ( Y_\mbf{n} - \mu ) } { \sq{ \mcal{V} (\mu) } }  \leq - \mcal{Z} ( 1 + \eta) \ri \}
= \Phi ( - \mcal{Z} ( 1 + \eta)  ). \eee  From (\ref{consistency1}) of Lemma \ref{use9668}, we have
\[
\lim_{\vep \downarrow 0} \Pr \li \{
 \f{ \sq{ \mbf{n} } [ Y_\mbf{n} - \bs{\mcal{U}}_{\mbf{n}} ]
} { \sq{ \mcal{V} (\mu) } } \notin [ - \mcal{Z} ( 1 + \eta), - \mcal{Z} ( 1 - \eta) ]  \ri \} = 0.
\]
Therefore,
\[
\Phi ( - \mcal{Z} ( 1 + \eta)  ) \leq  \lim_{\vep \downarrow 0} \Pr \li \{ \bs{\mcal{U}}_{\mbf{n}} \leq \mu \ri \} \leq \Phi ( - \mcal{Z} ( 1 -
\eta)  )
\]
for any $\eta \in (0, 1)$.  Since the above inequalities hold for any $\eta \in (0, 1)$, it follows that {\small $\lim_{\vep \downarrow 0} \Pr
\{ \bs{\mcal{U}}_{\mbf{n}} \leq \mu \} = \Phi ( - \mcal{Z} ) = \f{\de}{2}$}.  By a similar argument, we can show that {\small $\lim_{\vep
\downarrow 0} \Pr \{ \bs{\mcal{L}}_{\mbf{n}} \geq \mu  \} = \f{\de}{2}$}.  It follows that {\small $\lim_{\vep \downarrow 0} \Pr \{
\bs{\mcal{L}}_{\mbf{n}} < \mu < \bs{\mcal{U}}_{\mbf{n}} \} = 1 - \de$}.

\subsection{Proof of Property (IV)}

To prove (\ref{prop5}), we need some preliminary results.

\beL  \la{howgood966} For any $\eta \in (0, 1)$, there exist $\ga > 0$ and $\vep^* > 0$ such that $\{ \mbf{n} > m \} \subseteq \{ | Y_m - \mu |
\geq \ga \}$ for all $m \geq (1 + \eta) N$ with $\vep \in (0, \vep^*)$, where $N = \mcal{N} (\vep, \de, \mu)$.

\eeL

\bpf

Clearly,  for given $\eta \in (0, 1)$,  there exists $\vsi \in (0, 1)$ such that $\li ( \f{1 - \vsi}{1 + \vsi} \ri )^2 ( 1 + \eta )
> 1$.  As a consequence of $\ka(\mu) > 0$ and the continuity assumption associated with (\ref{uniformasp}),
there exist $\ga_1 > 0$ and $\vep_1 > 0$ such that
\[
\mscr{L} (z, \vep) \in \Se, \qqu \mscr{U} (z, \vep) \in \Se
\]
and that
{\small \bel &  & ( 1 - \vsi) \ka( \mu ) < \f{ \mscr{U} (z, \vep) - z }{\vep} <  ( 1 + \vsi) \ka( \mu ), \la{cita} \\
&   & ( 1 - \vsi) \ka( \mu ) < \f{ z - \mscr{L} (z, \vep) }{\vep} < ( 1 + \vsi) \ka( \mu ) \la{citb} \eel} for all $z \in (\mu - \ga_1, \mu +
\ga_1) \subseteq \Se$ and $\vep \in (0, \vep_1)$.  Since $\mcal{V} \li ( \mscr{U} (z, \vep)  \ri )$ and $\mcal{V} \li ( \mscr{L} (z, \vep)  \ri
)$ are continuous functions of $z$ and $\vep$ which tend to $\mcal{V}(\mu) > 0$ as $(z, \vep) \to (\mu, 0)$, there exist $\ga \in (0, \ga_1)$
and $\vep_2 \in (0, \vep_1)$ such that \be \la{need9}
 \mcal{V} \li ( \mscr{U} (z, \vep)  \ri ) <  \li ( \f{1 - \vsi}{1 + \vsi} \ri )^2 ( 1 + \eta ) \mcal{V} (\mu), \qqu  \mcal{V} \li (
\mscr{L} (z, \vep) \ri ) < \li ( \f{1 - \vsi}{1 + \vsi} \ri )^2 ( 1 + \eta ) \mcal{V} (\mu) \ee for all $z \in (\mu - \ga, \mu + \ga)$ and $\vep
\in (0, \vep_2)$.  It follows from (\ref{cita}), (\ref{citb}) and (\ref{need9}) that \bel &  &
 \f{ \sq{{ (1 + \eta) N }} [ \mscr{U} (z, \vep) - z ]  }{ \sq{ \mcal{V} \li ( \mscr{U} (z, \vep) \ri ) } } \geq
 \f{ \sq{{ (1 + \eta) N }} ( 1 - \vsi) \ka( \mu ) \vep }{ \sq{ \mcal{V} \li ( \mscr{U} (z, \vep) \ri ) } } \geq (1 + \vsi) \mcal{Z}, \la{oka} \\
&  & \f{ \sq{{ (1 + \eta) N }} [ z - \mscr{L} (z, \vep) ]  }{ \sq{ \mcal{V} \li ( \mscr{L} (z, \vep) \ri ) } } \geq
  \f{ \sq{{ (1 + \eta) N }} ( 1 - \vsi) \ka( \mu )  \vep }{  \sq{ \mcal{V} \li ( \mscr{L} (z, \vep) \ri ) } } \geq
(1 + \vsi) \mcal{Z} \la{okb} \eel for all $z \in (\mu - \ga, \mu + \ga)$ and $\vep \in (0, \vep_2)$. Observe that
\[
\Phi^{-1} \li ( 1 -  \f{\de_m}{2} + \f{\mcal{B} ( \mscr{U} (z, \vep)  ) }{\sq{m}} \ri ) \qqu \tx{and} \qqu \Phi^{-1} \li ( 1 - \f{\de_m}{2} +
\f{\mcal{B} ( \mscr{L} (z, \vep)  ) }{\sq{m}} \ri )
\] are continuous functions of $z, \vep$ and $m$.  Since $\mscr{L} (z, \vep)$ and $\mscr{U} (z, \vep)$ are bounded and contained in $\Se$
for all $z \in (\mu - \ga, \mu + \ga)$ and $\vep \in (0, \vep_2)$, it follows from the continuity of $\mscr{L} (., \vep), \; \mscr{U} (., \vep)$
and $\mcal{B}(.)$ that $\mcal{B} ( \mscr{L} (z, \vep)  )$ and $\mcal{B} ( \mscr{U} (z, \vep)  )$ are bounded for all $z \in (\mu - \ga, \mu +
\ga)$ and $\vep \in (0, \vep_2)$. By such boundedness and the fact that $\de_m \to \de$ as $\vep \downarrow 0$,
 there exists $\vep^* \in (0, \vep_2)$ such that
 \bel &  & \Phi^{-1} \li ( 1 -  \f{\de_m}{2} + \f{\mcal{B} ( \mscr{U} (z, \vep)  ) }{\sq{m}} \ri )
  < (1 + \vsi) \mcal{Z}, \la{okc} \\
  &  &  \Phi^{-1} \li ( 1 - \f{\de_m}{2} + \f{\mcal{B} ( \mscr{L} (z, \vep)  ) }{\sq{m}} \ri )
   < (1 + \vsi) \mcal{Z} \la{okd}
\eel for all $z \in (\mu - \ga, \mu + \ga)$ and $m \geq (1 + \eta) N$ with $\vep \in (0, \vep^*)$.  Combing (\ref{oka}), (\ref{okb}),
(\ref{okc}) and (\ref{okd}) yields
\[
 \f{ \sq{{ m }} [\mscr{U} (z, \vep) - z]  }{ \sq{ \mcal{V} \li ( \mscr{U} (z, \vep) \ri ) } }
 \geq \Phi^{-1} \li ( 1 -  \f{\de_m}{2} + \f{\mcal{B} ( \mscr{U} (z, \vep)  ) }{\sq{m}} \ri ), \qqu
  \f{ \sq{{ m }} [z - \mscr{L} (z, \vep) ] }{  \sq{ \mcal{V} \li ( \mscr{L} (z, \vep) \ri ) } }
  \geq \Phi^{-1} \li ( 1 - \f{\de_m}{2} + \f{\mcal{B} ( \mscr{L} (z, \vep)  ) }{\sq{m}} \ri )
\]
for all $z \in (\mu - \ga, \mu + \ga)$ and $m \geq (1 + \eta) N$ with $\vep \in (0, \vep^*)$.  This implies that that {\small \bee \{ | Y_m -
\mu | < \ga \} &  \subseteq & \li \{ \f{ \sq{{ m }} [ \mscr{U} (Y_m, \vep) - Y_m   ] }{ \sq{ \mcal{V} \li ( \mscr{U} (Y_m, \vep) \ri ) } } \geq
\Phi^{-1} \li ( 1 -  \f{\de_m}{2} + \f{\mcal{B} ( \mscr{U} (Y_m, \vep)  ) }{\sq{m}} \ri ), \; \mscr{U} (Y_m, \vep) \in \Se \ri \}\\
&  & \bigcap \li \{   \f{ \sq{{ m }} [ Y_m - \mscr{L} (Y_m, \vep)  ] }{ \sq{ \mcal{V} \li ( \mscr{L} (Y_m, \vep) \ri ) } } \geq \Phi^{-1} \li (
1 -  \f{\de_m}{2} + \f{\mcal{B} ( \mscr{L} (Y_m, \vep)  ) }{\sq{m}} \ri ), \; \mscr{L} (Y_m, \vep) \in \Se \ri \} \eee} for all $m \geq (1 +
\eta) N$ with $\vep \in (0, \vep^*)$. Making use of this result and the observation that \bee \{ \mbf{n} \leq m \} & \supseteq  & \li \{ F_{Y_m}
(Y_m, \; \mscr{U} (Y_m, \vep) ) \leq \f{\de_m}{2}, \; G_{Y_m} (Y_m, \; \mscr{L} (Y_m, \vep) ) \leq \f{\de_m}{2}, \;
\mscr{L} (Y_m, \vep) \in \Se, \; \mscr{U} (Y_m, \vep) \in \Se \ri \}\\
 &  \supseteq  & \li \{  \f{ \sq{{ m }} [ \mscr{U}
(Y_m, \vep) - Y_m   ] }{  \sq{ \mcal{V} \li ( \mscr{U} (Y_m,
\vep)  \ri ) } }  \geq \Phi^{-1} \li ( 1 -  \f{\de_m}{2} + \f{\mcal{B} ( \mscr{U} (Y_m, \vep) ) }{\sq{m}} \ri ), \; \mscr{U} (Y_m, \vep) \in \Se \ri \}\\
&  & \bigcap \li \{ \f{ \sq{{ m }} [ Y_m - \mscr{L} (Y_m, \vep)  ]   }{ \sq{ \mcal{V} \li ( \mscr{L} (Y_m, \vep) \ri ) } } \geq \Phi^{-1} \li (
1 - \f{\de_m}{2} + \f{\mcal{B} ( \mscr{L} (Y_m, \vep) ) }{\sq{m}} \ri ), \; \mscr{L} (Y_m, \vep) \in \Se \ri \}, \eee we have that $\{ | Y_m -
\mu | < \ga \} \subseteq \{ \mbf{n} \leq m \}$ for all $m \geq (1 + \eta) N$ with $\vep \in (0, \vep^*)$.  The proof of the lemma is thus
completed.

\epf

\beL

\la{genhow89}

Assume that $\Pr \{ \lim_{\vep \downarrow 0} \f{ \mbf{n} } { \mcal{N} (\vep, \de,  \mu) }  = 1 \} = 1$ and \be \la{asump89}
 \lim_{\vep
\downarrow 0} \f{ \sum_{m \geq (1 + \eta) N} \Pr \{ \mbf{n} > m \} } { \mcal{N} (\vep, \de, \mu) } = 0 \qqu \tx{for any $\eta \in (0, 1)$}. \ee
Then, \[ \lim_{\vep \downarrow 0} \f{ \bb{E} [ \mbf{n} ] } { \mcal{N} (\vep, \de, \mu) } = 1.
\]
\eeL

\bpf For notational simplicity, let $N = \mcal{N} (\vep, \de, \mu)$ as before.  Let $\eta \in (0, 1)$.    As a consequence of the assumption
that $\Pr \{ \lim_{\vep \downarrow 0} \f{ \mbf{n} } { \mcal{N} (\vep, \de, \mu) } = 1 \} = 1$, we have \be
\lim_{\vep \downarrow 0} \Pr \{ (1 - \eta) N  \leq  \mbf{n} \leq (1 + \eta) N \} = 1.  \la{limeq6} \\
\ee  Noting that {\small \bee \bb{E} [ \mbf{n} ] = \sum_{m = 0}^\iy m \Pr \{ \mbf{n} = m \}  \geq \sum_{ (1 - \eta) N \leq m \leq (1 + \eta) N }
m \Pr \{ \mbf{n} = m \} \geq (1 - \eta) N  \sum_{ (1 - \eta) N \leq m \leq (1 + \eta) N   } \Pr \{ \mbf{n} = m \},  \eee} we have \be \la{ieq99}
\bb{E} [ \mbf{n} ]  \geq (1 - \eta) N \Pr \{ (1 - \eta) N \leq \mbf{n} \leq (1 + \eta) N \}. \ee Combining (\ref{limeq6}) and (\ref{ieq99})
yields
\[
\liminf_{\vep \downarrow 0} \f{ \bb{E} [ \mbf{n} ]  } { \mcal{N} (\vep, \de, \mu) } \geq (1 - \eta) \lim_{\vep \downarrow 0} \Pr \{ (1 - \eta) N
\leq \mbf{n} \leq (1 + \eta) N \} = 1 - \eta.
\]
Since the above inequality holds for any $\eta \in (0, 1)$, we have that \be \la{usegooda}
 \liminf_{\vep \downarrow 0} \f{ \bb{E} [ \mbf{n} ]  }
{ \mcal{N} (\vep, \de, \mu) } \geq 1. \ee
 On the other hand, using $\bb{E} [ \mbf{n} ] = \sum_{m = 0}^\iy \Pr \{ \mbf{n}
> m \}$, we can write {\small \bee \bb{E} [ \mbf{n} ] = \sum_{0 \leq m < (1 + \eta) N } \Pr \{ \mbf{n}
> m \}  + \sum_{m \geq (1 + \eta) N } \Pr \{ \mbf{n} > m \} \leq  \lc (1 + \eta) N  \rc + \sum_{m \geq (1 + \eta) N } \Pr \{ \mbf{n} > m \}.
\eee}  Observing that {\small $\limsup_{\vep \downarrow 0}  \f{ \lc (1 + \eta) N  \rc  } { \mcal{N} (\vep, \de, \mu) } =  1 + \eta$} and making
use of the assumption (\ref{asump89}), we have  that \[ \limsup_{\vep \downarrow 0} \f{ \bb{E} [ \mbf{n} ] } { \mcal{N} (\vep, \de, \mu) } \leq
1 + \eta \]
 holds for any $\eta \in (0, 1)$, which implies \be \la{usegoodb} \limsup_{\vep \downarrow 0} \f{ \bb{E} [ \mbf{n} ] } { \mcal{N}
(\vep, \de, \mu) } \leq 1. \ee Finally, the lemma is established by combining (\ref{usegooda}) and (\ref{usegoodb}).

\epf

We are now in a position to prove (\ref{prop5}).  According to Lemma \ref{howgood966}, we have that \bee \Pr \{ \mbf{n} > m \} \leq \Pr \{ | Y_m
- \mu | \geq \ga \} \eee for all $m \geq (1 + \eta) N$ with sufficiently small $\vep
> 0$.
Using (\ref{mostvip}) of Lemma \ref{nonuniformBE}, we have \bee  \limsup_{\vep \downarrow 0} \f{ \sum_{m \geq (1 + \eta) N} \Pr \{ \mbf{n} > m
\} } { \mcal{N} (\vep, \de, \mu) }
 & \leq & \limsup_{\vep \downarrow 0} \f{ \sum_{m \geq (1 + \eta) N} \Pr \{ | Y_m - \mu | \geq \ga \} } { \mcal{N} (\vep, \de, \mu) }\\
& \leq & \limsup_{\vep \downarrow 0} \sum_{m \geq (1 + \eta) N}  \f{1}{ \mcal{N} (\vep, \de, \mu) } \li [ \exp \li ( - \f{m}{2} \f{ \ga^2 }{
\mcal{V} (\mu) } \ri ) + \f{2 C }{m^2} \f{ \mscr{W}(\mu) }{ \ga^3 } \ri ] = 0,  \eee which implies that (\ref{asump89}) holds.  Since $\Pr \{
\lim_{\vep \downarrow 0} \f{ \mbf{n} } { \mcal{N} (\vep, \de, \mu) } = 1 \} = 1$ and  (\ref{asump89}) is true, it follows from Lemma
\ref{genhow89} that $\lim_{\vep \downarrow 0} \f{ \bb{E} [ \mbf{n} ] } { \mcal{N} (\vep, \de, \mu) } = 1$.

\subsection{Proof of Inequality (\ref{VIP886})}

By the established property that $\Pr \{ \mbf{n} < \iy \} = 1$, we have \be \la{ineq868a}
 \Pr \li \{ \mu \notin \mbf{I} \ri \} \leq \sum_{n =
m_\vep}^\iy \Pr \li \{ \mu \notin \mbf{I}, \; \mbf{n} = n \ri \}. \ee Clearly, \be \la{ineq868b} \Pr \li \{ \mu \notin \mbf{I}, \; \mbf{n} = n
\ri \} \leq \Pr \li \{  \mbf{n} = n \ri \} \qu \tx{for all $n \geq m_\vep$}. \ee We claim that \be \la{claim8889}
 \Pr \li \{ \mu \notin \mbf{I}, \; \mbf{n} = n \ri \} \leq \de_n \qu \tx{for all $n \geq m_\vep$}.
\ee To prove this claim, note that, as a consequence of the definitions of the stopping rule and the sequential random interval $\mbf{I}$,  \bel
\Pr \li \{ \mu \notin \mbf{I}, \; \mbf{n} = n \ri \} & \leq &  \Pr \li \{ \mu \geq \mscr{U} (Y_n,
\vep), \; \mbf{n} = n \ri \} + \Pr \li \{ \mu \leq \mscr{L} (Y_n, \vep), \; \mbf{n} = n \ri \} \nonumber \\
& \leq & \Pr \li \{ \mu \geq  \mscr{U} (Y_n, \vep), \; F_{Y_n} ( Y_n, \mscr{U} (Y_n, \vep) ) \leq \f{\de_n}{2} \ri \} \nonumber \\
&  &  + \Pr \li \{ \mu \leq \mscr{L} (Y_n, \vep), \; G_{Y_n} ( Y_n, \mscr{L} (Y_n, \vep) ) \leq \f{\de_n}{2} \ri \}  \la{ineq8896a} \eel for all
$n \geq m_\vep$. By the assumption that $\Pr \{ Y_n \leq z \}$ is non-increasing with respect to $\mu \in \Se$, we have that \bel  \Pr \li \{
\mu \geq \mscr{U} (Y_n, \vep), \; F_{Y_n} ( Y_n, \mscr{U} (Y_n, \vep) ) \leq  \f{\de_n}{2} \ri \} & \leq & \Pr \li \{ \mu \geq  \mscr{U} (Y_n,
\vep), \; F_{Y_n} ( Y_n, \mu ) \leq \f{\de_n}{2} \ri \} \nonumber\\
& \leq &  \Pr \li \{ F_{Y_n} ( Y_n, \mu ) \leq \f{\de_n}{2} \ri \} \nonumber \\
& \leq & \f{\de_n}{2} \la{ineq8896b} \eel for all $n \geq m_\vep$. Since $\Pr \{ Y_n \leq z \}$ is non-increasing with respect to $\mu \in \Se$,
it follows that $\Pr \{ Y_n \geq z \}$ is non-decreasing with respect to $\mu \in \Se$.  Hence, \bel  \Pr \li \{ \mu \leq \mscr{L} (Y_n, \vep),
\; G_{Y_n} ( Y_n, \mscr{L} (Y_n, \vep) ) \leq  \f{\de_n}{2} \ri \} & \leq & \Pr \li \{ \mu \leq  \mscr{L} (Y_n, \vep), \;
G_{Y_n} ( Y_n, \mu ) \leq \f{\de_n}{2} \ri \} \nonumber\\
& \leq &  \Pr \li \{ G_{Y_n} ( Y_n, \mu ) \leq \f{\de_n}{2} \ri \} \nonumber\\
& \leq & \f{\de_n}{2} \la{ineq8896c} \eel for all $n \geq m_\vep$. Combining (\ref{ineq8896a}), (\ref{ineq8896b}) and (\ref{ineq8896c}) yields
(\ref{claim8889}).  The claim is thus proved.  Finally, the inequality (\ref{VIP886}) is established by virtue of (\ref{ineq868a}),
(\ref{ineq868b}) and (\ref{claim8889}).

\section{Proof of Theorem \ref{ThmCHB} } \la{ThmCHB_app}

\subsection{Proof of Property (I)}

To establish (\ref{prop1}), we can make use of a similar method as that of the counterpart of Theorem \ref{ThmCDF} and Lemma \ref{asn889222} as
follows.

\beL \la{asn889222} There exist a real number $\ga > 0$ and an integer $m > 0$ such that \be \la{makeuse222}
 \Pr \{ \mbf{n} > n \}
\leq \Pr \{ | Y_n - \mu | \geq \ga \} \qqu \tx{for all $n > m$}. \ee \eeL

\bpf From Lemma \ref{remove asp}, we have  that there exists a number $\ga > 0$ such that $(\mu - \ga, \mu + \ga) \subseteq \Se$ and that \[
\sup_{\se \in (\mu - \ga, \mu + \ga) } \mscr{M} (\se, \mscr{U} ( \se, \vep) ) < 0.
\]
Since $\lim_{n \to \iy} \Phi \li (  \sq{ 2 \ln \f{\de_n}{2}  } \ri ) = 1 - \f{\de}{2}$, we have that $\f{ \ln \f{\de_n}{2} }{n} \to 0$ as $n \to
\iy$.  It follows that there exists $m
> 0$ such that
\[
\f{ \ln \f{\de_n}{2} }{n} > \sup_{\se \in (\mu - \ga, \mu + \ga) } \mscr{M} (\se, \mscr{U} ( \se, \vep) ) \qu \tx{for $n > m$}.
\]
This implies that
\[
\li \{ \mscr{M} ( Y_n, \mscr{U} ( Y_n, \vep) )  > \f{ \ln \f{\de_n}{2} }{n}, \; | Y_n - \mu | < \ga \ri \} = \emptyset \qu \tx{for all $n
> m$}.
\]
It follows that {\small \bee  \li \{ \mscr{M} ( Y_n, \mscr{U} ( Y_n, \vep) )  > \f{ \ln \f{\de_n}{2} }{n} \ri \} & \subseteq & \li \{ \mscr{M} (
Y_n, \mscr{U} ( Y_n, \vep) )  > \f{ \ln \f{\de_n}{2} }{n}, \; | Y_n - \mu | < \ga \ri
\} \cup \{ | Y_n - \mu | \geq \ga \}\\
&  = & \{ | Y_n - \mu | \geq \ga \} \eee} for all $n > m$.  So, we have shown that there exist a real number $\ga
> 0$ and an integer $m > 0$ such that $\li \{ \mscr{M} ( Y_n, \mscr{U} ( Y_n, \vep) )  > \f{ \ln \f{\de_n}{2} }{n} \ri \}  \subseteq \{ | Y_n - \mu |
\geq \ga \}$ for all $n > m$.  In a similar manner, we can show that there exist a real number $\ga
> 0$ and an integer $m > 0$ such that
\[
\li \{ \mscr{M} ( Y_n, \mscr{L} ( Y_n, \vep) )  > \f{ \ln \f{\de_n}{2} }{n} \ri \} \subseteq  \{ | Y_n - \mu | \geq \ga \} \qu \tx{for all $n >
m$}.
\]
Therefore, there exist a real number $\ga > 0$ and an integer $m > 0$ such that {\small \[ \li \{ \mscr{M} ( Y_n, \mscr{U} (Y_n, \vep) ) \leq
\f{\ln \f{\de_n}{2} }{n}, \; \mscr{M} ( Y_n, \mscr{L} (Y_n, \vep) ) \leq \f{ \ln \f{\de_n}{2} }{n} \ri \} \supseteq   \{ | Y_n - \mu | < \ga \}
\]}
for all $n > m$.  Making use of this result and the definition of the stopping rule, we have that {\small $\{ \mbf{n} > n \} \subseteq \{ | Y_n
- \mu | \geq \ga \}$} for all $n > m$.  This completes the proof of the lemma. \epf

\subsection{Proof of Property (II)}

To show (\ref{prop3}), define \[ \bs{L}_{\mbf{n}} = \mscr{L} ( Y_\mbf{n}, \vep ), \qqu \bs{L}_{\mbf{n} - 1} = \mscr{L} ( Y_{\mbf{n} - 1}, \vep
), \qqu \bs{U}_{\mbf{n}} = \mscr{U} ( Y_\mbf{n}, \vep ), \qqu \bs{U}_{\mbf{n} - 1} = \mscr{U} ( Y_{\mbf{n} - 1}, \vep )
\]
and {\small \bee &  & \mscr{E}_1 = \li \{  \mscr{M} (Y_\mbf{n}, \bs{U}_{\mbf{n}} ) \leq \f{ \ln (\f{\de_\mbf{n}}{2}) }{\mbf{n}}, \; \mscr{M}
(Y_\mbf{n}, \bs{L}_{\mbf{n}} )
\leq \f{ \ln (\f{\de_\mbf{n}}{2}) }{\mbf{n}} \ri \},\\
&  & \mscr{E}_2 = \{ \mbf{n} = m_\vep \} \bigcup  \li \{ \mscr{M} (Y_{\mbf{n} - 1}, \bs{U}_{\mbf{n} - 1} ) > \f{ \ln (\f{\de_{\mbf{n} - 1} }{2})
}{\mbf{n} - 1}, \; \mbf{n} > m_\vep  \ri \}   \bigcup \li \{ \mscr{M} (Y_{\mbf{n} -
1}, \bs{L}_{\mbf{n} - 1} ) > \f{ \ln (\f{\de_{\mbf{n}-1}}{2}) }{\mbf{n} - 1}, \;  \mbf{n} > m_\vep  \ri \},\\
&  & \mscr{E}_3 = \li \{  \lim_{\vep \downarrow 0} \mbf{n} = \iy, \; \lim_{\vep \downarrow 0} Y_{\mbf{n}} = \mu, \; \lim_{\vep \downarrow 0}
Y_{\mbf{n} - 1} = \mu \ri \},\\
&  & \mscr{E}_4 = \mscr{E}_1 \cap \mscr{E}_2 \cap \mscr{E}_3, \\
&  & \mscr{E}_5 = \li \{ \lim_{\vep \downarrow 0} \f{ \mbf{n} }{ \mcal{N} (\vep, \de, \mu) } = 1 \ri \}.  \eee} Recalling the assumption that
$m_\vep \to \iy$ as $\vep \downarrow 0$, we have that $\{ \lim_{\vep \downarrow 0} \mbf{n} = \iy \}$ is a sure event.  By the strong law of
large numbers, $\mscr{E}_3$ is an almost sure event.  By the definition of the stopping rule, we have that $\mscr{E}_1 \cap \mscr{E}_2$ is an
almost sure event.  It follows that $\mscr{E}_4$ is an almost sure event.  Hence, to show (\ref{prop3}), it suffices to show $\mscr{E}_4
\subseteq \mscr{E}_5$.  For this purpose, we let $\om \in \mscr{E}_4$ and expect to show $\om \in \mscr{E}_5$. For simplicity of notations, let
\[
n = \mbf{n} (\om), \qu y_n = Y_\mbf{n} (\om), \qu y_{n-1} = Y_{\mbf{n} - 1} (\om) \] \[ L_{n-1} = \mscr{L} ( y_{n-1}, \vep), \qu U_{n-1}=
\mscr{U} ( y_{n-1}, \vep).
\]
We restrict $\vep > 0$ to be sufficiently small so that \[ L_n \in \Se, \qu U_n\in \Se, \qu L_{n - 1} \in \Se, \qu U_{n-1} \in \Se.
\]
As a consequence of $\om \in \mscr{E}_1$,  \bel & & n \geq \f{ \ln (\f{\de_n}{2}) }{\mscr{M} (y_n, U_n )}, \la{imp18} \\
&  & n \geq \f{ \ln (\f{\de_n}{2}) }{\mscr{M} (y_n, L_n
 )}.  \la{imp2}
\eel Note that (\ref{imp18}) can be written as \be \la{comb8963a}
  \f{ n }{ \mcal{N} (\vep, \de, \mu) }
 \geq \f{ 2 \ln (\f{\de_n}{2}) } { \mcal{Z}^2  } \f{  [ U_n- y_n ]^2 }
{ 2 \mcal{V} ( U_n) \mscr{M} (y_n, U_n ) }  \f{ \mcal{V} ( U_n)  } { [ U_n- y_n ]^2 }  \f{ [\ka (\mu) \vep]^2 }{ \mcal{V} (\mu) }.  \ee Invoking
Lemma \ref{uniformcon}, we have \be \la{comb8963b}
 \lim_{\vep \downarrow 0} \f{  [ U_n-
y_n ]^2 } { 2 \mcal{V} ( U_n) \mscr{M} (y_n, U_n ) }  = - 1. \ee Combining (\ref{comb8963a}) and (\ref{comb8963b}) yields \be \la{imp3}
\liminf_{\vep \downarrow 0} \f{ n }{ \mcal{N} (\vep, \de, \mu) } \geq 1. \ee Similarly, using Lemma \ref{uniformcon},  we can deduce
(\ref{imp3}) from (\ref{imp2}).  As a consequence of (\ref{imp3}) and the assumption that $\lim_{\vep \downarrow 0} \vep^2 m_\vep = 0$, it must
be true that $n \neq m_\vep$ for small enough $\vep > 0$.  In view of this fact and $\om \in \mscr{E}_2$, we have that either \bee & & n - 1 <
\f{ \ln (\f{\de_{n-1}}{2}) }{\mscr{M} (y_{n-1}, U_{n-1}
 )} \qu \tx{or} \qu  n - 1 < \f{ \ln (\f{\de_{n-1}}{2}) }{\mscr{M} (y_{n-1}, L_{n - 1}
 )}
\eee must be true for small enough $\vep > 0$.  Hence, either
\[
\f{ n }{ \mcal{N} (\vep, \de, \mu) }  < \f{ n } { n - 1 } \f{ 2 \ln (\f{\de_{n-1}}{2}) } { \mcal{Z}^2  } \f{  [ U_{n-1} - y_{n-1} ]^2 } { 2
\mcal{V} ( U_{n-1} ) \mscr{M} (y_{n-1}, U_{n-1}  ) }  \f{ \mcal{V} ( U_{n-1} ) } { [ U_{n-1} - y_{n-1} ]^2 } \f{ [\ka (\mu) \vep]^2 }{ \mcal{V}
(\mu) }
\]
or
\[
\f{ n }{ \mcal{N} (\vep, \de, \mu) }  < \f{ n } { n - 1 } \f{ 2 \ln (\f{\de_{n-1}}{2}) } { \mcal{Z}^2  } \f{  [ L_{n-1} - y_{n-1} ]^2 } { 2
\mcal{V} ( L_{n-1} ) \mscr{M} (y_{n-1}, L_{n-1}  ) }  \f{ \mcal{V} ( L_{n-1} ) } { [ L_{n-1} - y_{n-1} ]^2 } \f{ [\ka (\mu) \vep]^2 }{ \mcal{V}
(\mu) }
\]
must be true for small enough $\vep > 0$. Making use of Lemma \ref{uniformcon} and these two inequalities, we have \be \la{imp88} \limsup_{\vep
\downarrow 0} \f{ n }{ \mcal{N} (\vep, \de, \mu) } \leq 1. \ee Combining (\ref{imp3}) and (\ref{imp88}) yields \[ \lim_{\vep \downarrow 0} \f{ n
}{ \mcal{N} (\vep, \de, \mu) } = 1.
\]
Thus, we have shown that $\mscr{E}_4 \subseteq \mscr{E}_5$, which implies that (\ref{prop3}) is true.

\subsection{Proof of Property (III)}

To establish (\ref{prop4}), we can make use of a similar method as that of the counterpart of Theorem \ref{ThmCDF} and Lemma \ref{covch} as
follows.

\beL

\la{covch}

\bel &  &  \Pr \li \{
 \lim_{\vep \downarrow 0} \f{ \sq{ \mbf{n} } [ Y_\mbf{n} - \bs{\mcal{U}}_\mbf{n} ]
} { \sq{ \mcal{V} (\mu) } } = - \mcal{Z}   \ri \}  = 1, \la{consistency1ch}\\
&  & \Pr \li \{
 \lim_{\vep \downarrow 0} \f{ \sq{ \mbf{n} } [  \bs{\mcal{L}}_\mbf{n} - Y_\mbf{n} ]
} { \sq{ \mcal{V} (\mu) } } = - \mcal{Z}   \ri \}  = 1, \la{consistency2ch} \eel where $\bs{\mcal{L}}_{\mbf{n}} = \mcal{L} ( Y_\mbf{n}, \vep)$
and $ \bs{\mcal{U}}_{\mbf{n}} = \mcal{U} ( Y_\mbf{n}, \vep)$.

\eeL

\bpf Let $\mscr{E}_i, \; i = 1, \cd, 5$ be defined as before in the proof of property (II).  Define \[ \mscr{E}_6 = \li \{
 \lim_{\vep \downarrow 0} \f{ \sq{ \mbf{n} } [ Y_\mbf{n} - \bs{\mcal{U}}_\mbf{n} ]
} { \sq{ \mcal{V} (\mu) } } = - \mcal{Z}   \ri \}.
\]
To show  (\ref{consistency1ch}), it suffices to show $\mscr{E}_4 \subseteq \mscr{E}_6$.   For this purpose, we let $\om \in \mscr{E}_4$ and
attempt to show $\om \in \mscr{E}_6$.  For simplicity of notations, let $n = \mbf{n} (\om), \;  y_n = Y_\mbf{n} (\om), \; y_{n-1} = Y_{\mbf{n} -
1} (\om)$,
\[
L_n = \mscr{L} ( y_n, \vep), \qu U_n = \mscr{U} ( y_n, \vep), \qu \mcal{L}_n = \mcal{L} ( y_n, \vep), \qu \mcal{U}_n = \mcal{U} ( y_n, \vep),
\]
\[
L_{n-1} = \mscr{L} ( y_{n-1}, \vep), \qu U_{n-1} = \mscr{U} ( y_{n-1}, \vep), \qu \mcal{L}_{n-1} = \mcal{L} ( y_{n-1}, \vep), \qu \mcal{U}_{n-1}
= \mcal{U} ( y_{n-1}, \vep).
\]
We restrict $\vep > 0$ to be sufficiently small so that \[ L_n \in \Se, \qu U_n\in \Se, \qu L_{n - 1} \in \Se, \qu U_{n-1} \in \Se.
\]
As a consequence of $\om \in \mscr{E}_1$,  we have $\mscr{M} (y_n, U_n) \leq \f{\ln \f{\de_n}{2}}{n}$, which can be written as
\[
\f{ \sq{n} [ \mcal{U}_n - y_n ] }{  \sq{ \mcal{V} (\mu) } } \geq \mcal{Z} \sq{ \f{2 \ln (\f{\de_n}{2}) }{\mcal{Z}^2} \f{ (U_n- y_n)^2  } { 2
\mcal{V} (U_n) \mscr{M} (y_n, U_n)} \f{ \mcal{V} (U_n)   } {  \mcal{V} (\mu)  } }  \f{ \mcal{U}_n - y_n }{ U_n- y_n }.
\]
By virtue of this inequality and Lemma \ref{uniformcon}, we have
\[
\liminf_{\vep \downarrow 0} \f{ \sq{n} [ \mcal{U}_n - y_n ] }{  \sq{ \mcal{V} (\mu) } } \geq \mcal{Z}.
\]
On the other hand, since  $n \neq m_\vep$ for small enough $\vep > 0$ and $\om \in \mscr{E}_2$, it follows that either $\mscr{M} (y_{n-1},
U_{n-1} ) > \f{\ln \f{\de_{n-1}}{2}}{n-1}$ or $\mscr{M} (y_{n-1}, L_{n-1} ) > \f{\ln \f{\de_{n-1}}{2}}{n-1}$ must be true for small enough $\vep
> 0$.  This implies that either
\[
\f{ \sq{n} [ \mcal{U}_n - y_{ n  } ] }{  \sq{ \mcal{V} (\mu) }  }  < \mcal{Z} \sq{ \f{ n } { n - 1  }  \f{2 \ln (\f{\de_{n-1}}{2}) }{\mcal{Z}^2}
\f{ (U_{n-1} - y_{n-1} )^2  } {2 \mcal{V} (U_{n-1} ) \mscr{M} (y_{n - 1}, U_{n-1} )} \f{ \mcal{V} (U_{n-1})   } {  \mcal{V} (\mu)  } }  \f{
\mcal{U}_n - y_n }{ U_{n-1} - y_{n-1} } \qu \tx{or}
\]
\[
\f{ \sq{n} [ \mcal{U}_n - y_{ n  } ] }{  \sq{ \mcal{V} (\mu) }  }  < \mcal{Z} \sq{ \f{ n } { n - 1  }  \f{2 \ln (\f{\de_{n-1}}{2}) }{\mcal{Z}^2}
\f{ (L_{n-1} - y_{n-1} )^2  } {2 \mcal{V} (L_{n-1} ) \mscr{M} (y_{n - 1}, L_{n-1} )} \f{ \mcal{V} (L_{n-1})   } {  \mcal{V} (\mu)  } }  \f{
\mcal{U}_n - y_n }{ y_{n-1} - L_{n-1} }
\]
must be true for small enough $\vep > 0$.  Making use of these inequalities and Lemma \ref{uniformcon}, we have
\[
\limsup_{\vep \downarrow 0} \f{ \sq{n} [ \mcal{U}_n - y_n ] }{  \sq{ \mcal{V} (\mu) } } \leq \mcal{Z}.
\]
Hence, {\small $\lim_{\vep \downarrow 0} \f{ \sq{n} [ \mcal{U}_n - y_n ] }{  \sq{ \mcal{V} (\mu) } } = \mcal{Z}$}.  This shows $\om \in
\mscr{E}_6$ and thus $\mscr{E}_4 \subseteq \mscr{E}_6$. It follows that $\Pr \{  \mscr{E}_6  \} = 1$, which implies (\ref{consistency1ch}).
Similarly, we can show (\ref{consistency2ch}).  This completes the proof of the lemma.

\epf

\subsection{Proof of Property (IV)}

To establish (\ref{prop5}), we can make use of a similar method as that of the counterpart of Theorem \ref{ThmCDF} and  Lemma \ref{asnch899} as
follows.

\beL \la{asnch899}  For any $\eta \in (0, 1)$, there exist $\ga > 0$ and $\vep^* > 0$ such that $\{ \mbf{n} > m \} \subseteq  \{ | Y_m - \mu |
\geq \ga \}  $ for all $m \geq (1 + \eta) N$ with  $\vep \in (0, \vep^*)$, where $N = \mcal{N} (\vep, \de, \mu)$.

\eeL

\bpf

Clearly,  for $\eta \in (0, 1)$,  there exists $\vsi \in (0, 1)$ such that $( 1 - \vsi )^4 ( 1 + \eta )
> 1$.  Since $\Phi( \sq{ 2 \ln \f{2}{\de_m}  }  ) \to 1 - \f{\de}{2}$ as $m \to \iy$, there exists $\vep_0 > 0$ such that
\be \la{see88} \ln \f{2}{\de_m} < \f{\mcal{Z}^2}{2(1 - \vsi)} \ee  for all $m \geq ( 1 + \eta) N$ with $\vep \in (0, \vep_0)$. As a consequence
of $\ka(\mu) > 0$ and the continuity assumption associated with (\ref{uniformasp}), there exist $\ga_1 > 0$ and $\vep_1 \in (0, \vep_0)$ such
that \be \la{use89a}
 [ \mscr{U} ( z, \vep) - z ]^2 > (1 - \vsi)^2 [\ka(\mu) \vep]^2, \qqu [ \mscr{L} ( z, \vep) - z ]^2 > (1 - \vsi)^2
[\ka(\mu) \vep]^2 \ee for all $z \in (\mu - \ga_1, \mu + \ga_1)$ and
$\vep \in (0, \vep_1)$.  Since $\mcal{V} \li ( \mscr{U} (z, \vep)
\ri )$ and $\mcal{V} \li ( \mscr{L} (z, \vep)  \ri )$ are continuous
functions of $z$ and $\vep$ which tend to $\mcal{V}(\mu) > 0$ as
$(z, \vep) \to (\mu, 0)$, there exist $\ga_2 \in (0, \ga_1)$ and
$\vep_2 \in (0, \vep_1)$ such that \[ \mscr{L} (z, \vep) \in \Se,
\qqu \mscr{U} (z, \vep) \in \Se
\]
and that \be \la{use89b}
 \mcal{V} \li ( \mscr{U} (z, \vep)
\ri ) <  \li ( 1 - \vsi \ri )^4 ( 1 + \eta ) \mcal{V} (\mu), \qqu  \mcal{V} \li ( \mscr{L} (z, \vep) \ri ) < \li ( 1 - \vsi \ri )^4 ( 1 + \eta )
\mcal{V} (\mu) \ee for all $z \in (\mu - \ga_2, \mu + \ga_2)$ and $\vep \in (0, \vep_2)$.  Using (\ref{use89a}), (\ref{use89b}) and the fact
that $\mcal{V} (\mu) = \f{ N [ \ka(\mu) \vep ]^2 } { \mcal{Z}^2  }$, we have
 \[ \mcal{V} \li ( \mscr{U} (z, \vep)  \ri ) <  \li ( 1 - \vsi \ri )^4 ( 1 + \eta ) \f{ N [
\ka(\mu) \vep ]^2 } { \mcal{Z}^2 } < \li ( 1 - \vsi \ri )^2 ( 1 + \eta ) \f{ N [ \mscr{U} ( z, \vep) - z ]^2 } { \mcal{Z}^2 },
\]
\[
  \mcal{V} \li ( \mscr{L} (z, \vep) \ri ) < \li ( 1 - \vsi \ri )^4 ( 1 + \eta ) \f{ N [
\ka(\mu) \vep ]^2 } { \mcal{Z}^2  } < \li ( 1 - \vsi \ri )^2 ( 1 + \eta ) \f{ N [ \mscr{L} ( z, \vep) - z ]^2 } { \mcal{Z}^2 }
\]
for all $z \in (\mu - \ga_2, \mu + \ga_2)$ and $\vep \in (0, \vep_2)$.   It follows that  \be \la{see8}
 \f{ \mcal{Z}^2 }{2}  <   ( 1 + \eta ) N
\f{ \li ( 1 - \vsi \ri )^2 [ \mscr{U} ( z, \vep) - z ]^2 } { 2  \mcal{V} \li ( \mscr{U} (z, \vep) \ri ) }, \qqu  \f{ \mcal{Z}^2 }{2}  <   ( 1 +
\eta ) N \f{ \li ( 1 - \vsi \ri )^2 [ \mscr{L} ( z, \vep) - z ]^2 } { 2  \mcal{V} \li ( \mscr{L} (z, \vep) \ri ) } \ee for all $z \in (\mu -
\ga_2, \mu + \ga_2)$ and $\vep \in (0, \vep_2)$.   According to Lemma \ref{uniformcon}, there exist $\ga \in (0, \ga_2)$ and $\vep^* \in (0,
\vep_2)$ such that \be \la{see8b}
 \mscr{M} (z, \mscr{U} ( z, \vep) ) <  \f{ - (1 - \vsi) [ \mscr{U} ( z, \vep) - z ]^2 } { 2 \mcal{V} ( \mscr{U}
( z, \vep) ) }, \qqu \mscr{M} (z, \mscr{L} ( z, \vep) ) <  \f{ - (1 - \vsi) [\mscr{L} ( z, \vep) - z] ^2 } { 2 \mcal{V} ( \mscr{L} ( z, \vep) )
} \ee for $z \in (\mu - \ga, \mu + \ga)$ and $\vep \in (0, \vep^*)$.   Combining (\ref{see8}) and (\ref{see8b}) yields  \be \la{see8c}
\f{\mcal{Z}^2}{2 (1 - \vsi) } < ( 1 + \eta ) N [ - \mscr{M} (z, \mscr{U} ( z, \vep) ) ] , \qqu \f{\mcal{Z}^2}{2(1 - \vsi)} <   ( 1 + \eta ) N [
- \mscr{M} (z, \mscr{L} ( z, \vep) ) ] \ee for all $z \in (\mu - \ga, \mu + \ga)$ and $\vep \in (0, \vep^*)$.  It follows from (\ref{see88}) and
(\ref{see8c}) that
\[ \ln \f{2}{\de_m}  <   m [ - \mscr{M} (z, \mscr{U} ( z, \vep) ) ] , \qqu \ln \f{2}{\de_m}  <   m  [ - \mscr{M} (z, \mscr{L} ( z, \vep) ) ]
\]
for all $z \in (\mu - \ga, \mu + \ga)$ and $m \geq ( 1 + \eta ) N$
with $\vep \in (0, \vep^*)$.  This implies that {\small \bee &  & \{ | Y_m - \mu | < \ga \}\\
&  & \subseteq \li \{ \mscr{M} (Y_m, \mscr{U} ( Y_m, \vep) )  \leq \f{ \ln (\f{\de_m}{2}) } {m}, \; \mscr{M} (Y_m, \mscr{L} ( Y_m, \vep) ) \leq
\f{ \ln (\f{\de_m}{2}) } {m}, \; \mscr{L} ( Y_m, \vep) \in \Se, \; \mscr{U} ( Y_m, \vep)
\in \Se  \ri \}\\
&  & \subseteq \{ \mbf{n} \leq m \} \eee} for all $m \geq ( 1 + \eta
) N$ with $\vep \in (0, \vep^*)$.  The proof of the lemma is thus
completed.

\epf

\subsection{Proof of Inequality (\ref{VIP886b})}

By the established property that $\Pr \{ \mbf{n} < \iy \} = 1$, we have \be \la{ineq868a889}
 \Pr \li \{ \mu \notin \mbf{I} \ri \} \leq \sum_{n =
m_\vep}^\iy \Pr \li \{ \mu \notin \mbf{I}, \; \mbf{n} = n \ri \}. \ee Clearly, \be \la{ineq868b889} \Pr \li \{ \mu \notin \mbf{I}, \; \mbf{n} =
n \ri \} \leq \Pr \li \{  \mbf{n} = n \ri \} \qu \tx{for all $n \geq m_\vep$}. \ee We claim that \be \la{claim8889889}
 \Pr \li \{ \mu \notin \mbf{I}, \; \mbf{n} = n \ri \} \leq \de_n \qu \tx{for all $n \geq m_\vep$}.
\ee To prove this claim, note that, as a consequence of the definitions of the stopping rule and the sequential random interval $\mbf{I}$,  \bel
\Pr \li \{ \mu \notin \mbf{I}, \; \mbf{n} = n \ri \} & \leq &  \Pr \li \{ \mu \geq \mscr{U} (Y_n,
\vep), \; \mbf{n} = n \ri \} + \Pr \li \{ \mu \leq \mscr{L} (Y_n, \vep), \; \mbf{n} = n \ri \} \nonumber \\
& \leq & \Pr \li \{ \mu \geq  \mscr{U} (Y_n, \vep), \; \mscr{M} ( Y_n, \mscr{U} (Y_n, \vep) ) \leq \f{1}{n} \ln \f{\de_n}{2} \ri \} \nonumber \\
&  &  + \Pr \li \{ \mu \leq \mscr{L} (Y_n, \vep), \; \mscr{M} ( Y_n, \mscr{L} (Y_n, \vep) ) \leq \f{1}{n} \ln \f{\de_n}{2} \ri \}  \nonumber\\
& \leq & \Pr \li \{ \mu \geq  \mscr{U} (Y_n, \vep), \; F_{Y_n} ( Y_n, \mscr{U} (Y_n, \vep) ) \leq \f{\de_n}{2} \ri \} \nonumber \\
&  &  + \Pr \li \{ \mu \leq \mscr{L} (Y_n, \vep), \; G_{Y_n} ( Y_n, \mscr{L} (Y_n, \vep) ) \leq \f{\de_n}{2} \ri \}  \la{ineq8896a889}
 \eel for
all $n \geq m_\vep$.

By the assumption that $\Pr \{ Y_n \leq z \}$ is non-increasing with respect to $\mu \in \Se$, we have that \bel  \Pr \li \{ \mu \geq \mscr{U}
(Y_n, \vep), \; F_{Y_n} ( Y_n, \mscr{U} (Y_n, \vep) ) \leq  \f{\de_n}{2} \ri \} & \leq & \Pr \li \{ \mu \geq  \mscr{U}
(Y_n, \vep), \; F_{Y_n} ( Y_n, \mu ) \leq \f{\de_n}{2} \ri \} \nonumber\\
& \leq &  \Pr \li \{ F_{Y_n} ( Y_n, \mu ) \leq \f{\de_n}{2} \ri \} \nonumber \\
& \leq & \f{\de_n}{2} \la{ineq8896b889} \eel for all $n \geq m_\vep$. Since $\Pr \{ Y_n \leq z \}$ is non-increasing with respect to $\mu \in
\Se$, it follows that $\Pr \{ Y_n \geq z \}$ is non-decreasing with respect to $\mu \in \Se$.  Hence, \bel  \Pr \li \{ \mu \leq \mscr{L} (Y_n,
\vep), \; G_{Y_n} ( Y_n, \mscr{L} (Y_n, \vep) ) \leq  \f{\de_n}{2} \ri \} & \leq & \Pr \li \{ \mu \leq  \mscr{L} (Y_n, \vep), \;
G_{Y_n} ( Y_n, \mu ) \leq \f{\de_n}{2} \ri \} \nonumber\\
& \leq &  \Pr \li \{ G_{Y_n} ( Y_n, \mu ) \leq \f{\de_n}{2} \ri \} \nonumber\\
& \leq & \f{\de_n}{2} \la{ineq8896c889} \eel for all $n \geq m_\vep$. Combining (\ref{ineq8896a889}), (\ref{ineq8896b889}) and
(\ref{ineq8896c889}) yields (\ref{claim8889889}).  The claim is thus proved.  Finally, the inequality (\ref{VIP886b}) is established by virtue
of (\ref{ineq868a889}), (\ref{ineq868b889}) and (\ref{claim8889889}).

\sect{Proof of Theorem \ref{ThmNPL}} \la{ThmNPL_app}

\subsection{Proof of Property (I)}

To establish (\ref{prop1}), we can make use of a similar method as that of the counterpart of Theorem \ref{ThmCDF} and Lemma \ref{finitech899}
as follows.

\beL \la{finitech899}

For $\mu \in \Se$, there exist a number $\ga > 0$ and an integer $m > 0$ such that $\{ \mbf{n} > n \} \subseteq \{ | Y_n - \mu | \geq \ga \}$
for all $n > m$. \eeL

\bpf  To show the lemma, we need to prove two claims as follows.

(i): For $\mu \in \Se$, there exist a number $\ga > 0$ and  a positive integer $m > 0$  such that
\[
\li \{  \mcal{V} \li ( Y_n + \ro [\mscr{U} (Y_n, \vep) - Y_n ] \ri )  > \f{ n [ Y_n - \mscr{U} (Y_n, \vep) ]^2 } { \ln \f{1}{\de_n}  } \ri \}
\subseteq \{  | Y_n - \mu | \geq \ga   \} \qu \tx{for all $n > m$}.
\]

(ii): For $\mu \in \Se$, there exist a number $\ga > 0$ and  a positive integer $m > 0$  such that
\[
\li \{  \mcal{V} \li ( Y_n + \ro [\mscr{L} (Y_n, \vep) - Y_n ] \ri )  > \f{ n [ Y_n - \mscr{L} (Y_n, \vep) ]^2 } { \ln \f{1}{\de_n}  } \ri \}
\subseteq \{  | Y_n - \mu | \geq \ga   \} \qu \tx{for all $n > m$}.
\]

To show the first claim, note that there exists $\ga > 0$ such that  $\mu - \ga, \; \mu + \ga \in \Se$ and that $\mscr{U} (\se, \vep) - \se \geq
\f{ \mscr{U} (\mu, \vep) - \mu }{2} > 0$ for all $\se \in (\mu - \ga, \mu + \ga) \subseteq \Se$.  Define $S = \{ \se \in (\mu - \ga, \mu + \ga):
\se + \ro [ \mscr{U} (\se, \vep) - \se ] \in \Se \}$. If $S = \emptyset$, then \[ \li \{ \mcal{V} \li ( Y_n + \ro [\mscr{U} (Y_n, \vep) - Y_n ]
\ri ) > \f{n}{\ln \f{1}{\de_n}} [ Y_n - \mscr{U} (Y_n, \vep) ]^2, \; | Y_n - \mu | < \ga \ri \} = \emptyset
\]
 for $n \in \bb{N}$.  In the other case that $S \neq \emptyset$, we have that
there exists a constant $D < \iy$ such that
\[
\sup_{ z \in S } \mcal{V}(z + \ro [\mscr{U} (z, \vep) - z ] ) < D.
\]
It follows that there exists an integer $m > 0$ such that {\small \bee &  &  \li \{ \mcal{V} \li ( Y_n + \ro [\mscr{U} (Y_n,
\vep) - Y_n ] \ri ) > \f{n}{\ln \f{1}{\de_n}} [ Y_n - \mscr{U} (Y_n, \vep) ]^2, \; | Y_n - \mu | < \ga  \ri \}\\
&  & = \li \{ n [ Y_n - \mscr{U} (Y_n, \vep) ]^2 < \mcal{V} \li ( Y_n + \ro [\mscr{U} (Y_n, \vep) - Y_n ] \ri ) \ln \f{1}{\de_n}, \; Y_n + \ro
[\mscr{U} (Y_n, \vep) - Y_n ]  \in \Se, \; |
Y_n - \mu | < \ga  \ri \}\\
&  & \subseteq \li \{ n [ Y_n - \mscr{U} (Y_n, \vep) ]^2 < D \ln \f{1}{\de_n}, \;  |
Y_n - \mu | < \ga  \ri \}\\
&  & \subseteq \li \{  \mscr{U} (Y_n, \vep) - Y_n < \f{ \mscr{U} (\mu, \vep) - \mu  }{2}, \; | Y_n - \mu | < \ga \ri \} = \emptyset \eee} for
all $n > m$.   It follows that {\small \bee &  & \li \{ \mcal{V} \li ( Y_n + \ro [\mscr{U} (Y_n, \vep) -
Y_n ] \ri ) > \f{n}{\ln \f{1}{\de_n}} [ Y_n - \mscr{U} (Y_n, \vep) ]^2  \ri \}\\
 & \subseteq & \li \{ n [
Y_n - \mscr{U} (Y_n, \vep) ]^2 < \mcal{V} \li ( Y_n + \ro [\mscr{U} (Y_n, \vep) - Y_n ] \ri ) \ln \f{1}{\de_n}, \; | Y_n - \mu | < \ga \ri \}
\cup \{ |
Y_n - \mu | \geq \ga \}\\
& = & \{ | Y_n - \mu | \geq \ga \}  \eee} for all $n > m$. This establishes the first claim. In a similar manner, we can show the second claim.
Finally, making use of the two established claims and the definition of the stopping rule completes the proof of the lemma.  \epf

\subsection{Proof of Property (II)}

To show (\ref{prop3}), define \[ \bs{L}_{\mbf{n}} = \mscr{L} ( Y_\mbf{n}, \vep ), \qqu \bs{L}_{\mbf{n} - 1} = \mscr{L} ( Y_{\mbf{n} - 1}, \vep
), \qqu \bs{U}_{\mbf{n}} = \mscr{U} ( Y_\mbf{n}, \vep ), \qqu \bs{U}_{\mbf{n} - 1} = \mscr{U} ( Y_{\mbf{n} - 1}, \vep )
\]
and {\small \bee &  & \mscr{E}_1 = \li \{  \mcal{V} \li ( Y_\mbf{n} + \ro [ \bs{U}_{\mbf{n}} - Y_\mbf{n} ] \ri ) \leq \f{ \mbf{n} [ Y_\mbf{n} -
\bs{U}_{\mbf{n}} ]^2  }{ \ln \f{1}{ \de_{\mbf{n}} }  }, \;  \mcal{V} \li ( Y_\mbf{n} + \ro [ \bs{L}_{\mbf{n}} - Y_\mbf{n} ] \ri )  \leq \f{
\mbf{n} [ Y_\mbf{n} -
\bs{L}_{\mbf{n}} ]^2  }{ \ln \f{1}{ \de_{\mbf{n}} }  } \ri \}, \\
&  & \mscr{E}_2 =  \li \{  \mcal{V} \li ( Y_{\mbf{n} - 1} + \ro [ \bs{U}_{\mbf{n}-1} - Y_{\mbf{n}-1} ] \ri ) > \f{ (\mbf{n}-1) [ Y_{\mbf{n}-1} -
\bs{U}_{\mbf{n}-1} ]^2  }{ \ln \f{1}{ \de_{\mbf{n}-1} }  }  \ri \}\\
&  & \qqu \; \bigcup  \li \{  \mcal{V} \li ( Y_{\mbf{n} - 1} + \ro [ \bs{L}_{\mbf{n}-1} - Y_{\mbf{n}-1} ] \ri ) > \f{ (\mbf{n}-1) [
Y_{\mbf{n}-1} - \bs{L}_{\mbf{n}-1} ]^2  }{ \ln \f{1}{ \de_{\mbf{n}-1} }  }  \ri \},\\
&  & \mscr{E}_3 =   [  \mscr{E}_2 \cap \{ \mbf{n} > m_\vep  \} ] \cup \{ \mbf{n} = m_\vep  \},\\
&  & \mscr{E}_4 = \li \{ \lim_{\vep \downarrow 0} \mbf{n} = \iy, \;  \lim_{\vep \downarrow 0} Y_{\mbf{n}}
= \mu, \; \lim_{\vep \downarrow 0} Y_{\mbf{n} - 1} = \mu \ri \},\\
&  &  \mscr{E}_5 = \mscr{E}_1 \cap \mscr{E}_3 \cap \mscr{E}_4,\\
&  & \mscr{E}_6 = \li \{ \lim_{\vep \downarrow 0} \f{ \mbf{n} } { \mcal{N} (\vep, \de, \mu) } = 1 \ri \}.  \eee} Recalling the assumption that
$m_\vep \to \iy$ almost surely as $\vep \downarrow 0$, we have that $\{ \lim_{\vep \downarrow 0} \mbf{n} = \iy \}$ is a sure event. By the
strong law of large numbers, $\mscr{E}_4$ is an almost sure event. By the definition of the stopping rule, $\mscr{E}_1 \cap \mscr{E}_3$ is an
almost sure event. Hence, $ \mscr{E}_5$ is an almost sure event. To show (\ref{prop3}), it suffices to show that $\mscr{E}_5 \subseteq
\mscr{E}_6$. For this purpose, we let $\om \in \mscr{E}_5$ and expect to show $\om \in \mscr{E}_6$.  For simplicity of notations, let $n =
\mbf{n} (\om), \;  y_n = Y_\mbf{n} (\om), \; y_{n-1} = Y_{\mbf{n} - 1} (\om)$,
\[
L_n = \mscr{L} ( y_n, \vep), \qu U_n = \mscr{U} ( y_n, \vep), \qu L_{n-1} = \mscr{L} ( y_{n-1}, \vep), \qu U_{n-1} = \mscr{U} ( y_{n-1}, \vep).
\]
We restrict $\vep > 0$ to be sufficiently small such that
\[ y_n + \ro [U_n- y_n ] \in \Se, \qqu y_n + \ro [L_n - y_n ] \in \Se.  \]  As a consequence of $\om \in \mscr{E}_1$,
\bee &  & \f{n} { \mcal{N} (\vep, \de, \mu) }  \geq  \f{ \ln \f{1}{\de_n} } { \mcal{Z}^2 }  \f{ \mcal{V} \li ( y_n + \ro
[U_n- y_n ] \ri ) } { \mcal{V} (\mu) }  \li [ \f{  \ka (\mu) \vep  }{ U_n- y_n } \ri ]^2,\\
&  & \f{n} { \mcal{N} (\vep, \de, \mu) }  \geq \f{ \ln \f{1}{\de_n} } { \mcal{Z}^2 } \f{ \mcal{V} \li ( y_n + \ro [L_n - y_n ] \ri ) } {
\mcal{V} (\mu) } \li [ \f{  \ka (\mu) \vep  }{ y_n - L_n } \ri ]^2. \eee These two inequalities imply that $\liminf_{\vep \downarrow 0} \f{n} {
\mcal{N} (\vep, \de, \mu) } \geq 1$.  As a consequence of this inequality and the assumption that $\lim_{\vep \downarrow 0} \vep^2 m_\vep = 0$,
it must be true that $n \neq m_\vep$ for small enough $\vep > 0$.  Making use of this observation and noting that $\om \in \mscr{E}_3$, we have
that either
\[
\f{n} { \mcal{N} (\vep, \de, \mu) }  < \f{ \ln \f{1}{\de_{n-1} } } { \mcal{Z}^2 }  \f{ n } { n - 1  } \f{ \mcal{V} \li ( y_{ n - 1 } + \ro
[U_{n-1} - y_{ n - 1 } ] \ri )   } { \mcal{V} (\mu)  } \li [ \f{ \ka (\mu) \vep }{ U_{n-1} - y_{ n - 1 } } \ri ]^2 \qu \tx{or}
\]
\[
\f{n} { \mcal{N} (\vep, \de, \mu) } < \f{ \ln \f{1}{\de_{n-1} } } { \mcal{Z}^2 } \f{ n } { n - 1  } \f{ \mcal{V} \li ( y_{ n - 1 } + \ro [L_{n -
1} - y_{ n - 1 } ] \ri )   } { \mcal{V} (\mu) } \li [ \f{ \ka (\mu) \vep }{ y_{ n - 1 } - L_{n - 1} } \ri ]^2
\]
must be true for small enough $\vep > 0$.  These two inequalities imply that $\limsup_{\vep \downarrow 0} \f{n} { \mcal{N} (\vep, \de, \mu) }
\leq 1$. Therefore, it must be true that $\lim_{\vep \downarrow 0} \f{n} { \mcal{N} (\vep, \de, \mu) } = 1$. This proves $\mscr{E}_5 \subseteq
\mscr{E}_6$. So, $\mscr{E}_6$ is an almost sure event, which implies (\ref{prop3}).

\subsection{Proof of Property (III)}

To establish (\ref{prop4}), we can make use of a similar method as that of the counterpart of Theorem \ref{ThmCDF} and Lemma \ref{covnpl} as
follows.

\beL

\la{covnpl}

\bel &  &  \Pr \li \{
 \lim_{\vep \downarrow 0} \f{ \sq{ \mbf{n} } [ Y_\mbf{n} - \bs{\mcal{U}}_\mbf{n} ]
} { \sq{ \mcal{V} (\mu) } } = - \mcal{Z}   \ri \}  = 1, \la{covnpl889a} \\
&  &  \Pr \li \{
 \lim_{\vep \downarrow 0} \f{ \sq{ \mbf{n} } [  \bs{\mcal{L}}_\mbf{n} - Y_\mbf{n} ]
} { \sq{ \mcal{V} (\mu) } } = - \mcal{Z}   \ri \}  = 1, \la{covnpl889b} \eel where $\bs{\mcal{L}}_{\mbf{n}} = \mcal{L} ( Y_\mbf{n}, \vep)$ and $
\bs{\mcal{U}}_{\mbf{n}} = \mcal{U} ( Y_\mbf{n}, \vep)$.

\eeL

\bpf  Let $\mscr{E}_i, \; i = 1, \cd, 6$ be defined as before in the proof of property (II).
  Define \[ \mscr{E}_7 = \li \{
 \lim_{\vep \downarrow 0} \f{ \sq{ \mbf{n} } [ Y_\mbf{n} - \bs{\mcal{U}}_\mbf{n} ]
} { \sq{ \mcal{V} (\mu) } } = - \mcal{Z}   \ri \}.
\]
To show the first equation, it suffices to show that $\mscr{E}_5 \subseteq \mscr{E}_7$. For this purpose, we let $\om \in \mscr{E}_5$ and
attempt to show that $\om \in \mscr{E}_7$.  For simplicity of notations, let $n = \mbf{n} (\om), \;  y_n = Y_\mbf{n} (\om), \; y_{n-1} =
Y_{\mbf{n} - 1} (\om)$,
\[
L_n = \mscr{L} ( y_n, \vep), \qu U_n = \mscr{U} ( y_n, \vep), \qu \mcal{L}_n = \mcal{L} ( y_n, \vep), \qu \mcal{U}_n = \mcal{U} ( y_n, \vep),
\]
\[
L_{n-1} = \mscr{L} ( y_{n-1}, \vep), \qu U_{n-1} = \mscr{U} ( y_{n-1}, \vep), \qu \mcal{L}_{n-1} = \mcal{L} ( y_{n-1}, \vep), \qu \mcal{U}_{n-1}
= \mcal{U} ( y_{n-1}, \vep).
\]
We restrict $\vep > 0$ to be sufficiently small such that
\[ y_n + \ro [U_n- y_n ] \in \Se, \qqu
y_n + \ro [L_n - y_n ] \in \Se.  \] As a consequence of $\om \in \mscr{E}_1$, we have
\[
\mcal{V} \li ( y_n + \ro [U_n- y_n ] \ri )  \leq \f{n}{ \ln \f{1}{\de_n} } [ U_n - y_n ]^2,
\]
which can be written as
\[
\f{ \sq{n} [ \mcal{U}_n - y_n ] }{  \sq{ \mcal{V} (\mu) } } \geq  \sq{ \f{ \mcal{V} \li ( y_n + \ro [U_n- y_n ] \ri )   } { \mcal{V} (\mu) } }
\f{ \mcal{U}_n - y_n  }{ U_n- y_n } \sq{ \ln \f{1}{\de_n} }.
\]
This implies that \be \la{use369a}
 \liminf_{\vep \downarrow 0} \f{ \sq{n} [ \mcal{U}_n - y_n ] }{  \sq{ \mcal{V} (\mu) } } \geq \mcal{Z}.
\ee On the other hand, note that $n \neq m_\vep$ for small enough $\vep > 0$. As a consequence of this observation and $\om \in \mscr{E}_3$, we
have that either
\[ \mcal{V} \li ( y_{n-1} + \ro [U_{n-1} - y_{n-1} ] \ri )  > \f{n-1}{ \ln \f{1}{\de_{n-1}} } [ U_{n-1} -
y_{n-1} ]^2
\]
or
\[
\mcal{V} \li ( y_{n-1} + \ro [L_{n-1} - y_{n-1} ] \ri )  > \f{n-1}{ \ln \f{1}{\de_{n-1}} } [ L_{n-1} - y_{n-1} ]^2
\]
must be true for small enough $\vep > 0$.  This implies that either
\[
\f{ \sq{n} [ \mcal{U}_n - y_{ n  } ] }{  \sq{ \mcal{V} (\mu) }  }  < \sq{ \f{ n } { n - 1 } \f{ \mcal{V} \li ( y_{ n - 1 } + \ro [U_{n-1} - y_{
n - 1 } ] \ri )    } { \mcal{V} (\mu)    } } \f{ U_n- y_{ n } }{ U_{n-1} - y_{ n - 1 } } \f{ \mcal{U}_n - y_n }{ U_n- y_n } \sq{ \ln
\f{1}{\de_{n-1}} } \qu \tx{or}
\]
\[
\f{ \sq{n} [ \mcal{U}_n - y_{ n  } ] }{  \sq{ \mcal{V} (\mu) }  }  < \sq{ \f{ n } { n - 1 } \f{ \mcal{V} \li ( y_{ n - 1 } + \ro [L_{n-1} - y_{
n - 1 } ] \ri )  } {  \mcal{V} (\mu) } } \f{ U_n- y_{ n } }{ y_{ n - 1 } -  L_{n-1} } \f{ \mcal{U}_n - y_n }{ U_n- y_n } \sq{ \ln
\f{1}{\de_{n-1}} }
\]
must be true for small enough $\vep > 0$.  Making use of these inequalities, we have \be \la{use369b} \limsup_{\vep \downarrow 0} \f{ \sq{n} [
\mcal{U}_n - y_n ] }{  \sq{ \mcal{V} (\mu) } } \leq \mcal{Z}. \ee Combining (\ref{use369a}) and (\ref{use369b}) yields $\lim_{\vep \downarrow 0}
\f{ \sq{n} [ \mcal{U}_n - y_n ] }{  \sq{ \mcal{V} (\mu) } } = \mcal{Z}$.   This shows $\om \in \mscr{E}_7$ and thus $\mscr{E}_5 \subseteq
\mscr{E}_7$. It follows that $\Pr \{  \mscr{E}_7  \} = 1$, which implies (\ref{covnpl889a}). Similarly, we can show (\ref{covnpl889b}).  This
completes the proof of the lemma.

\epf

\subsection{Proof of Property (IV)}

To establish (\ref{prop5}), we can make use of a similar method as that of the counterpart of Theorem \ref{ThmCDF} and Lemma \ref{asnnpl899} as
follows.

\beL \la{asnnpl899} For any $\eta \in (0, 1)$, there exist $\ga > 0$ and $\vep^* > 0$ such that $\{ \mbf{n} > m \} \subseteq \{ | Y_m - \mu |
\geq \ga \}$ for all $m \geq (1 + \eta) N$ with $\vep \in (0, \vep^*)$, where $N = \mcal{N} (\vep, \de, \mu)$.

\eeL

\bpf

 Clearly, for $\eta \in (0, 1)$, there exists $\vsi \in (0, 1)$ such that $(1 - \vsi)^4 ( 1 + \eta ) > 1$.
 Since $\Phi \li ( \sq{ \ln \f{1}{\de_m} } \ri )  \to 1 - \f{\de}{2}$ as $m \to \iy$, there exists $\vep_0 > 0$ such that \be \la{convde} \sq{
\ln \f{1}{\de_m} }  < \f{\mcal{Z}}{1 - \vsi} \ee for all $m > (1 + \eta) N$ with $\vep \in (0, \vep_0)$. As a consequence of $\ka(\mu) > 0$ and
the continuity assumption associated with (\ref{uniformasp}), there exist $\ga_1 > 0$ and $\vep_1 \in (0, \vep_0)$ such that
 (\ref{cita}) and (\ref{citb})  hold for all $z \in (\mu - \ga_1,
\mu + \ga_1)$ and $\vep \in (0, \vep_1)$.  Since $\mcal{V} \li ( z + \ro [ \mscr{L} (z, \vep) - z]  \ri )$ and $\mcal{V} \li ( z + \ro [
\mscr{U} (z, \vep) - z]  \ri )$ are continuous functions of $z$ and $\vep$ which tend to $\mcal{V}(\mu)
> 0$ as $(z, \vep) \to (\mu, 0)$, there exist $\ga \in (0, \ga_1)$ and $\vep^* \in (0, \vep_1)$ such that
\[
z + \ro [ \mscr{U} (z, \vep) - z] \in \Se, \qqu z + \ro [ \mscr{L}
(z, \vep) - z] \in \Se
\]
and that \be \la{herewhere}
 \mcal{V} \li ( z + \ro [ \mscr{U} (z, \vep) - z]  \ri ) <  (1 - \vsi)^4 ( 1 + \eta ) \mcal{V} (\mu), \qqu  \mcal{V} \li ( z +
\ro [ \mscr{L} (z, \vep) - z] \ri ) < (1 - \vsi)^4 ( 1 + \eta ) \mcal{V} (\mu) \ee for all $z \in (\mu - \ga, \mu + \ga)$ and $\vep \in (0,
\vep^*)$.  It follows from (\ref{herewhere}) and the definition of $N$ that \be \la{nowhere}
 \f{ { (1 + \eta) N } [ ( 1 - \vsi) \ka( \mu ) \vep ]^2 }{ \mcal{V} \li ( z + \ro [ \mscr{U} (z, \vep) - z] \ri ) } > \li ( \f{\mcal{Z}}{1 - \vsi} \ri )^2, \qqu
  \f{ (1 + \eta) N [ ( 1 - \vsi) \ka( \mu )  \vep ]^2 }{  \mcal{V} \li ( z + \ro [ \mscr{L} (z, \vep) - z] \ri ) } >
 \li ( \f{\mcal{Z}}{1 - \vsi} \ri )^2
\ee for all $z \in (\mu - \ga, \mu + \ga)$ and $\vep \in (0, \vep^*)$.  Making use of (\ref{convde}) and (\ref{nowhere}), we have \be
\la{nowhere88}
 \f{  m [( 1 - \vsi) \ka( \mu ) \vep]^2 }{ \mcal{V} \li ( z + \ro [ \mscr{U} (z, \vep) - z] \ri ) } > \ln \f{1}{\de_m}, \qqu
  \f{ m [( 1 - \vsi) \ka( \mu )  \vep]^2 }{ \mcal{V} \li ( z + \ro [ \mscr{L} (z, \vep) - z] \ri ) } > \ln \f{1}{\de_m}
\ee for all $z \in (\mu - \ga, \mu + \ga)$ and $m \geq (1 + \eta) N$ with $\vep \in (0, \vep^*)$.  Combining (\ref{cita}), (\ref{citb}) and
(\ref{nowhere88}) yields
\[
 \f{ m  [\mscr{U} (z, \vep) - z]^2  }{ \mcal{V} \li ( z + \ro [ \mscr{U} (z, \vep) - z] \ri ) } > \ln \f{1}{\de_m}, \qqu
  \f{ m [z - \mscr{L} (z, \vep) ]^2 }{  \mcal{V} \li ( z + \ro [ \mscr{L} (z, \vep) - z] \ri ) } > \ln \f{1}{\de_m}
\]
for all $z \in (\mu - \ga, \mu + \ga)$ and $m \geq (1 + \eta) N$ with $\vep \in (0, \vep^*)$.   It follows that {\small \bee  \{ | Y_m - \mu | <
\ga \}  &   \subseteq & \li \{ \f{ m [ \mscr{U} (Y_m, \vep) - Y_m   ]^2 }{ \mcal{V} \li ( Y_m + \ro [ \mscr{U} (Y_m, \vep) - Y_m ] \ri )  } >
\ln \f{1}{\de_m}, \; \; Y_m + \ro
[ \mscr{U} (Y_m, \vep) - Y_m ] \in \Se \ri \}\\
&  & \bigcap \li \{ \f{ m [ Y_m - \mscr{L} (Y_m, \vep) ]^2 }{ \mcal{V} \li ( Y_m + \ro [ \mscr{L} (Y_m, \vep) - Y_m ]  \ri ) }
> \ln \f{1}{\de_m}, \; Y_m + \ro
[ \mscr{L} (Y_m, \vep) - Y_m ] \in \Se \ri \}\\
&  \subseteq &  \{ \mbf{n} \leq m \} \eee} for all $m \geq (1 + \eta) N$ with $\vep \in (0, \vep^*)$.  This implies that $\{ | Y_m - \mu | < \ga
\} \subseteq \{ \mbf{n} \leq m \}$ for all $m \geq (1 + \eta) N$ with $\vep \in (0, \vep^*)$.  The proof of the lemma is thus completed.

\epf

\sect{Proof of Theorem \ref{Seqgenfull}}  \la{Seqgenfull_app}

\subsection{Proof of Property (I)}

To show (\ref{prop1}), we need the following preliminary result.

\beL \la{finitefull}

If the random variable $X$ has mean $\mu$ and variance $\nu$, then there exist a number $\ga > 0$ and an integer $m
> 0$ such that $\{ \mbf{n} > n \} \subseteq \{ | Y_n - \mu | \geq \ga \; \tx{or} \; | V_{n} - \nu | \geq \ga  \}$ for all $n > m$.

\eeL

\bpf

   Note that there exists $\ga > 0$ such that
 \[
\mscr{U}( \se, \vep  ) - \se > \f{  \mscr{U}( \mu, \vep  ) - \mu }{2} > 0
 \]
 for all $\se \in (\mu - \ga, \mu + \ga) \subseteq \Se$.  Since $\lim_{n \to \iy} \Phi \li ( \sq{ \ln \f{1}{\de_n}  }  \ri ) = 1 - \f{\de}{2}$,
 it must be true that $\f{ \ln \f{1}{\de_n} }{n} \to 0$ as $n \to \iy$.  Therefore, there exists an integer $m > 0$ such that
\bee &  &  \li \{ V_{n} + \f{\ro}{n} > \f{ n [ \mscr{U}( Y_n, \vep  ) - Y_n ]^2}{  \ln \f{1}{\de_n} }, \;  | Y_n - \mu |
< \ga, \; | V_{n} - \nu | < \ga \ri \}\\
&  & \subseteq  \li \{ V_{n} + \f{\ro}{n} > \f{ n [ \f{  \mscr{U}( \mu, \vep  ) - \mu }{2} ]^2}{  \ln \f{1}{\de_n} }, \;  | Y_n - \mu | < \ga,
\; | V_{n} - \nu | < \ga \ri \} = \emptyset \eee for all $n > m$.  It follows that \be \la{cite9886}  \li \{ V_{n} + \f{\ro}{n} > \f{ n [
\mscr{U}( Y_n, \vep  ) - Y_n ]^2}{  \ln \f{1}{\de_n} } \ri \} \subseteq  \li \{ | Y_n - \mu | \geq \ga \; \tx{or} \; | V_{n} - \nu | \geq \ga
\ri \} \qu \tx{for all $n > m$}. \ee Thus, we have shown that there exist a number $\ga > 0$ and an integer $m > 0$ such that (\ref{cite9886})
holds.  In a similar manner, we can show that there exist a number $\ga > 0$ and an integer $m > 0$ such that \[ \li \{ V_{n} + \f{\ro}{n} > \f{
n [ \mscr{L}( Y_n, \vep  ) - Y_n ]^2}{  \ln \f{1}{\de_n} } \ri \} \subseteq  \li \{ | Y_n - \mu | \geq \ga \; \tx{or} \; | V_{n} - \nu | \geq
\ga \ri \} \qu \tx{for all $n > m$}. \] Finally, making use of these results and the definition of the stopping rule completes the proof of the
lemma.

\epf

We are now in a position to prove (\ref{prop1}). According to Lemma \ref{finitefull}, there exist an integer $m$ and a positive number $\ga > 0$
such that \be \la{makeusefull}
 \Pr \{ \mbf{n} > n \}
\leq \Pr \{ | Y_n - \mu | \geq \ga \} + \Pr \{ | V_n - \nu | \geq \ga  \} \qqu \tx{for all $n > m$}. \ee Making use of (\ref{makeusefull}),  the
weak law of large numbers, and Lemma \ref{vipnow},  we have $\lim_{n \to \iy} \Pr \{ \mbf{n} > n \} = 0$, which implies $\Pr \{ \mbf{n} < \iy \}
= 1$.  By virtue of (\ref{mostvip}) of Lemma \ref{nonuniformBE}, Lemma \ref{vipnow},  and (\ref{makeusefull}), we have \bee  \bb{E} [
\mbf{n} ] & \leq &  m + \sum_{n > m} [ \Pr \{ | Y_n - \mu | \geq \ga \} + \Pr \{ | V_n - \nu | \geq \ga  \} ]\\
&  \leq &  m + \sum_{n > m}  \li [ \exp \li ( - \f{n}{2} \f{ \ga^2 }{ \nu } \ri ) + \f{2 C }{n^2} \f{ \mscr{W} }{ \ga^3 } \ri ]\\
&   & + \sum_{n > m} \li [ \exp \li ( - \f{n}{4} \f{ \ga}{ \nu } \ri ) +  \exp \li ( - \f{n}{8} \f{ \ga^2 }{ \varpi } \ri ) + \f{4 C }{n^2
\ga^3} \li (  \sq{2} \mscr{W} \ga^{3\sh 2} +  4 \mscr{V} \ri ) \ri ]\\
& < & \iy. \eee

\bsk

\subsection{Proof of Property (II)}

To show (\ref{prop3}), define \[ \bs{L}_{\mbf{n}} = \mscr{L} ( Y_\mbf{n}, \vep ), \qqu \bs{L}_{\mbf{n} - 1} = \mscr{L} ( Y_{\mbf{n} - 1}, \vep
), \qqu \bs{U}_{\mbf{n}} = \mscr{U} ( Y_\mbf{n}, \vep ), \qqu \bs{U}_{\mbf{n} - 1} = \mscr{U} ( Y_{\mbf{n} - 1}, \vep )
\] and {\small \bee &  & \mscr{E}_1 = \li \{ V_\mbf{n} + \f{\ro}{\mbf{n}} \leq \f{ \mbf{n} [ Y_\mbf{n} - \bs{U}_{\mbf{n}}
]^2 }{ \ln \f{1}{\de_{\mbf{n}}} }, \;  V_\mbf{n} + \f{\ro}{\mbf{n}} \leq \f{ \mbf{n} [  \bs{L}_{\mbf{n}}  - Y_\mbf{n} ]^2 }{ \ln
\f{1}{\de_{\mbf{n}}} } \ri \}, \\
&  & \mscr{E}_2 = \li \{ V_{\mbf{n} - 1} + \f{\ro}{\mbf{n} - 1} > \f{ (\mbf{n} - 1) [ Y_{\mbf{n} - 1} -  \bs{U}_{\mbf{n} - 1} ]^2 }{ \ln
\f{1}{\de_{\mbf{n}-1}} } \ri \}  \bigcup \li \{ V_{\mbf{n}-1} + \f{\ro}{\mbf{n}-1} > \f{ (\mbf{n} - 1) [  \bs{L}_{\mbf{n}-1}  - Y_{\mbf{n}-1}
]^2 }{ \ln
\f{1}{\de_{\mbf{n}-1}} } \ri \},\\
&   & \mscr{E}_3 = \li \{ \lim_{\vep \downarrow 0} \mbf{n} = \iy, \; \lim_{\vep \downarrow 0} Y_{\mbf{n}} = \mu, \; \lim_{\vep \downarrow 0}
V_{\mbf{n}} = \nu, \; \lim_{\vep \downarrow 0} Y_{\mbf{n} - 1} = \mu, \; \lim_{\vep \downarrow
0} V_{\mbf{n} - 1} = \nu \ri \}, \\
&   & \mscr{E}_4 = \mscr{E}_1 \cap \mscr{E}_2 \cap \mscr{E}_3, \\
&  & \mscr{E}_5 = \li \{ \lim_{\vep \downarrow 0} \f{ \mbf{n} } { \mcal{N} (\vep, \de, \mu, \nu) } = 1 \ri \}. \eee} From the definition of the
stopping rule, it is obvious that  $\mbf{n} \to \iy$ almost surely as $\vep \downarrow 0$, which implies that $\{ \lim_{\vep \downarrow 0}
\mbf{n} = \iy \}$ is an almost sure event. By the strong law of large numbers, $\mscr{E}_3$ is an almost sure event.  By the definition of the
stopping rule, $\mscr{E}_1 \cap \mscr{E}_2$ is an almost sure event.  Hence, $\mscr{E}_4$ is an almost sure event. To show (\ref{prop3}), it
suffices to show $\mscr{E}_4 \subseteq \mscr{E}_5$.  For this purpose, we let $\om \in \mscr{E}_4$ and expect to show $\om \in \mscr{E}_5$.  For
simplicity of notations, let $n = \mbf{n} (\om), \;  y_n = Y_\mbf{n} (\om), \; y_{n-1} = Y_{\mbf{n} - 1} (\om), \; v_n = V_{\mbf{n}} (\om), \;
v_{n-1} = V_{\mbf{n} - 1} (\om)$,
\[
L_n = \mscr{L} ( y_n, \vep), \qu U_n = \mscr{U} ( y_n, \vep), \qu L_{n-1} = \mscr{L} ( y_{n-1}, \vep), \qu U_{n-1} = \mscr{U} ( y_{n-1}, \vep).
\]
For small $\vep > 0$, $y_n$ and $y_{n - 1}$ will be bounded in a neighborhood of $\mu$.  For small $\vep > 0$, $v_n$ and $v_{n - 1}$ will be
bounded in a neighborhood of $\nu$.   As a consequence of $\om \in \mscr{E}_1$, \bee &  &  \f{ n }{ \mcal{N} (\vep, \de, \mu, \nu) } \geq \f{
\ln \f{1}{\de_n} } { \mcal{Z}^2 }  \f{ v_{n} + \f{\ro}{n} } { \nu } \li [ \f{ \ka (\mu) \vep }{ U_n-
y_n } \ri ]^2, \\
&  & \f{ n }{ \mcal{N} (\vep, \de, \mu, \nu) }  \geq  \f{ \ln \f{1}{\de_n} } { \mcal{Z}^2 } \f{ v_{n} + \f{\ro}{n} } { \nu } \li [ \f{ \ka (\mu)
\vep }{ y_n - L_n } \ri ]^2. \eee These two inequalities imply that $\liminf_{\vep \downarrow 0} \f{ n }{ \mcal{N} (\vep, \de, \mu, \nu) } \geq
1$. On the other hand, as a consequence of $\om \in \mscr{E}_2$, we have that either
\[
\f{ n }{ \mcal{N} (\vep, \de, \mu, \nu) }   < \f{ \ln \f{1}{\de_{n-1}} } { \mcal{Z}^2 }   \f{ n } { n - 1 } \f{ v_{n - 1} + \f{\ro}{n - 1} } {
\nu } \li [ \f{ \ka (\mu) \vep }{ U_{n-1} - y_{ n - 1 } } \ri ]^2 \qu \tx{or}
\]
\[
\f{ n }{ \mcal{N} (\vep, \de, \mu, \nu) }   < \f{ \ln \f{1}{\de_{n-1}} } { \mcal{Z}^2 }  \f{ n } { n - 1 } \f{ v_{n - 1} + \f{\ro}{n - 1} } {
\nu } \li [ \f{ \ka (\mu) \vep }{ y_{ n - 1 } -  L_{n-1}  } \ri ]^2
\]
must be true. These two inequalities imply that $\limsup_{\vep \downarrow 0} \f{ n }{ \mcal{N} (\vep, \de, \mu, \nu) }  \leq 1$.  Hence,
$\lim_{\vep \downarrow 0} \f{ n }{ \mcal{N} (\vep, \de, \mu, \nu) }  = 1$.  This proves that $\mscr{E}_4 \subseteq \mscr{E}_5$. So, $\mscr{E}_5$
is an almost sure event, which implies (\ref{prop3}).

\subsection{Proof of Property (III)}

To establish (\ref{prop4}), we can make use of a similar method as that of the counterpart of Theorem \ref{ThmCDF} and Lemma \ref{covfull89} as
follows.

\beL \la{covfull89} \bel &  &  \Pr \li \{
 \lim_{\vep \downarrow 0} \f{ \sq{ \mbf{n} } ( Y_\mbf{n} - \bs{\mcal{U}}_\mbf{n} )
} { \sq{ \nu } } = - \mcal{Z}   \ri \}  = 1, \la{consistencyfull1}\\
&  &  \Pr \li \{  \lim_{\vep \downarrow 0} \f{ \sq{ \mbf{n} } (  \bs{\mcal{L}}_\mbf{n} - Y_\mbf{n} ) } { \sq{ \nu } } = - \mcal{Z} \ri \} = 1,
\la{consistencyfull2} \eel where $\bs{\mcal{L}}_{\mbf{n}} = \mcal{L} ( Y_\mbf{n}, \vep)$ and $ \bs{\mcal{U}}_{\mbf{n}} = \mcal{U} ( Y_\mbf{n},
\vep)$.

\eeL

\bpf  Let $\mscr{E}_i, \; i = 1, \cd, 5$ be defined as before in the proof of property (II).   Define \[ \mscr{E}_6 = \li \{
 \lim_{\vep \downarrow 0} \f{ \sq{ \mbf{n} } ( Y_\mbf{n} - \bs{\mcal{U}}_\mbf{n} )
} { \sq{ \nu } } = - \mcal{Z}   \ri \}.
\]
To show the first equation, it suffices to show $\mscr{E}_4 \subseteq \mscr{E}_6$.  For this purpose, we let $\om \in \mscr{E}_4$ and expect to
show $\om \in \mscr{E}_6$.  For simplicity of notations, let $n = \mbf{n} (\om), \;  y_n = Y_\mbf{n} (\om), \; y_{n-1} = Y_{\mbf{n} - 1} (\om),
\; v_n = V_{\mbf{n}} (\om), \; v_{n-1} = V_{\mbf{n} - 1} (\om)$,
\[
L_n = \mscr{L} ( y_n, \vep), \qu U_n = \mscr{U} ( y_n, \vep), \qu L_{n-1} = \mscr{L} ( y_{n-1}, \vep), \qu U_{n-1} = \mscr{U} ( y_{n-1}, \vep),
\]
\[
\mcal{L}_n = \mcal{L} ( y_n, \vep), \qu \mcal{U}_n = \mcal{U} ( y_n, \vep), \qu \mcal{L}_{n-1} = \mcal{L} ( y_{n-1}, \vep), \qu \mcal{U}_{n-1} =
\mcal{U} ( y_{n-1}, \vep).
\]
As a consequence of $\om \in \mscr{E}_1$,
\[
v_{n} + \f{\ro}{n} \leq \f{ n }{ \ln \f{1}{\de_n} } ( U_n - y_n )^2,
\]
which can be written as
\[
\f{ \sq{n} ( \mcal{U}_n - y_n ) }{ \sq{ \nu  }  } \geq \f{ \sq{ v_{n} + \f{\ro}{n} } } { \sq{ \nu  } } \f{ \mcal{U}_n - y_n  } { U_n- y_n  }
\sq{ \ln \f{1}{\de_n} }.  \] This implies that
\[
\liminf_{\vep \downarrow 0} \f{ \sq{n} ( \mcal{U}_n - y_n ) }{  \sq{ \mcal{V} (\mu) } } \geq \mcal{Z}.
\]
On the other hand, as a consequence of $\om \in \mscr{E}_2$,  we have that either
\[
v_{n-1} + \f{\ro}{n-1} > \f{ n -1}{ \ln \f{1}{\de_{n-1}} } ( U_{n-1} - y_{n-1} )^2 \qu \tx{or} \qu v_{n-1} + \f{\ro}{n-1} > \f{ n -1}{\ln
\f{1}{\de_{n-1}} } ( L_{n-1} - y_{n-1} )^2
\]
must be true.  Hence,  either
\[ \f{ \sq{n} ( \mcal{U}_n - y_{ n  } ) }{  \sq{ \nu  } } <  \sq{ \f{ n } { n - 1 } \f{ v_{n - 1} +
\f{\ro}{n - 1}  } { \nu } } \f{ U_n- y_{ n } }{ U_{n-1} - y_{ n - 1 } } \f{ \mcal{U}_n - y_n  } { U_n- y_n } \sq{ \ln \f{1}{\de_{n-1}} } \qu
\tx{or}
\]
\[
\f{ \sq{n} ( \mcal{U}_n - y_{ n  } ) }{  \sq{ \nu  } } < \sq{ \f{ n } { n - 1 } \f{ v_{n - 1} + \f{\ro}{n - 1} } { \nu
 } } \f{ U_n- y_{ n  }  }{ y_{ n - 1 } - L_{n-1}  } \f{ \mcal{U}_n - y_n  }
{  U_n- y_n  } \sq{ \ln \f{1}{\de_{n-1}} }
\]
must be true. Making use of these inequalities, we have
\[
\limsup_{\vep \downarrow 0} \f{ \sq{n} ( \mcal{U}_n - y_n ) }{  \sq{ \mcal{V} (\mu) } } \leq \mcal{Z}.
\]
It follows that
\[
\lim_{\vep \downarrow 0} \f{ \sq{n} ( \mcal{U}_n - y_n ) }{  \sq{ \mcal{V} (\mu) } } = \mcal{Z}.
\]
This shows $\om \in \mscr{E}_6$ and thus $\mscr{E}_4 \subseteq \mscr{E}_6$. It follows that $\Pr \{  \mscr{E}_6  \} = 1$, which implies
(\ref{consistencyfull1}).  Similarly, we can show (\ref{consistencyfull2}).  This completes the proof of the lemma.

\epf

\subsection{Proof of Property (IV)}

To establish (\ref{prop5}), we can make use of a similar method as that of the counterpart of Theorem \ref{ThmCDF} and Lemma \ref{asn8full89} as
follows.

\beL  \la{asn8full89} For any $\eta \in (0, 1)$, there exist $\ga > 0$ and $\vep^* > 0$ such that $\{ \mbf{n} > m \} \subseteq \{ | Y_m - \mu |
\geq \ga \; \tx{or}  \; | V_m - \nu | \geq \ga \}$ for all $m \geq (1 + \eta) N$ with $\vep \in (0, \vep^*)$, where $N = \mcal{N} (\vep, \de,
\mu, \nu)$.

\eeL

\bpf

Clearly,  for $\eta \in (0, 1)$, there exists $\vsi \in (0, 1)$ such that $(1 - \vsi)^4 ( 1 + \eta ) > 1$.  Since $\Phi \li ( \sq{ \ln
\f{1}{\de_m} } \ri )  \to 1 - \f{\de}{2}$ as $m \to \iy$, there exists $\vep_0 > 0$ such that \be \la{dede} \sq{ \ln \f{1}{\de_m} } <
\f{\mcal{Z}}{1 - \vsi} \ee for all $m \geq (1 + \eta) N$ with $\vep \in (0, \vep_0)$.  As a consequence of $\ka(\mu)
> 0$ and the continuity assumption associated with (\ref{uniformasp}), there exist $\ga_1 > 0$ and $\vep_1 \in (0, \vep_0)$ such that
 (\ref{cita}) and (\ref{citb})  hold for all $z \in (\mu - \ga_1,
\mu + \ga_1)$ and $\vep \in (0, \vep_1)$.  Note that there exist $\ga \in (0, \ga_1)$ and $\vep^* \in (0, \vep_1)$ such that \be \la{secon} \se
+ \f{\ro}{(1 + \eta) N}  <  (1 - \vsi)^4 ( 1 + \eta ) \nu \ee for all $\se \in (\nu - \ga, \nu + \ga)$ and $\vep \in (0, \vep^*)$. It follows
from (\ref{secon}) and the definition of $N$ that \be \la{hownow}
 \f{ (1 + \eta) N [ ( 1 - \vsi) \ka( \mu ) \vep ]^2 }{ \se +  \f{\ro}{(1 + \eta) N}    } > \li ( \f{\mcal{Z}}{1 - \vsi} \ri )^2
\ee for all $\se \in (\nu - \ga, \nu + \ga)$ and $\vep \in (0, \vep^*)$.  Making use of (\ref{dede}) and (\ref{hownow}), we have \be
\la{hownow88}
 \f{ m [ ( 1 - \vsi) \ka( \mu ) \vep ]^2 }{ \se +  \f{\ro}{m} } > \ln \f{1}{\de_m}
\ee for all $\se \in (\nu - \ga, \nu + \ga)$ and $m \geq (1 + \eta) N$ with $\vep \in (0, \vep^*)$. Combining (\ref{cita}), (\ref{citb}) and
(\ref{hownow88}) yields
\[
 \f{ m [\mscr{U} (z, \vep) - z]^2  }{ \se + \f{\ro}{m} } \geq \ln \f{1}{\de_m}, \qqu
  \f{ m [z - \mscr{L} (z, \vep) ]^2 }{ \se + \f{\ro}{m} } \geq \ln \f{1}{\de_m}
\]
for all $z \in (\mu - \ga, \mu + \ga), \; \se \in (\nu - \ga, \nu + \ga)$ and $m \geq (1 + \eta) N$ with $\vep \in (0, \vep^*)$.   It follows
that {\small \bee \{ | Y_m - \mu | < \ga, \; | V_m - \nu | < \ga \} & \subseteq & \li \{ V_m +  \f{\ro}{m} \leq \f{ m [ \mscr{U} (Y_m, \vep) -
Y_m   ]^2 }{ \ln \f{1}{\de_m} }, \; \; V_m +
 \f{\ro}{m} \leq  \f{m [ Y_m - \mscr{L} (Y_m, \vep) ]^2 }{ \ln \f{1}{\de_m} } \ri \}\\
& \subseteq  & \{ \mbf{n} \leq m \}  \eee} for all $m \geq (1 + \eta) N$ with $\vep \in (0, \vep^*)$.   The proof of the lemma is thus
completed.

\epf

\section{ General Preliminary Results for Multistage Sampling }

To investigate properties of the multistage schemes in this paper, we to establish some general preliminary results.

\beL

\la{pro689a} Let $\jm$ be defined by (\ref{defim89}).  Assume that {\small \[ \Pr \li \{ \liminf_{\vep \downarrow 0} \Lm_{\bs{l}} (\mu, \nu)
\geq 1 \ri \} = 1, \qqu \Pr \li \{ \tx{Both} \; \limsup_{\vep \downarrow 0} \Lm_{\bs{l} - 1} (\mu, \nu) \leq 1 \; \tx{and} \;  \liminf_{\vep
\downarrow 0} \bs{l} > 1 \; \tx{hold or} \; \lim_{\vep \downarrow 0} \bs{l}  = 1 \ri \} = 1.
\]}
Then, \be \la{lemvip689}
 \Pr \li \{  \jm - 1 \leq \liminf_{ \vep \downarrow 0} \bs{l} \leq \limsup_{ \vep \downarrow 0} \bs{l} \leq \jm \ri \} = 1. \ee
Moreover, \be \la{lemvip689b} \Pr \li \{  \lim_{ \vep \downarrow 0} \bs{l} = \jm \ri \} = 1 \qu \tx{provided that $\jm = 1$ or both $\Lm_{\jm-1}
(\mu, \nu) < 1$ and $\jm
> 1$ hold}.
\ee

\eeL

\bpf

Define
\[
\mscr{E} = \li \{ \liminf_{\vep \downarrow 0} \Lm_{\bs{l}} (\mu, \nu) \geq 1 \ri \} \bigcap \li \{ \tx{Both} \; \limsup_{\vep \downarrow 0}
\Lm_{\bs{l} - 1} (\mu, \nu) \leq 1 \; \tx{and} \; \liminf_{\vep \downarrow 0} \bs{l} > 1 \; \tx{hold or} \; \lim_{\vep \downarrow 0} \bs{l}  = 1
\ri \}
\]
and  $ \mscr{E}^\prime = \li \{  \jm - 1 \leq \liminf_{ \vep \downarrow 0} \bs{l} \leq \limsup_{ \vep \downarrow 0} \bs{l} \leq \jm \ri \}$.  To
show (\ref{lemvip689}), it suffices to show that $\mscr{E} \subseteq \mscr{E}^\prime$.  For this purpose, we let $\om \in \mscr{E}$ and attempt
to show $\om \in \mscr{E}^\prime$.   For simplicity of notations, let $\ell = \bs{l} (\om)$.  Since $\om \in \mscr{E}$, it must be true that \be
\la{3388966ieq}
 \liminf_{\vep \downarrow 0} \Lm_\ell (\mu, \nu) \geq 1.
\ee Moreover,  either $\lim_{\vep \downarrow 0} \ell = 1$ or
\[
\limsup_{\vep \downarrow 0} \Lm_{\ell - 1} (\mu, \nu) \leq 1, \qqu \liminf_{\vep \downarrow 0} \ell > 1 \]
 must be true.   On the other hand, it
follows from the definition of $\jm$ that either $\jm = 1$ or
\[
\Lm_{\jm - 1} (\mu, \nu) \leq 1 < \Lm_\jm (\mu, \nu), \qqu \jm > 1 \] must be true.  Therefore, one of the following $4$ cases must occur.

Case (a): \[ \lim_{\vep \downarrow 0} \ell = 1 = \jm.
\]

Case (b):
\[
\lim_{\vep \downarrow 0} \ell = 1 < \jm, \qqu \Lm_{\jm - 1} (\mu, \nu) \leq 1 < \Lm_\jm (\mu, \nu).
\]

Case (c):
\[
\limsup_{\vep \downarrow 0} \Lm_{\ell - 1} (\mu, \nu) \leq 1, \qqu \liminf_{\vep \downarrow 0} \ell > 1 = \jm.
\]

Case (d):
\[ \limsup_{\vep \downarrow 0} \Lm_{\ell - 1} (\mu, \nu) \leq 1, \qqu  \liminf_{\vep \downarrow 0} \ell > 1, \qqu \Lm_{\jm - 1} (\mu, \nu)
\leq 1 < \Lm_\jm (\mu, \nu), \qqu \jm > 1.
\]

Actually,  Case (c) is impossible, since  $\liminf_{\vep \downarrow 0} (\ell - 1) \geq \jm$ and it follows that $\limsup_{\vep \downarrow 0}
\Lm_{\ell - 1} (\mu, \nu) \geq \Lm_\jm (\mu, \nu) > 1$, which contradicts to $\limsup_{\vep \downarrow 0} \Lm_{\ell - 1} (\mu, \nu) \leq 1$.  By
virtue of (\ref{3388966ieq}), it can be verified that
\[ \limsup_{\vep \downarrow 0} \Lm_{\ell - 1} (\mu, \nu) \leq \Lm_{\jm - 1} (\mu, \nu) \leq \liminf_{\vep \downarrow 0} \Lm_\ell (\mu, \nu)
\]
holds in Cases (a), (b) and (d).  It follows from the monotonicity of $\Lm_\ell (\mu, \nu)$ with respect to $\ell$ that $\jm - 1 \leq \liminf_{
\vep \downarrow 0} \ell \leq \limsup_{ \vep \downarrow 0} \ell \leq \jm$.  This implies that $\om \in \mscr{E}^\prime$ and consequently
$\mscr{E} \subseteq \mscr{E}^\prime$.  Thus, $\mscr{E}^\prime$ is also an almost sure event and accordingly (\ref{lemvip689}) holds.

To show (\ref{lemvip689b}), define $\mscr{E}^{\prime \prime} = \li \{  \lim_{ \vep \downarrow 0} \bs{l} = \jm \ri \}$.  It suffices to show $\om
\in \mscr{E}^{\prime \prime}$ for $\om \in \mscr{E}$.  By the assumption that either $\jm = 1$ or both $\Lm_{\jm - 1} (\mu, \nu) < 1$ and $\jm
> 1$ hold, one of the following $4$ cases must happen.

Case (i): \[ \lim_{\vep \downarrow 0} \ell = 1 = \jm.
\]

Case (ii):
\[
\lim_{\vep \downarrow 0} \ell = 1 < \jm, \qqu \Lm_{\jm - 1} (\mu, \nu) < 1 < \Lm_\jm (\mu, \nu).
\]

Case (iii):
\[
\limsup_{\vep \downarrow 0} \Lm_{\ell - 1} (\mu, \nu) \leq 1, \qqu \liminf_{\vep \downarrow 0} \ell > 1 = \jm.
\]

Case (iv):
\[ \limsup_{\vep \downarrow 0} \Lm_{\ell - 1} (\mu, \nu) \leq 1, \qqu  \liminf_{\vep \downarrow 0} \ell > 1, \qqu \Lm_{\jm - 1} (\mu, \nu) < 1 < \Lm_\jm (\mu, \nu), \qqu \jm > 1.
\]

As discussed earlier, Case (iii) is impossible.  Moreover, Case (ii) is also impossible, since $\lim_{\vep \downarrow 0} \ell \leq \jm - 1$ and
it follows that  $\Lm_{\jm - 1} (\mu, \nu) \geq \liminf_{\vep \downarrow 0} \Lm_\ell (\mu, \nu) \geq 1$, which contradicts to $\Lm_{\jm - 1}
(\mu, \nu) < 1$. By virtue of (\ref{3388966ieq}), it can be verified that
\[ \limsup_{\vep \downarrow 0} \Lm_{\ell - 1} (\mu, \nu) \leq \Lm_{\jm - 1} (\mu, \nu) < \Lm_{\jm} (\mu, \nu) \leq \liminf_{\vep \downarrow 0} \Lm_\ell (\mu, \nu)
\]
holds in Cases (i) and (iv). It follows from the monotonicity of $\Lm_\ell (\mu, \nu)$ with respect to $\ell$ that $\lim_{ \vep \downarrow 0}
\ell = \jm$.  This implies that $\om \in \mscr{E}^{\prime \prime}$ and consequently $\mscr{E} \subseteq \mscr{E}^{\prime \prime}$. Thus,
$\mscr{E}^{\prime \prime}$ is also an almost sure event and accordingly (\ref{lemvip689b}) holds.

\epf

\beL

\la{lem20888}

For $\mu \in \Se$,  if there exist $\ga > 0, \; \vep^* > 0$ and an index $j \geq 1$ such that $\{ \bs{D}_\ell = 1  \} \subseteq \{ | Z_\ell -
\mu | \geq \ga \}$ holds for $\tau_\vep \leq \ell \leq j$ and  all $\vep \in (0, \vep^*)$, then $\lim_{\vep \downarrow 0} \sum_{\ell =
\tau_\vep}^j n_\ell \Pr \{ \bs{D}_\ell = 1 \} = 0$.

Similarly, for $\mu \in \Se$,  if there exist $\ga > 0, \; \vep^* > 0$ and an index $j \geq 1$ such that $\{ \bs{D}_\ell = 0  \} \subseteq \{ |
Z_\ell - \mu | \geq \ga \}$ holds for $\ell \geq j \geq \tau_\vep$ and all $\vep \in (0, \vep^*)$, then $\lim_{\vep \downarrow 0} \sum_{\ell =
j}^\iy n_\ell \Pr \{ \bs{D}_\ell = 0 \} = 0$.

 \eeL

\bpf First, we shall show that for arbitrary $\ga > 0$ and $C > 0$, \be \la{imp89630}
 \lim_{\vep \downarrow 0} \sum_{\ell = \tau_\vep}^\iy  \li [ n_\ell
\exp \li ( - \f{n_\ell}{2} \f{ \ga^2 }{  \nu } \ri ) + \f{2 C }{n_\ell} \f{ \mscr{W} }{ \ga^3 } \ri ] = 0. \ee By the assumption that
$n_{\tau_\vep} \to \iy$ as $\vep \downarrow 0$, we have that
\[
0 \leq \limsup_{\vep \downarrow 0} \sum_{\ell = \tau_\vep}^\iy  n_\ell \exp \li ( - \f{n_\ell}{2} \f{ \ga^2 }{  \nu } \ri ) \leq \limsup_{\vep
\downarrow 0} \sum_{k = n_{\tau_\vep}}^\iy  k \exp \li ( - \f{k}{2} \f{ \ga^2 }{  \nu } \ri ) = 0,
\]
which implies that \be \la{imp8963a} \lim_{\vep \downarrow 0} \sum_{\ell = \tau_\vep}^\iy  n_\ell \exp \li ( - \f{n_\ell}{2} \f{ \ga^2 }{ \nu }
\ri )  = 0. \ee By the assumption on $\{ n_\ell \}_{\ell \geq 1}$, we have that there exists a number $\eta > 0$ such that
\[
\inf_{\ell \geq 1} \f{n_{\ell + 1}}{n_\ell} > 1 + \eta,
\]
which implies that
\[
n_\ell \geq n_{\tau_\vep} (1 + \eta)^{\ell - \tau_\vep} \qu \tx{for all} \; \ell \geq \tau_\vep.
\]
Hence, \bee \limsup_{\vep \downarrow 0} \sum_{\ell = \tau_\vep}^\iy   \f{2 C }{n_\ell} \f{ \mscr{W} }{ \ga^3 } & \leq & \limsup_{\vep
\downarrow 0} \sum_{\ell = \tau_\vep}^\iy   \f{2 C }{n_{\tau_\vep} (1 + \eta)^{\ell - \tau_\vep} } \f{ \mscr{W} }{ \ga^3 }\\
& = & \limsup_{\vep \downarrow 0} \f{2 C }{n_{\tau_\vep} } \f{ \mscr{W} }{ \ga^3 }  \sum_{\ell = \tau_\vep}^\iy \f{1 }{(1 + \eta)^{\ell -
\tau_\vep} } \\
& \leq & \limsup_{\vep \downarrow 0} \f{2 C }{n_{\tau_\vep} } \f{ \mscr{W} }{ \ga^3 }  \sum_{\ell = 0}^\iy \f{1 }{(1 + \eta)^{\ell} }\\
& = & \li ( 1 + \f{1}{\eta} \ri )  \f{ 2 C \mscr{W} }{ \ga^3 } \limsup_{\vep \downarrow 0} \f{1 }{n_{\tau_\vep} } = 0,  \eee from which it
follows that \be \la{imp8963b} \lim_{\vep \downarrow 0} \sum_{\ell = \tau_\vep}^\iy   \f{2 C }{n_\ell} \f{ \mscr{W} }{ \ga^3 } = 0.  \ee
Combining (\ref{imp8963a}) and (\ref{imp8963b}) yields (\ref{imp89630}).

For $\mu \in \Se$,  if there exist $\ga > 0, \; \vep^* > 0$ and an index $j \geq 1$ such that $\{ \bs{D}_\ell = 1  \} \subseteq \{ | Z_\ell -
\mu | \geq \ga \}$ holds for $\tau_\vep \leq \ell \leq j$ and all $\vep \in (0, \vep^*)$, then,  using (\ref{mostvip}) of Lemma
\ref{nonuniformBE} and (\ref{imp89630}), we have, for some constant $C > 0$,   \bee \limsup_{\vep \downarrow 0} \sum_{\ell = \tau_\vep}^j n_\ell
\Pr \{ \bs{D}_\ell = 1 \} & \leq & \limsup_{\vep \downarrow 0} \sum_{\ell = \tau_\vep}^j n_\ell \Pr \{ | Z_\ell - \mu  | \geq
\ga \}\\
& \leq & \limsup_{\vep \downarrow 0} \sum_{\ell = \tau_\vep}^j n_\ell \li [ \exp \li ( - \f{n_\ell}{2} \f{ \ga^2 }{  \nu } \ri ) +
\f{2 C }{n_\ell^2} \f{ \mscr{W} }{ \ga^3 } \ri ]\\
& = &  \limsup_{\vep \downarrow 0} \sum_{\ell = \tau_\vep}^j  \li [ n_\ell \exp \li ( - \f{n_\ell}{2} \f{ \ga^2 }{  \nu } \ri ) +
\f{2 C }{n_\ell} \f{ \mscr{W} }{ \ga^3 } \ri ]\\
& \leq &  \lim_{\vep \downarrow 0} \sum_{\ell = \tau_\vep}^\iy  \li [ n_\ell \exp \li ( - \f{n_\ell}{2} \f{ \ga^2 }{  \nu } \ri ) + \f{2 C
}{n_\ell} \f{ \mscr{W} }{ \ga^3 } \ri ] = 0, \eee which implies $\lim_{\vep \downarrow 0} \sum_{\ell = \tau_\vep}^j n_\ell \Pr \{ \bs{D}_\ell =
1 \} = 0$.

For $\mu \in \Se$,  if there exist $\ga > 0, \; \vep^* > 0$ and an index $j \geq 1$ such that $\{ \bs{D}_\ell = 0  \} \subseteq \{ | Z_\ell -
\mu | \geq \ga \}$ holds for $\ell \geq j$ and all $\vep \in (0, \vep^*)$, then,  using (\ref{mostvip}) of Lemma \ref{nonuniformBE} and
(\ref{imp89630}), we have,  for some constant $C > 0$,   \bee \limsup_{\vep \downarrow 0} \sum_{\ell = j}^\iy n_\ell \Pr \{ \bs{D}_\ell = 0 \} &
\leq & \limsup_{\vep \downarrow 0} \sum_{\ell = j}^\iy n_\ell \Pr \{ | Z_\ell - \mu  | \geq \ga \}\\
& \leq & \limsup_{\vep \downarrow 0} \sum_{\ell = j}^\iy n_\ell \li [ \exp \li ( - \f{n_\ell}{2} \f{ \ga^2 }{  \nu } \ri ) +
\f{2 C }{n_\ell^2} \f{ \mscr{W} }{ \ga^3 } \ri ]\\
& = &  \limsup_{\vep \downarrow 0} \sum_{\ell = j}^\iy  \li [ n_\ell \exp \li ( - \f{n_\ell}{2} \f{ \ga^2 }{  \nu } \ri ) +
\f{2 C }{n_\ell} \f{ \mscr{W} }{ \ga^3 } \ri ]\\
& \leq &  \lim_{\vep \downarrow 0} \sum_{\ell = \tau_\vep}^\iy  \li [ n_\ell \exp \li ( - \f{n_\ell}{2} \f{ \ga^2 }{  \nu } \ri ) + \f{2 C
}{n_\ell} \f{ \mscr{W} }{ \ga^3 } \ri ] = 0, \eee which implies $\lim_{\vep \downarrow 0} \sum_{\ell = j}^\iy n_\ell \Pr \{ \bs{D}_\ell = 0 \} =
0$.  This completes the proof of the lemma.

\epf

\beL

\la{lem21888}

For any $\mu \in \Se$, if there exists an index $i \geq 1$ such that
\bel &  & \limsup_{\vep \downarrow 0} \sum_{\ell = \tau_\vep}^{i - 2} n_\ell \Pr \{ \bs{D}_\ell = 1 \} < \iy, \la{same869}\\
&  & \limsup_{\vep \downarrow 0} \sum_{\ell = i}^\iy n_\ell \Pr \{ \bs{D}_{\ell} = 0 \} < \iy, \la{same8} \eel then,  \bel & & \limsup_{\vep
\downarrow 0} \f{ \bb{E} [ \mbf{n} ] } { \mcal{N} (\vep, \de, \mu, \nu)  }  = \f{1} { \mcal{Z}^2 } \li [ \xi_{i-1}^2 + (\xi_i^2 - \xi_{i-1}^2)
\limsup_{\vep \downarrow 0} \Pr \{ \bs{l} = i \}  \ri ] \leq \li ( \f{ \xi_{i} }{ \mcal{Z} } \ri )^2, \la{use986a}\\
 &  & \liminf_{\vep \downarrow 0} \f{ \bb{E} [ \mbf{n} ] } { \mcal{N} (\vep, \de, \mu, \nu)  }  = \f{1} { \mcal{Z}^2 } \li [ \xi_{i-1}^2 + (\xi_i^2 - \xi_{i-1}^2)
\liminf_{\vep \downarrow 0} \Pr \{ \bs{l} = i \} \ri ]  \geq \li ( \f{ \xi_{i-1} }{ \mcal{Z} } \ri )^2. \la{use986b} \eel For any $\mu \in \Se$,
if there exists an index $i \geq 1$ such that (\ref{same8}) holds and that \be \la{same6c8} \limsup_{\vep \downarrow 0} \sum_{\ell =
\tau_\vep}^{i - 1} n_\ell \Pr \{ \bs{D}_\ell = 1 \} < \iy,  \ee then,
\[
\lim_{\vep \downarrow 0 } \f{ \bb{E} [ \mbf{n} ] } { \mcal{N} (\vep, \de, \mu, \nu) }  =  \li ( \f{ \xi_i  }{ \mcal{Z} } \ri )^2.
\]

\eeL

\bpf

First, consider $\mu$ such that (\ref{same869}) and (\ref{same8}) hold. By the definition of the sampling scheme, we have {\small \bee \bb{E} [
\mbf{n} ]  =  \sum_{\ell = \tau_\vep }^{i - 2} n_\ell \Pr \{ \bs{l} = \ell \} + \sum_{\ell = i + 1}^\iy n_\ell \Pr \{ \bs{l} = \ell \} + n_{i -
1} \Pr \{ \bs{l} = i - 1 \} + n_i  \Pr \{ \bs{l} = i \} \eee} and {\small $\bb{E} [ \mbf{n} ] \geq n_{i - 1} \Pr \{ \bs{l} = i - 1 \} + n_i
 \Pr \{ \bs{l} = i \} $}.  Note that \be \la{citelatera}
 \limsup_{\vep
\downarrow 0} \sum_{\ell = \tau_\vep}^{i - 2}  n_\ell \Pr \{ \bs{l} = \ell  \} \leq \limsup_{\vep \downarrow 0} \sum_{\ell = \tau_\vep}^{i - 2}
n_\ell \Pr \{ \bs{D}_\ell = 1 \} < \iy, \qqu \ee \bel  \limsup_{\vep \downarrow 0} \sum_{\ell = i + 1}^\iy  n_\ell \Pr \{ \bs{l} = \ell \} &
\leq & \limsup_{\vep \downarrow 0} \sum_{\ell =
i}^\iy n_{\ell +1} \Pr \{ \bs{l} > \ell \} \nonumber\\
&  \leq & \limsup_{\vep \downarrow 0} \sum_{\ell = i}^\iy n_{\ell +1} \Pr \{ \bs{D}_\ell = 0 \} \nonumber\\
&  \leq &  \limsup_{\vep \downarrow 0} \li [ \li ( \sup_{ \ell \geq i } \f{ n_{\ell +1} }{n_\ell} \ri ) \li ( \sum_{\ell = i}^\iy n_{\ell} \Pr
\{ \bs{D}_\ell = 0 \} \ri ) \ri ] \nonumber\\
&  \leq & \li [ \sup_{ \ell \geq i } \f{ C_{\ell +1}  \Up_{\ell +1} }{C_\ell \Up_\ell} \ri ] \limsup_{\vep \downarrow 0}  \sum_{\ell = i}^\iy
n_{\ell} \Pr \{
\bs{D}_\ell = 0 \} \nonumber\\
& < & \iy,  \la{citelaterb} \eel
From (\ref{citelatera}) and (\ref{citelaterb}), we have \be \la{comb89930}
 \limsup_{\vep \downarrow 0} \li [ \sum_{\ell = \tau_\vep }^{i - 2} n_\ell \Pr \{ \bs{l} = \ell \} + \sum_{\ell = i + 1}^\iy
n_\ell \Pr \{ \bs{l} = \ell \} \ri ] < \iy. \ee From (\ref{same8}), we have $\limsup_{\vep \downarrow 0} n_i \Pr \{ \bs{D}_i = 0 \} < \iy$,
which implies that \be \la{now8936a}
 \lim_{\vep
\downarrow 0} \Pr \{ \bs{D}_i = 0 \} = 0, \ee since $n_i \to \iy$ as $\vep \downarrow 0$.  From (\ref{same869}), we have \be \la{now8936b}
\lim_{\vep \downarrow 0} \sum_{\ell = \tau_\vep}^{i - 2}  \Pr \{ \bs{D}_\ell = 1 \} = 0, \ee since $n_{\tau_\vep}  \to \iy$ as $\vep \downarrow
0$.  It follows from (\ref{now8936a}) and (\ref{now8936b}) that {\small \[ \limsup_{\vep \downarrow 0} \Pr \{ \bs{l} \neq i - 1, \; \bs{l} \neq
i \} \leq  \limsup_{\vep \downarrow 0} [ \Pr \{ \bs{l} > i \} + \Pr \{  \bs{l} < i - 1 \} ] \leq  \limsup_{\vep \downarrow 0} \li [ \Pr \{
\bs{D}_i = 0 \} + \sum_{\ell = \tau_\vep}^{i - 2}  \Pr \{ \bs{D}_\ell = 1 \} \ri ] = 0,
\]}
which implies that \be \la{comb8993a}
 \lim_{\vep \downarrow 0} \Pr \{ \bs{l} \neq i - 1, \; \bs{l} \neq i \} = 0.
\ee By the definition of $n_{i-1}$, we have that \be \la{comb8993b}
 \limsup_{\vep \downarrow 0} \f{n_{i-1}} { \mcal{N} (\vep, \de, \mu, \nu)
} < \iy. \ee Combing (\ref{comb8993a}) and (\ref{comb8993b})  yields \be \la{comb8993c} \lim_{\vep \downarrow 0} \f{n_{i-1} \Pr \{ \bs{l} \neq i
- 1, \; \bs{l} \neq i \} } { \mcal{N} (\vep, \de, \mu, \nu) }  = 0. \ee By virtue of (\ref{comb89930}) and (\ref{comb8993c}), we have  {\small
\bee \limsup_{\vep \downarrow 0} \f{ \bb{E} [ \mbf{n} ] } { \mcal{N} (\vep, \de, \mu, \nu)  } & =  & \limsup_{\vep \downarrow
0} \f{ \sum_{\ell = \tau_\vep }^{i - 2} n_\ell \Pr \{ \bs{l} = \ell \} + \sum_{\ell = i + 1}^\iy n_\ell \Pr \{ \bs{l} = \ell \}  } { \mcal{N} (\vep, \de, \mu, \nu) }\\
&  & - \lim_{\vep \downarrow 0} \f{n_{i-1} \Pr \{ \bs{l} \neq i - 1, \; \bs{l} \neq i \} } { \mcal{N} (\vep, \de, \mu, \nu) }\\
&  & + \limsup_{\vep \downarrow 0} \f{ n_{i - 1} [ 1 - \Pr \{ \bs{l} = i \} ] + n_i \Pr \{ \bs{l} = i \}  } { \mcal{N} (\vep, \de, \mu, \nu) }\\
& = & \limsup_{\vep \downarrow 0} \f{ n_{i - 1} [1 - \Pr \{ \bs{l} = i \} ] + n_i \Pr \{ \bs{l} = i \} } { \mcal{N} (\vep, \de, \mu, \nu) }\\
& = & \limsup_{\vep \downarrow 0} \f{ n_{i - 1} + ( n_i - n_{i - 1} ) \Pr \{ \bs{l} = i \} } { \mcal{N} (\vep, \de, \mu, \nu) }\\
& = & \f{1}{ \mcal{Z}^2 } \li [ \xi_{i-1}^2 + (\xi_i^2 - \xi_{i-1}^2) \limsup_{\vep \downarrow 0} \Pr \{ \bs{l} = i \} \ri ] \leq \li ( \f{
\xi_i }{ \mcal{Z} } \ri )^2. \eee} On the other hand, by virtue of (\ref{comb89930}) and (\ref{comb8993c}), we have {\small \bee \liminf_{\vep
\downarrow 0} \f{ \bb{E} [ \mbf{n} ] } { \mcal{N} (\vep, \de, \mu, \nu)  } & =  & \liminf_{\vep \downarrow
0} \f{ \sum_{\ell = \tau_\vep }^{i - 2} n_\ell \Pr \{ \bs{l} = \ell \} + \sum_{\ell = i + 1}^\iy n_\ell \Pr \{ \bs{l} = \ell \}  } { \mcal{N} (\vep, \de, \mu, \nu) }\\
&  & - \lim_{\vep \downarrow 0} \f{n_{i-1} \Pr \{ \bs{l} \neq i - 1, \; \bs{l} \neq i \} } { \mcal{N} (\vep, \de, \mu, \nu) }\\
&  & + \liminf_{\vep \downarrow 0} \f{ n_{i - 1} [ 1 - \Pr \{ \bs{l} = i \} ] + n_i \Pr \{ \bs{l} = i \}  } { \mcal{N} (\vep, \de, \mu, \nu) }\\
& = & \liminf_{\vep \downarrow 0} \f{ n_{i - 1} [1 - \Pr \{ \bs{l} = i \} ] + n_i \Pr \{ \bs{l} = i \} } { \mcal{N} (\vep, \de, \mu, \nu) }\\
& = & \liminf_{\vep \downarrow 0} \f{ n_{i - 1} + ( n_i - n_{i - 1} ) \Pr \{ \bs{l} = i \} } { \mcal{N} (\vep, \de, \mu, \nu) }\\
& = & \f{1} { \mcal{Z}^2 } \li [ \xi_{i-1}^2 + (\xi_i^2 - \xi_{i-1}^2) \liminf_{\vep \downarrow 0} \Pr \{ \bs{l} = i \} \ri ] \geq \li ( \f{
\xi_{i-1} }{ \mcal{Z} } \ri )^2. \eee}

Now, consider $\mu$ such that (\ref{same8}) and (\ref{same6c8}) hold.  Note that \[ \bb{E} [ \mbf{n} ] \leq \sum_{\ell = \tau_\vep}^{i - 1}
n_\ell \Pr \{ \bs{l} = \ell \} + \sum_{\ell = i + 1}^\iy n_\ell \Pr \{ \bs{l} = \ell \} + n_i. \] From preceding analysis, we have \[ \sum_{\ell
= i + 1}^\iy n_\ell \Pr \{ \bs{l} = \ell \} < \iy.  \]  From (\ref{same6c8}), we have
\[
\sum_{\ell = \tau_\vep}^{i - 1} n_\ell \Pr \{ \bs{l} = \ell \} \leq \sum_{\ell =  \tau_\vep}^{i - 1} n_\ell \Pr \{ \bs{D}_\ell = 1 \} < \iy.
\]
 Therefore,  {\small \[ \limsup_{\vep \downarrow 0} \f{ \bb{E} [ \mbf{n} ]
} { \mcal{N} (\vep, \de, \mu, \nu)  } \leq \lim_{\vep \downarrow 0} \f{ \sum_{\ell = \tau_\vep}^{i - 1} n_\ell \Pr \{ \bs{D}_\ell = 1 \} +
\sum_{\ell = i}^\iy n_{\ell + 1} \Pr \{ \bs{D}_\ell = 0 \} + n_i } { \mcal{N} (\vep, \de, \mu, \nu) } = \lim_{\vep \downarrow 0}  \f{ n_i } {
\mcal{N} (\vep, \de, \mu, \nu) }.
\]}
On the other hand, \be \la{otherhand}
 \bb{E} [ \mbf{n} ] \geq n_i \Pr \{ \bs{l} = i \} \geq n_i \li ( 1 - \Pr \{ \bs{D}_i = 0 \} - \sum_{\ell =
\tau_\vep}^{i - 1} \Pr \{ \bs{D}_\ell = 1 \} \ri ). \ee Making use of (\ref{now8936a}),  (\ref{otherhand}), and the observation
\[
\lim_{\vep \downarrow 0} \sum_{\ell =  \tau_\vep}^{i - 1}  \Pr \{ \bs{D}_\ell = 1 \} = 0
\]
derived from (\ref{same6c8}),  we have {\small \[ \liminf_{\vep \downarrow 0} \f{ \bb{E} [ \mbf{n} ] } { \mcal{N} (\vep, \de, \mu, \nu)  } \geq
\lim_{\vep \downarrow 0} \f{ n_i \li ( 1 - \Pr \{ \bs{D}_i = 0 \} - \sum_{\ell = \tau_\vep}^{i - 1} \Pr \{ \bs{D}_\ell = 1 \} \ri ) } { \mcal{N}
(\vep, \de, \mu, \nu) } = \lim_{\vep \downarrow 0} \f{ n_i } { \mcal{N} (\vep, \de, \mu, \nu) }. \] } So,
\[
\lim_{\vep \downarrow 0} \f{ \bb{E} [ \mbf{n} ] } { \mcal{N} (\vep, \de, \mu, \nu)  } = \lim_{\vep \downarrow 0} \f{ n_i } { \mcal{N} (\vep,
\de, \mu, \nu) } = \li (  \f{\xi_i}{ \mcal{Z}  } \ri )^2
\] for $\mu$ such that (\ref{same8}) and (\ref{same6c8}) hold.  This completes the proof of the lemma.

\epf

\beL

\la{lem22888}

If there exists an index $i$ such that $\Pr \{ i - 1 \leq \liminf_{\vep \downarrow 0} \bs{l} \leq \limsup_{\vep \downarrow 0} \bs{l} \leq i \} =
1$, then $\liminf_{\vep \downarrow 0 } \Pr \{  \mu \in \bs{\mcal{I}} \} \geq 2 \sum_{\ell = i - 1}^i \Phi \li ( \xi_\ell \ri ) - 3$. Similarly,
if there exists an index $i$ such that $\Pr \li \{ \lim_{\vep \downarrow 0} \bs{l} = i \ri \} = 1$, then $\lim_{\vep \downarrow 0 } \Pr \{  \mu
\in \bs{\mcal{I}} \} = 2 \Phi \li ( \xi_i \ri ) - 1$.

\eeL

\bpf

If $\Pr \li \{ i - 1 \leq \liminf_{\vep \downarrow 0} \bs{l} \leq \limsup_{\vep \downarrow 0} \bs{l}  \leq i \ri \} = 1$, then $\lim_{\vep
\downarrow 0 } \Pr \{  \bs{l} < i - 1 \; \tx{or} \; \bs{l} > i \} = 0$.  It follows that \bee \limsup_{\vep \downarrow 0 } \Pr \{  \mu \notin
\bs{\mcal{I}} \} & \leq & \limsup_{\vep \downarrow 0 } \Pr \{  \mu \notin
\bs{\mcal{I}}, \; i - 1 \leq \bs{l} \leq i \} +  \limsup_{\vep \downarrow 0 } \Pr \{  \bs{l} < i - 1 \; \tx{or} \; \bs{l} > i \}\\
& = &  \limsup_{\vep \downarrow 0 } \Pr \{  \mu \notin \bs{\mcal{I}}, \; i - 1 \leq \bs{l} \leq i \} \\
& \leq & \limsup_{\vep \downarrow 0 } \Pr \{ \mu \notin \bs{\mcal{I}}_{i-1} \} + \limsup_{\vep \downarrow 0 } \Pr \{ \mu \notin \bs{\mcal{I}}_i \}\\
& = & \lim_{\vep \downarrow 0 } \Pr \{ \mu \notin \bs{\mcal{I}}_{i-1} \} + \lim_{\vep \downarrow 0 } \Pr \{ \mu \notin \bs{\mcal{I}}_i \}\\
& = &  2 [ 1 - \Phi (\xi_{i-1}) ] + 2 [ 1 - \Phi (\xi_i) ], \eee which implies that $\liminf_{\vep \downarrow 0 } \Pr \{  \mu \in \bs{\mcal{I}}
\} \geq 2 \sum_{\ell = i - 1}^i \Phi \li ( \xi_\ell \ri ) - 3$.

If $\Pr \li \{ \lim_{\vep \downarrow 0} \bs{l} = i \ri \} = 1$, then $\lim_{\vep \downarrow 0} \Pr \li \{ \bs{l} \neq i \ri \} = 0$.  It follows
that \bee \limsup_{\vep \downarrow 0} \Pr \{ \mu \in \bs{\mcal{I}} \} & \leq & \limsup_{\vep \downarrow 0} \Pr \{ \mu \in \bs{\mcal{I}},
\; \bs{l} = i \} + \limsup_{\vep \downarrow 0} \Pr \{ \bs{l} \neq i \}\\
& = &  \limsup_{\vep \downarrow 0} \Pr \{ \mu \in \bs{\mcal{I}}_i, \; \bs{l} = i \} \leq  \limsup_{\vep \downarrow 0} \Pr \{ \mu \in
\bs{\mcal{I}}_i \},  \eee \bee \liminf_{\vep \downarrow 0} \Pr \{ \mu \in \bs{\mcal{I}} \}
& \geq & \liminf_{\vep \downarrow 0} \Pr \{ \mu \in \bs{\mcal{I}}, \; \bs{l} = i \} =  \liminf_{\vep \downarrow 0} \Pr \{ \mu \in \bs{\mcal{I}}_i, \; \bs{l} = i \}\\
& \geq &  \liminf_{\vep \downarrow 0} \Pr \{ \mu \in \bs{\mcal{I}}_i \} - \limsup_{\vep \downarrow 0} \Pr \{ \bs{l} \neq i \}\\
& = &  \liminf_{\vep \downarrow 0} \Pr \{ \mu \in \bs{\mcal{I}}_i \}.   \eee Since
\[
\liminf_{\vep \downarrow 0} \Pr \{ \mu \in \bs{\mcal{I}}_i \} = \limsup_{\vep \downarrow 0} \Pr \{ \mu \in \bs{\mcal{I}}_i \} = \lim_{\vep
\downarrow 0} \Pr \{ \mu \in \bs{\mcal{I}}_i \} = 2 \Phi (\xi_i) - 1,
\]
it must be true that
\[
\lim_{\vep \downarrow 0} \Pr \{ \mu \in \bs{\mcal{I}} \} = \lim_{\vep \downarrow 0} \Pr \{ \mu \in \bs{\mcal{I}}_i \} = 2 \Phi (\xi_i) - 1.
\]
This completes the proof of the lemma.  \epf

\beL

\la{lem23888}

For $\mu \in \Se$, if there exists an index $i \geq 1$ such that \bee &  & \lim_{\vep \downarrow 0} \sum_{\ell = \tau_\vep}^{i - 2} \Pr \{
\bs{D}_\ell = 1 \} = 0, \qqu \lim_{\vep \downarrow 0} \Pr \{ \bs{D}_{\ell} = 0 \} = 0 \qu \tx{for all $\ell \geq i$}, \eee then,  for all $\ell
\in \bb{Z}$,
 \[
 0 \leq \Pr \{ \bs{D}_\ell = 0 \} - \Pr \{ \bs{l} > \ell \}  \to 0 \qqu \tx{as $\vep \downarrow 0$}.
\]

\eeL

\bpf

By the definition of the stopping rule,
\[
\Pr \{ \bs{D}_\ell = 0 \} \geq \Pr \{  \bs{l} > \ell  \}  \geq \Pr \{ \bs{D}_\ell = 0 \}  - \sum_{j = \tau_\vep}^{\ell - 1} \Pr \{  \bs{D}_j = 1
\}. \] For $\ell < i$, we have
\[
0 \leq \sum_{j = \tau_\vep}^{\ell - 1} \Pr \{  \bs{D}_j = 1 \} \leq \sum_{j = \tau_\vep}^{i - 2} \Pr \{  \bs{D}_j = 1 \} \to 0 \qu \tx{as $\vep
\downarrow 0$},
\]
which implies $\lim_{ \vep \downarrow 0  } [ \Pr \{ \bs{D}_\ell = 0 \} - \Pr \{ \bs{l} > \ell \} ] = 0$. For $\ell \geq i$, we have
\[
0 \leq \Pr \{  \bs{l} > \ell  \} \leq \Pr \{ \bs{D}_\ell = 0 \} \to 0 \qu \tx{as $\vep \downarrow 0$},
\]
which also implies $\lim_{ \vep \downarrow 0  } [ \Pr \{ \bs{D}_\ell = 0 \} - \Pr \{ \bs{l} > \ell \} ] = 0$. This completes the proof of the
lemma.

\epf

\beL  \la{lem24888} For any $\mu \in \Se$, \be \la{obvious89}
 \udl{P} \leq \Pr \{ \mu \notin \bs{\mcal{I}} \} \leq \ovl{P} \leq Q.
\ee For any $\mu \in \Se$, if there exists an index $i \geq 1$ such that \bel &  &  \lim_{\vep \downarrow 0} \sum_{\ell = \tau_\vep}^{i -
2} \Pr \{ \bs{D}_\ell = 1 \} = 0, \la{howgooda} \\
&  & \lim_{\vep \downarrow 0} \sum_{\ell = i}^\iy \Pr \{ \bs{D}_{\ell} = 0 \} = 0,  \la{howgoodb} \eel then,  \be \la{display89} \lim_{\vep
\downarrow 0} | \Pr \{ \mu \notin \bs{\mcal{I}} \} - \ovl{P} | = \lim_{\vep \downarrow 0}  | \Pr \{ \mu \notin \bs{\mcal{I}} \} - \udl{P} | = 0.
\ee For any $\mu \in \Se$, if there exists an index $i \geq 1$ such that (\ref{howgoodb}) holds and \be  \lim_{\vep \downarrow 0} \sum_{\ell =
\tau_\vep}^{i - 1} \Pr \{ \bs{D}_\ell = 1 \} = 0, \la{asp896a} \ee then, \be  \la{ineqone98a}
 \lim_{\vep \downarrow 0} [ Q -  \Pr \{ \mu
\notin \bs{\mcal{I}} \} ] = 0. \ee

\eeL

\bpf

Clearly, (\ref{obvious89}) follows from the definition of the stopping rule.  To show (\ref{display89}), it suffices to show \be \la{show1}
 \lim_{\vep \downarrow 0}  \sum_{\ell
= \tau_\vep}^\iy \Pr \{ \bs{D}_{\ell - 1} = 0, \; \bs{D}_\ell = 1 \} = 1. \ee This is because $\udl{P} \leq \Pr \{ \mu  \notin \bs{\mcal{I}} \}
\leq \ovl{P}$ and $\ovl{P} - \udl{P} = \sum_{\ell = \tau_\vep}^\iy \Pr \{ \bs{D}_{\ell - 1} = 0, \; \bs{D}_\ell = 1 \} - 1$. Observing that
{\small \bee & & \sum_{\ell = \tau_\vep}^{i - 2}  \Pr \{
\bs{D}_{\ell - 1} = 0, \; \bs{D}_\ell = 1 \} \leq  \sum_{\ell = \tau_\vep}^{i - 2} \Pr \{ \bs{D}_\ell = 1 \},\\
&   & \sum_{\ell = i + 1}^\iy  \Pr \{ \bs{D}_{\ell - 1} = 0, \; \bs{D}_\ell = 1 \} \leq \sum_{\ell = i + 1}^\iy  \Pr \{ \bs{D}_{\ell - 1} = 0 \}
= \sum_{\ell = i}^\iy  \Pr \{ \bs{D}_{\ell} = 0 \} \eee} and using (\ref{howgooda}) and (\ref{howgoodb}), we have {\small \bee \lim_{\vep
\downarrow 0} \sum_{\ell = \tau_\vep}^{i - 2}  \Pr \{ \bs{D}_{\ell - 1} = 0, \; \bs{D}_\ell = 1 \} = 0, \qqu  \lim_{\vep \downarrow 0}
\sum_{\ell = i + 1}^\iy \Pr \{ \bs{D}_{\ell - 1} = 0, \; \bs{D}_\ell = 1 \} = 0. \eee}  Hence, to show (\ref{show1}), it suffices to show \be
\la{vip988} \lim_{\vep \downarrow 0} \sum_{\ell = i - 1}^i \Pr \{ \bs{D}_{\ell - 1} = 0, \; \bs{D}_\ell = 1 \} = 1. \ee Observing that  \bee \Pr
\{ \bs{D}_{i - 2} = \bs{D}_{i - 1} = 1 \} + \Pr \{ \bs{D}_{i - 1} = \bs{D}_i = 0 \} + \sum_{\ell = i - 1}^i \Pr \{ \bs{D}_{\ell - 1} = 0, \;
\bs{D}_\ell = 1 \} =  1, \eee we have \be \la{cite8963}
  \sum_{\ell = i - 1}^i \Pr \{ \bs{D}_{\ell - 1} = 0,  \; \bs{D}_\ell = 1 \} = 1 - \Pr \{
\bs{D}_{i - 2} = \bs{D}_{i - 1} = 1 \} - \Pr \{ \bs{D}_{i - 1} = \bs{D}_i = 0 \}. \ee So, to show (\ref{vip988}), it is sufficient  to show \bel
& & \lim_{\vep \downarrow 0} \Pr \{ \bs{D}_{i - 2} = \bs{D}_{i - 1} = 1 \} = 0,
\la{use8a9936} \\
 &  & \lim_{\vep \downarrow 0}  \Pr \{ \bs{D}_{i - 1} = \bs{D}_i = 0 \} = 0. \la{use8b9936}
\eel  Invoking  (\ref{howgooda}), we have {\small $\lim_{\vep \downarrow 0} \Pr \{ \bs{D}_{i - 2} = \bs{D}_{i - 1} = 1 \} \leq \lim_{\vep
\downarrow 0} \Pr \{ \bs{D}_{i - 2} =  1 \} = 0$}, which implies (\ref{use8a9936}). Using (\ref{howgoodb}), we have {\small $\lim_{\vep
\downarrow 0} \Pr \{ \bs{D}_{i - 1} = \bs{D}_i = 0 \} \leq \lim_{\vep \downarrow 0} \Pr \{ \bs{D}_i = 0 \} = 0$},  which implies
(\ref{use8b9936}).  It follows that (\ref{display89}) holds.

Now we shall show (\ref{ineqone98a}).  For simplicity of notations, define
\[
Q_\ell =  \min \li [ \Pr \{ \mu \notin  \bs{\mcal{I}}_\ell \}, \; \Pr \{ \bs{D}_{\ell - 1} = 0 \}, \; \Pr \{ \bs{D}_\ell = 1 \} \ri ] \qu
\tx{for $\ell \geq \tau_\vep$}.
\]
Clearly,  from (\ref{asp896a}),
\[
\sum_{\ell = \tau_\vep}^{i - 1} Q_\ell \leq \sum_{\ell = \tau_\vep}^{i - 1} \Pr \{ \bs{D}_\ell = 1 \}  \to 0 \qu \tx{as $\vep \downarrow 0$}
\]
and,  from (\ref{howgoodb}),
\[
\sum_{\ell = i + 1}^\iy Q_\ell \leq \sum_{\ell = i + 1}^\iy \Pr \{ \bs{D}_{\ell - 1} = 0 \}  = \sum_{\ell = i}^\iy \Pr \{ \bs{D}_{\ell} = 0 \}
\to 0 \qu \tx{as $\vep \downarrow 0$}.
\]
On the other hand, from (\ref{asp896a}),
\[
\sum_{\ell = \tau_\vep}^{i - 1} \Pr \{ \mu \notin \bs{\mcal{I}}, \; \bs{l} = \ell \} \leq \sum_{\ell = \tau_\vep}^{i - 1} \Pr \{ \bs{D}_\ell = 1
\}  \to 0 \qu \tx{as $\vep \downarrow 0$}
\]
and,  from (\ref{howgoodb}),
\[
\sum_{\ell = i + 1}^\iy \Pr \{ \mu \notin \bs{\mcal{I}}, \; \bs{l} = \ell \} \leq \sum_{\ell = i + 1}^\iy \Pr \{ \bs{D}_{\ell - 1} = 0 \} =
\sum_{\ell = i}^\iy \Pr \{ \bs{D}_{\ell} = 0 \} \to 0 \qu \tx{as $\vep \downarrow 0$}.
\]
Hence, to show (\ref{ineqone98a}), it suffices  to show that \[
 \lim_{\vep \downarrow 0} [ Q_i -  \Pr \{ \mu
\notin \bs{\mcal{I}}, \; \bs{l} = i \} ] = 0. \]  This can be accomplished by establishing \be \la{est89a} \lim_{\vep \downarrow 0} Q_i = \Pr \{
\mu \notin \bs{\mcal{I}}_i \} \ee and \be \la{est89b}
 \lim_{\vep \downarrow 0} \Pr \{ \mu \notin \bs{\mcal{I}}, \; \bs{l} = i \} = \lim_{\vep
\downarrow 0} \Pr \{ \mu \notin \bs{\mcal{I}}_i\}. \ee From (\ref{asp896a}), \be \la{ok89a}
 \lim_{\vep \downarrow 0} \Pr \{ \bs{D}_{i - 1} = 0
\} = 1, \ee From (\ref{howgoodb}), \be \la{ok89b} \lim_{\vep \downarrow 0} \Pr \{ \bs{D}_i = 1 \} = 1. \ee Combing (\ref{ok89a}) and
(\ref{ok89b}) yields (\ref{est89a}).  On the other hand,
\[
\limsup_{\vep \downarrow 0} \Pr \{ \bs{l} \neq i  \} \leq \limsup_{\vep \downarrow 0} \li [ \Pr \{ \bs{D}_i = 0 \} + \sum_{\ell = \tau_\vep}^{i
- 1} \Pr \{ \bs{D}_\ell = 1 \} \ri ] = 0
\]
as a consequence of (\ref{howgoodb}), (\ref{asp896a}) and the definition of the stopping rule.  Hence, \be \la{vip8963}
 \lim_{\vep \downarrow 0} \Pr \{ \bs{l} = i  \} = 1
\ee and it follows that
\[
\lim_{\vep \downarrow 0} \Pr \{ \mu \notin \bs{\mcal{I}}, \; \bs{l} = i \}  = \lim_{\vep \downarrow 0} \Pr \{ \mu \notin \bs{\mcal{I}}_i, \;
\bs{l} = i \} = \lim_{\vep \downarrow 0} \Pr \{ \mu \notin \bs{\mcal{I}}_i\}.
\]
So, (\ref{est89b}) holds and thus (\ref{ineqone98a}) is established. The proof of the lemma is thus completed.

\epf

\beL

\la{lem25888}

 For $\mu \in \Se$, if there exists an index $i \geq 1$ such that \be  \limsup_{\vep \downarrow 0} \sum_{\ell = \tau_\vep}^{i -
2} n_\ell \Pr \{ \bs{D}_\ell = 1 \} < \iy, \qqu \limsup_{\vep \downarrow 0} \sum_{\ell = i}^\iy n_\ell \Pr \{ \bs{D}_{\ell} = 0 \} < \iy,
\la{use8936} \ee then, \[ \lim_{\vep \downarrow 0} \f{ 1 }{\bb{E} [ \mbf{n} ] } \li [ n_{\tau_\vep} + \sum_{\ell = \tau_\vep}^\iy (n_{\ell + 1}
- n_\ell)  \Pr \{ \bs{D}_\ell = 0 \} \ri ] = 1. \]

\eeL

\bpf

By the definition of the stopping rule, \bee &  & \bb{E} [ \mbf{n} ] \leq  n_{\tau_\vep} + \sum_{\ell = \tau_\vep}^\iy (n_{\ell + 1} - n_\ell)  \Pr \{ \bs{D}_\ell = 0 \}, \\
&  & \bb{E} [ \mbf{n} ] \leq \sum_{\ell = \tau_\vep}^\iy n_\ell \Pr \{ \bs{D}_{\ell - 1} = 0, \; \bs{D}_\ell = 1 \}. \eee Since \bee & & \li [
n_{\tau_\vep} + \sum_{\ell = \tau_\vep}^\iy (n_{\ell + 1} - n_\ell)  \Pr \{ \bs{D}_\ell = 0  \} \ri ] - \li [ \sum_{\ell = \tau_\vep}^\iy n_\ell
\Pr \{ \bs{D}_{\ell - 1} = 0, \; \bs{D}_\ell = 1  \} \ri ]\\
&  & =  \sum_{\ell = \tau_\vep}^\iy n_\ell \Pr \{ \bs{D}_{\ell - 1} = 0, \; \bs{D}_\ell = 0  \}  - \sum_{\ell = \tau_\vep}^\iy n_\ell \Pr \{
\bs{D}_{\ell} = 0  \} \\
&  & = - \sum_{\ell = \tau_\vep}^\iy n_\ell \Pr \{ \bs{D}_{\ell - 1} = 1, \; \bs{D}_\ell = 0  \} \leq 0, \eee it suffices to show that \be
\la{suf89}
 \lim_{\vep \downarrow 0} R_\vep  = 1, \ee   where \bee R_\vep = \f { \bb{E} [ \mbf{n} ]  }{ \sum_{\ell = \tau_\vep}^\iy n_\ell \Pr \{ \bs{D}_{\ell - 1} = 0, \;
\bs{D}_\ell = 1  \} } = \f {  \sum_{\ell = \tau_\vep}^\iy n_\ell \Pr \{ \bs{l} = \ell \} }{ \sum_{\ell = \tau_\vep}^\iy n_\ell \Pr \{
\bs{D}_{\ell - 1} = 0, \; \bs{D}_\ell = 1  \} } \leq 1. \eee   To show (\ref{suf89}), it suffices to show $\liminf_{\vep \downarrow 0} R_\vep =
1$. Note that \be \la{alpha89}
 \limsup_{\vep \downarrow 0} \sum_{\ell = \tau_\vep}^{i - 2}  n_\ell \Pr \{ \bs{D}_{\ell - 1} = 0, \; \bs{D}_\ell
= 1 \} \leq \limsup_{\vep \downarrow 0} \sum_{\ell = \tau_\vep}^{i - 2}  n_\ell \Pr \{  \bs{D}_\ell = 1  \} < \iy, \ee {\small \bel
\limsup_{\vep \downarrow 0} \sum_{\ell = i + 1}^\iy n_\ell \Pr \{ \bs{D}_{\ell - 1} = 0, \; \bs{D}_\ell = 1  \} & \leq & \limsup_{\vep
\downarrow 0} \sum_{\ell = i}^\iy n_{\ell + 1} \Pr \{ \bs{D}_{\ell} = 0  \} \nonumber\\
& \leq &  \limsup_{\vep \downarrow 0} \; \li [ \sup_{\ell \geq i} \f{n_{\ell+1}}{n_\ell} \sum_{\ell = i}^\iy n_\ell \Pr \{ \bs{D}_{\ell} = 0  \}
\ri ] < \iy. \la{beta89} \eel} From (\ref{use8936}), we have
\[
\limsup_{\vep \downarrow 0} n_{i-2} \Pr \{ \bs{D}_{i-2} = 1 \} < \iy, \qqu \limsup_{\vep \downarrow 0} n_i \Pr \{ \bs{D}_i = 0 \} < \iy,
\]
which implies that \be \la{nowcite66}
 \lim_{\vep \downarrow 0}  \Pr \{ \bs{D}_{i-2} = 1 \} = 0, \qqu \lim_{\vep \downarrow 0} \Pr \{ \bs{D}_i =
0 \} = 0, \ee since $n_{i-2} \to \iy$ as $\vep \downarrow 0$. By virtue of (\ref{cite8963}) and (\ref{nowcite66}), we have \be \la{quot83}
 \lim_{\vep
\downarrow 0} \sum_{\ell = i - 1}^i \Pr \{ \bs{D}_{\ell - 1} = 0, \; \bs{D}_\ell = 1  \} = 1, \ee from which we have \be \la{gama89}
 \sum_{\ell = i - 1}^i n_\ell \Pr \{ \bs{D}_{\ell - 1} = 0, \; \bs{D}_\ell = 1  \} \geq n_{i-1} \sum_{\ell = i - 1}^i \Pr \{ \bs{D}_{\ell - 1} = 0, \;
\bs{D}_\ell = 1 \} \to \iy \ee as $\vep \downarrow 0$.  In view of (\ref{comb89930}), (\ref{alpha89}),  (\ref{beta89}) and (\ref{gama89}), to
show (\ref{suf89}), it suffices to show \be \la{pass89}
 \liminf_{\vep \downarrow 0} \f{  \sum_{\ell = i - 1}^i n_\ell \Pr \{ \bs{l} = \ell \} } {  \sum_{\ell = i - 1}^i n_\ell \Pr
\{ \bs{D}_{\ell - 1} = 0, \; \bs{D}_\ell = 1  \} }  = 1. \ee Observing that
\[
\Pr \{ \bs{l} = j \}  \geq \Pr \{ \bs{D}_{j - 1} = 0, \; \bs{D}_j = 1  \} - \sum_{\ell = \tau_\vep}^{j-2} \Pr \{  \bs{D}_\ell = 1 \} \] for $j
\in \{ i - 1, \; i \}$, we have \bee \sum_{\ell = i - 1}^i  n_\ell \Pr \{ \bs{l} = \ell \}  & \geq & \sum_{\ell = i - 1}^i  n_\ell \Pr \{
\bs{D}_{\ell - 1} = 0, \; \bs{D}_\ell = 1  \} - n_i  \sum_{\ell = \tau_\vep}^{i-2} \Pr \{ \bs{D}_\ell
= 1 \} -  n_{i-1}   \sum_{\ell = \tau_\vep}^{i-3} \Pr \{ \bs{D}_\ell = 1 \}\\
& \geq & \sum_{\ell = i - 1}^i  n_\ell \Pr \{ \bs{D}_{\ell - 1} = 0, \; \bs{D}_\ell = 1  \} - 2 n_i  \sum_{\ell = \tau_\vep}^{i-2} \Pr \{
\bs{D}_\ell = 1 \}.  \eee So, to show (\ref{pass89}), it suffices to show $\lim_{\vep \downarrow 0} Q_\vep = 0$, where
\[
Q_\vep = \f{ 2 n_i \sum_{\ell = \tau_\vep}^{i-2} \Pr \{  \bs{D}_\ell = 1 \} } { \sum_{\ell = i - 1}^i  n_\ell \Pr \{ \bs{D}_{\ell - 1} = 0, \;
\bs{D}_\ell = 1  \}  }.
\]
From (\ref{use8936}), we have
\[
0 \leq \limsup_{\vep \downarrow 0} \sum_{\ell = \tau_\vep}^{i-2} \Pr \{  \bs{D}_\ell = 1 \} \leq \limsup_{\vep \downarrow 0}
\f{1}{n_{\tau_\vep}} \sum_{\ell = \tau_\vep}^{i-2} n_\ell \Pr \{  \bs{D}_\ell = 1 \} = 0,
\]
which implies \be \la{qto968}
 \lim_{\vep \downarrow 0} \sum_{\ell = \tau_\vep}^{i-2} \Pr \{  \bs{D}_\ell = 1 \} = 0. \ee
Making use of (\ref{quot83}), (\ref{qto968}) and the observation that
\[
\f{ 2 n_i  \sum_{\ell = \tau_\vep}^{i-2} \Pr \{  \bs{D}_\ell = 1 \} } { n_i \sum_{\ell = i - 1}^i  \Pr \{ \bs{D}_{\ell - 1} = 0, \; \bs{D}_\ell
= 1 \} } \leq Q_\vep \leq \f{ 2 n_i  \sum_{\ell = \tau_\vep}^{i-2} \Pr \{  \bs{D}_\ell = 1 \} } { n_{i-1} \sum_{\ell = i - 1}^i \Pr \{
\bs{D}_{\ell - 1} = 0, \; \bs{D}_\ell = 1  \}  },
\]
we have $\lim_{\vep \downarrow 0} Q_\vep = 0$.  Therefore, we have shown (\ref{suf89}).  The proof of the lemma is thus completed.  \epf

\sect{Proof of Theorem \ref{CHCDF}}  \la{CHCDF_app}

The proof of (\ref{VIP886c}) is similar to that of (\ref{VIP886}).  The other properties are proved in the sequel.

\subsection{Proof of Property (I)}

To show (\ref{perty891}), we need the following preliminary result.

\beL \la{asn889finite} For $\mu \in \Se$, there exist a number $\ga > 0$ and a positive integer $m > \tau_\vep$  such that $\{ \bs{l} > l \}
\subseteq \{ | Z_\ell - \mu | \geq \ga   \}$ for all $\ell > m$.

\eeL

\bpf

To show the lemma, we need to prove two claims as follows.

(i): For $\mu \in \Se$, there exist a number $\ga > 0$ and  a positive integer $m > \tau_\vep$  such that
\[
\li \{  F_{Z_\ell} ( Z_\ell, \mscr{U} (Z_\ell, \vep) ) > \f{\de_\ell}{2} \ri \} \subseteq \{  | Z_\ell - \mu | \geq \ga \} \qu \tx{for all $\ell
> m$}.
\]

(ii): For $\mu \in \Se$, there exist a number $\ga > 0$ and  a positive integer $m > \tau_\vep$  such that
\[
\li \{  G_{Z_\ell} ( Z_\ell, \mscr{L} (Z_\ell, \vep) ) > \f{\de_\ell}{2} \ri \} \subseteq \{  | Z_\ell - \mu | \geq \ga \} \qu \tx{for all $\ell
> m$}.
\]

To prove the first claim, note that there exists $\ga > 0$ such that  $\mu - \ga, \; \mu + \ga \in \Se$ and that $\mscr{U} (\se, \vep) - \se
\geq \f{ \mscr{U} (\mu, \vep) - \mu }{2} > 0$ for all $\se \in (\mu - \ga, \mu + \ga) \subseteq \Se$.  Define $S = \{ \se \in (\mu - \ga, \mu +
\ga): \mscr{U} (\se, \vep) \in \Se \}$. If $S = \emptyset$, then $\li \{  F_{Z_\ell} ( Z_\ell, \mscr{U} (Z_\ell, \vep) )
> \f{\de_\ell}{2}, \; | Z_\ell - \mu | < \ga \ri \} = \emptyset$ for $\ell \geq \tau_\vep$.  In the other case that $S \neq
\emptyset$, we have that there exists a constant $D < \iy$ such that
\[
\sup_{ z \in S } \f{\mscr{W}(\mscr{U} (z, \vep) )}{| z - \mscr{U} ( z, \vep)  |^3} < D
\]
By the assumptions on $\{n_\ell \}$ and $\{\de_\ell \}$, we have that $n_\ell^2 \de_\ell \to \iy$ as $\ell \to \iy$. Hence, for any constant $C
> 0$, \[ \li \{ \mscr{U} (Z_\ell, \vep) \in \Se, \; | Z_\ell - \mu | < \ga  \ri \} = \li \{ \f{C \mscr{W}(\mscr{U} (Z_\ell, \vep) )}{n_\ell^2 | Z_\ell
- \mscr{U} (Z_\ell, \vep)  |^3} < \f{\de_\ell}{4}, \;  \mscr{U} (Z_\ell, \vep) \in \Se, \; | Z_\ell - \mu | < \ga \ri \}
\]
for large enough $\ell$.  Since $n_\ell^2 \de_\ell \to \iy$ as $\ell \to \iy$, we have $\f{1}{n_\ell} \ln \f{2}{\de_\ell} \to 0$ as $\ell \to
\iy$.  It follows from (\ref{beryb1}) of Lemma \ref{nonuniformBE} that
\bee &  &  \li \{  F_{Z_\ell} ( Z_\ell, \mscr{U} (Z_\ell, \vep) ) > \f{\de_\ell}{2}, \; | Z_\ell - \mu | < \ga  \ri \}\\
&  & = \li \{ F_{Z_\ell} ( Z_\ell, \mscr{U} (Z_\ell, \vep) ) > \f{\de_\ell}{2}, \; \mscr{U} (Z_\ell, \vep) \in \Se, \; |
Z_\ell - \mu | < \ga  \ri \}\\
&  & \subseteq \li \{ \f{1}{2} \exp \li ( - \f{ n_\ell | Z_\ell - \mscr{U} (Z_\ell, \vep) |^2 }{  2 \mcal{V} (\mscr{U} (Z_\ell, \vep)) } \ri ) +
\f{C \mscr{W}(\mscr{U} (Z_\ell, \vep) )}{n_\ell^2 | Z_\ell - \mscr{U} (Z_\ell, \vep)  |^3} > \f{\de_\ell}{2}, \;
 \mscr{U} (Z_\ell, \vep) \in \Se, \; | Z_\ell - \mu | < \ga \ri \} \\
&  & \subseteq \li \{  \f{1}{2} \exp \li ( - \f{ n_\ell | Z_\ell - \mscr{U} (Z_\ell, \vep) |^2 }{  2 \mcal{V} (\mscr{U}
(Z_\ell, \vep))   } \ri ) > \f{\de_\ell}{4}, \; \mscr{U} (Z_\ell, \vep) \in \Se, \; | Z_\ell - \mu | < \ga \ri \} \\
&  & \subseteq \li \{  \exp \li ( - \f{ n_\ell | Z_\ell - \mscr{U} (Z_\ell, \vep) |^2 }{ 2 \mcal{V}^* } \ri )
> \f{\de_\ell}{2}, \; | Z_\ell - \mu | < \ga \ri \} \\
&  & \subseteq \li \{  \mscr{U} (Z_\ell, \vep) - Z_\ell < \f{ \mscr{U} (\mu, \vep) - \mu  }{2}, \; | Z_\ell - \mu | < \ga \ri \} = \emptyset
\eee for large enough $\ell$, where $\mcal{V}^* = \sup_{\se \in S} \mcal{V} (\se)$ and $C$ is some absolute constant.  It follows that, in both
cases that $S = \emptyset$ and $S \neq \emptyset$,  \bee \li \{ F_{Z_\ell} ( Z_\ell, \mscr{U} (Z_\ell, \vep) ) > \f{\de_\ell}{2} \ri \} &
\subseteq & \li \{ F_{Z_\ell} ( Z_\ell, \mscr{U} (Z_\ell, \vep) ) > \f{\de_\ell}{2}, \; | Z_\ell - \mu | < \ga \ri \} \cup \{ |
Z_\ell - \mu | \geq \ga \}\\
& = & \{ | Z_\ell - \mu | \geq \ga \}  \eee for large enough $\ell$.  This proves the first claim.  By a similar method, we can show the second
claim.

By the definition of the stopping rule, we have that {\small \be \la{observa8899finite}
 \{ \bs{l} > \ell \} \subseteq \li \{  F_{Z_\ell} (
Z_\ell, \mscr{U} (Z_\ell, \vep) ) > \f{\de_\ell}{2} \ri \} \cup \li \{ G_{Z_\ell} ( Z_\ell, \mscr{L} (Z_\ell, \vep) )
> \f{\de_\ell}{2} \ri \}. \ee} Finally, the proof of the lemma is completed by making use of (\ref{observa8899finite}) and our two established claims.

\epf

We are now in a position to prove (\ref{perty891}).  According to Lemma \ref{asn889finite}, there exists an integer $m > \tau_\vep$ and a
positive number $\ga
> 0$ such that \be \la{makeusefinite}
 \Pr \{ \bs{l} > \ell \}
\leq \Pr \{ | Z_\ell - \mu | \geq \ga \} \qqu \tx{for all $\ell > m$}. \ee Making use of (\ref{makeusefinite}) and the weak law of large
numbers, we have
\[ 0 \leq \limsup_{\ell \to \iy} \Pr \{ \bs{l}
> \ell \} \leq \lim_{\ell \to \iy} \Pr \{ | Z_\ell - \mu | \geq \ga \} = 0,
\]
which implies that $\lim_{\ell \to \iy} \Pr \{ \bs{l} > \ell \} = 0$.  It follows that
\[
\Pr \{ \bs{l} < \iy \} = 1 - \lim_{\ell \to \iy} \Pr \{ \bs{l} > \ell \} = 1
\]
and thus  $\Pr \{ \mbf{n} < \iy \} = 1$.  By virtue of (\ref{mostvip}) of Lemma \ref{nonuniformBE} and (\ref{makeusefinite}), we have \bee
\bb{E} [ \mbf{n}  ]  \leq  n_{m+1} + \sum_{\ell
> m} \Pr \{ | Z_\ell - \mu | \geq \ga \} \leq  n_{m+1}  + \sum_{\ell > m} (n_{\ell + 1} - n_\ell) \li [ \exp \li ( - \f{n_\ell}{2} \f{ \ga^2 }{  \mcal{V} (\mu) } \ri ) +
\f{2 C }{n_\ell^2} \f{ \mscr{W}(\mu) }{ \ga^3 } \ri ] < \iy.  \eee

\subsection{Proof of Property (II)}

To prove (\ref{propertyii1}) and (\ref{propertyii2}), we can apply Lemma \ref{pro689a} and Lemma \ref{pro689b} as follows.

\beL

\la{pro689b}

\bel &  & \Pr \li \{ \liminf_{\vep \downarrow 0} \Lm_{\bs{l}} (\mu) \geq 1 \ri \} = 1, \la{cite369a}\\
&  & \Pr \li \{ \tx{Both} \; \limsup_{\vep \downarrow 0} \Lm_{\bs{l} - 1} (\mu) \leq 1 \; \tx{and} \;  \liminf_{\vep \downarrow 0} \bs{l} > 1 \;
\tx{hold or} \; \lim_{\vep \downarrow 0} \bs{l}  = 1 \ri \} = 1.  \la{cite369b} \eel

\eeL

\bpf To simplify notations, define
\[
\bs{L}_{\bs{l}} = \mscr{L} ( Z_{\bs{l}}, \vep ), \qqu \bs{L}_{\bs{l} - 1} = \mscr{L} ( Z_{\bs{l} - 1}, \vep ), \qqu \bs{U}_{\bs{l}} = \mscr{U} (
Z_{\bs{l}}, \vep ), \qqu \bs{U}_{\bs{l} - 1} = \mscr{U} ( Z_{\bs{l} - 1}, \vep ).
\]
Define {\small \bee &  & \mscr{E}_1 = \li \{ F_{Z_{\bs{l}}} (Z_{\bs{l}}, \bs{U}_{\bs{l}} ) \leq \f{\de_{\bs{l}}}{2}, \; G_{Z_{\bs{l}}}
(Z_{\bs{l}}, \bs{L}_{\bs{l}} ) \leq \f{\de_{\mbf{n}}}{2}  \ri \}, \\
&  & \mscr{E}_2 = \{  \bs{l} = \tau_\vep  \} \bigcup \li \{  F_{Z_{\bs{l}-1}} (Z_{\bs{l}-1}, \bs{U}_{\bs{l}-1} ) > \f{\de_{\bs{l}-1}}{2}, \; \bs{l} > \tau_\vep \ri \}
\bigcup \li \{ G_{Z_{\bs{l}-1}} (Z_{\bs{l}-1}, \bs{L}_{\bs{l}-1} ) > \f{\de_{\bs{l}-1}}{2}, \; \bs{l} > \tau_\vep \ri \}, \\
&  & \mscr{E}_3 = \li \{ \lim_{\vep \downarrow 0} \mbf{n} = \iy, \; \lim_{\vep \downarrow 0} Z_{\bs{l}} = \mu, \; \lim_{\vep \downarrow 0}
Z_{\bs{l} - 1} = \mu \ri
\}, \\
&  & \mscr{E}_4 = \mscr{E}_1 \cap \mscr{E}_2 \cap  \mscr{E}_3, \\
&  &  \mscr{E}_5 =  \li \{ \liminf_{\vep \downarrow 0} \Lm_{\bs{l}} (\mu) \geq 1 \ri \},\\
&  &  \mscr{E}_6 =  \li \{ \limsup_{\vep \downarrow 0} \Lm_{\bs{l} - 1} (\mu) \leq 1, \;  \liminf_{\vep \downarrow 0} \bs{l} > 1 \; \tx{or} \;
\lim_{\vep \downarrow 0} \bs{l}  = 1 \ri \}. \eee} By the assumption that $n_{\tau_\vep} \to \iy$ as $\vep \downarrow 0$, we have that $\{
\lim_{\vep \downarrow 0} \mbf{n} = \iy \}$ is a sure event. It follows from the strong law of large numbers that $\mscr{E}_3$ is an almost sure
event. By the definition of the stopping rule, we have that  $\mscr{E}_1 \cap \mscr{E}_2$ is an almost sure event. Hence, $ \mscr{E}_4$ is an
almost sure event.  To show (\ref{cite369a}), it suffices to show that $\mscr{E}_4 \subseteq \mscr{E}_5$. For this purpose, we let $\om \in
\mscr{E}_4$ and attempt to show that $\om \in \mscr{E}_5$.  For simplicity of notations, let $\ell = \bs{l} (\om), \;  z_\ell = Z_{\bs{l}}
(\om), \; z_{\ell-1} = Z_{\bs{l} - 1} (\om)$,
\[
L_\ell = \mscr{L} ( z_\ell, \vep), \qu U_\ell = \mscr{U} ( z_\ell, \vep), \qu L_{\ell-1} = \mscr{L} ( z_{\ell-1}, \vep), \qu U_{\ell-1} =
\mscr{U} ( z_{\ell-1}, \vep).
\]
For small $\vep > 0$, $z_\ell$ and $z_{\ell - 1}$ will be bounded in a neighborhood of $\mu$.  Let $\eta > 0$ be small enough so that $\mu +
\eta, \; \mu - \eta \in \Se$.  Then, there exists $\vep^\star > 0$ such that
\[ | L_\ell - \mu | < \eta, \qu | U_\ell - \mu | < \eta, \qu | L_{\ell-1} - \mu | < \eta, \qu
| U_{\ell-1} - \mu | < \eta \] for all $\vep \in (0, \vep^\star)$.  It follows from the continuity of $\mcal{B} (.)$ that there exists a
constant $K > 0$ such that
\[
\mcal{B} (L_\ell) < K, \qu \mcal{B} (U_\ell) < K, \qu \mcal{B} (L_{\ell-1}) < K, \qu \mcal{B} (U_{\ell-1}) < K
\]
for all $\vep \in (0, \vep^\star)$.   As a consequence of $\om \in \mscr{E}_1$, we have
\[
F_{Z_\ell} (z_\ell, U_\ell ) \leq \f{\de_\ell}{2},  \qqu G_{Z_\ell} (z_\ell, L_\ell ) \leq \f{\de_\ell}{2}.
\]
Making use of these inequalities and (\ref{berry55}),  (\ref{berry56}) of Lemma \ref{nonuniformBE}, we have
\[
\Phi \li ( \f{ \sq{n_\ell} [  z_\ell - U_\ell  ]   }{  \sq{ \mcal{V} ( U_\ell ) } } \ri ) - \f{ \mcal{B} (  U_\ell ) }{\sq{n_\ell}} \leq
\f{\de_\ell}{2}, \qqu \Phi \li ( \f{ \sq{n_\ell} [ L_\ell - z_\ell  ]   }{  \sq{ \mcal{V} ( L_\ell ) } } \ri ) - \f{ \mcal{B} ( L_\ell )
}{\sq{n_\ell}}  \leq \f{\de_\ell}{2}
\]
from which we have
\[
\Lm_\ell (\mu) = \f{ \ka^2 (\mu) C_\ell }{ \mcal{V} (\mu) }  \geq \f{C_\ell}{ \vep^2 n_\ell} \li [ \Phi^{-1}   \li ( 1 - \f{\de_\ell}{2} - \f{
\mcal{B} ( U_\ell ) }{\sq{n_\ell}} \ri ) \ri ]^2  \li [ \f{ \ka (\mu) \vep }{  U_\ell - z_\ell  } \ri ]^2 \f{ \mcal{V} ( U_\ell ) } { \mcal{V}
(\mu) },
\]
\[
\Lm_\ell (\mu) = \f{ \ka^2 (\mu) C_\ell }{ \mcal{V} (\mu) }  \geq \f{C_\ell}{ \vep^2 n_\ell} \li [ \Phi^{-1}   \li ( 1 -  \f{\de_\ell}{2} - \f{
\mcal{B} ( L_\ell ) }{\sq{n_\ell}} \ri ) \ri ]^2  \li [ \f{ \ka (\mu) \vep }{  z_\ell - L_\ell } \ri ]^2  \f{ \mcal{V} ( L_\ell ) } { \mcal{V}
(\mu) }.
\]
These two inequalities imply that $\liminf_{\vep \downarrow 0} \Lm_\ell (\mu) \geq 1$. Hence, $\om \in \mscr{E}_5$ and thus (\ref{cite369a})
holds.

To show (\ref{cite369b}), it suffices to show that $\om \in \mscr{E}_6$ for $\om \in \mscr{E}_4$.  As a consequence of  $\om \in \mscr{E}_2$, we
have that either $\ell = \tau_\vep$ or
\[
F_{Z_{\ell-1}} (z_{\ell - 1}, U_{\ell-1} ) > \f{\de_{\ell - 1}}{2} \qu \tx{or} \qu  G_{Z_{\ell - 1}} (z_{\ell - 1}, L_{\ell-1} ) > \f{\de_{\ell
- 1}}{2}
\]
must be true for small enough $\vep \in (0, \vep^\star)$.  Making use of (\ref{berry55}) and (\ref{berry56}), we have that either $\ell =
\tau_\vep$ or
\[
\Phi \li ( \f{ \sq{{n_{\ell - 1}}} [  z_{\ell - 1} - U_{\ell-1}  ]   }{  \sq{ \mcal{V} ( U_{\ell-1} ) } } \ri ) + \f{\mcal{B} ( U_{\ell-1}  )
}{\sq{{n_{\ell - 1}}}} > \f{\de_{\ell-1}}{2} \qu \tx{or}
\]
\[
\Phi \li ( \f{ \sq{{n_{\ell - 1}}} [ L_{\ell-1} - z_{\ell - 1}  ]   }{  \sq{ \mcal{V} ( L_{\ell-1} ) } } \ri ) + \f{ \mcal{B} ( L_{\ell-1}  )
}{\sq{{n_{\ell - 1}}}} > \f{\de_{\ell-1}}{2}
\]
must be true for small enough $\vep \in (0, \vep^\star)$.   It follows that either $\ell = \tau_\vep$ or
\[
\Lm_{\ell - 1} (\mu) = \f{ \ka^2 (\mu) C_{\ell - 1} }{ \mcal{V} (\mu) }  < \f{C_{\ell - 1}}{ \vep^2 n_{\ell - 1}} \li [ \Phi^{-1}   \li ( 1 -
\f{\de_{\ell - 1}}{2} + \f{ \mcal{B} ( U_{\ell-1} ) }{\sq{n_{\ell - 1}}} \ri ) \ri ]^2  \li [ \f{ \ka (\mu) \vep }{  U_{\ell-1} - z_{\ell - 1} }
\ri ]^2 \f{ \mcal{V} ( U_{\ell-1} ) } { \mcal{V} (\mu) } \qqu \tx{or}
\]
\[
\Lm_{\ell - 1} (\mu) = \f{ \ka^2 (\mu) C_{\ell - 1} }{ \mcal{V} (\mu) }  < \f{C_{\ell - 1}}{ \vep^2 n_{\ell - 1}} \li [ \Phi^{-1}   \li ( 1 -
\f{\de_{\ell - 1}}{2} + \f{ \mcal{B} ( L_{\ell-1} ) }{\sq{n_{\ell - 1}}} \ri ) \ri ]^2  \li [ \f{ \ka (\mu) \vep }{  z_{\ell - 1} - L_{\ell-1} }
\ri ]^2 \f{ \mcal{V} ( L_{\ell-1} ) } { \mcal{V} (\mu) }
\]
must be true for small enough $\vep \in (0, \vep^\star)$.  This implies that either $\lim_{\vep \downarrow 0} \ell = 1$ or
\[
\limsup_{\vep \downarrow 0} \Lm_{\ell - 1} (\mu) \leq 1, \qqu \liminf_{\vep \downarrow 0} \ell > 1 \]
 must be true.  Hence, $\om \in \mscr{E}_6$ and thus (\ref{cite369b}) holds.

\epf

\subsection{Proof of Properties (III) -- (VII)}

Property (III) follows from Lemma \ref{lem22888} and the established Property (II).  A key result to establish Properties (IV) -- (VII) is Lemma
\ref{normboundary}, which will be developed in the sequel.  Once Lemma \ref{normboundary} is established,

Property (IV) follows from Lemmas \ref{lem20888}, \ref{lem21888} and \ref{normboundary};

Property (V) follows from Lemmas \ref{lem20888}, \ref{lem23888} and \ref{normboundary};

Property (VI) follows from Lemmas \ref{lem20888}, \ref{lem24888} and \ref{normboundary};

Property (VII) follows from Lemmas \ref{lem20888}, \ref{lem25888} and \ref{normboundary}.

\bsk

To establish Lemma \ref{normboundary}, we need the following result.

\beL  \la{lemvip98} Let $\vsi \in (0, 1)$ and $\eta \in (0, 1)$ such that  $\li ( \f{1 + \eta}{1 - \eta} \ri )^3 \f{1}{1 + \vsi} < 1 < \li (
\f{1 - \eta}{1 + \eta} \ri )^3  \f{1}{1 - \vsi}$. Let $\mscr{B}$ be a positive constant.  Then, there exists $\vep^* > 0$ such that \be
 \la{first8996}
 (1 - \eta) \f{C_\ell}{\vep^2}
 \li [ \Phi^{-1} \li (  1 - \f{\de_\ell}{2}   \ri ) \ri ]^2 < n_\ell < (1 + \eta) \f{C_\ell}{\vep^2} \li [ \Phi^{-1} \li (  1 - \f{\de_\ell}{2}   \ri ) \ri ]^2
 \ee
and \bel &  &  \li ( \f{1 + \eta}{1 - \eta} \ri )^3 \f{1}{1 + \vsi} \li [ \Phi^{-1} \li (  1 - \f{\de_\ell}{2}   \ri ) \ri ]^2 < \li [ \Phi^{-1}
 \li ( 1 - \f{\de_\ell}{2} - \f{ \mscr{B} }{n_\ell^2 \vep^3} \ri ) \ri ]^2, \la{note8}\\
&  & \li ( \f{1 - \eta}{1 + \eta} \ri )^3 \f{1}{1 - \vsi} \li [ \Phi^{-1} \li (  1 - \f{\de_\ell}{2}   \ri ) \ri ]^2 > \li [ \Phi^{-1}
 \li ( 1 - \f{\de_\ell}{2} + \f{ \mscr{B} }{n_\ell^2 \vep^3} \ri ) \ri ]^2 \la{note8b} \eel
for all $\vep \in (0, \vep^*)$ and $\ell \geq \tau_\vep$.

\eeL

\bpf

By the assumption that $\lim_{\vep \downarrow 0} \tau_\vep = 1$ and $\f{\vep^2 n_\ell}{C_\ell \Up_\ell} \to 1$ uniformly for $\ell \geq 1$ as
$\vep \downarrow 0$, we have that there exists $\vep_1 > 0$ such that $\tau_\vep = 1$ and (\ref{first8996}) holds for all $\vep \in (0, \vep_1)$
and $\ell \geq \tau_\vep$.  It follows that
\[
 \de_\ell n_\ell^2 \vep^4  > (1 - \eta)^2  \de_\ell ( C_\ell \Up_\ell )^2
\]
for $\ell \geq \tau_\vep = 1$ with $\vep \in (0, \vep_1)$.  By the assumption that $\de_\ell ( C_\ell \Up_\ell )^2 \to \iy$ as $\ell \to \iy$,
we have that there exists a constant $K > 0$ such that
\[
\f{ \mscr{B} }{ \de_\ell n_\ell^2 \vep^4 } < \f{ \mscr{B} }{ (1 - \eta)^2  \de_\ell ( C_\ell \Up_\ell )^2 } < \f{K}{2}
\]
for $\ell \geq \tau_\vep = 1$ with $\vep \in (0, \vep_1)$.  This implies that
\[
\f{ \mscr{B} }{ n_\ell^2 \vep^3 } < \f{\de_\ell}{2} K \vep
\]
for $\ell \geq \tau_\vep = 1$ with $\vep \in (0, \vep_2)$, where  \[ \vep_2 = \min \li \{ \vep_1, \f{1}{K} \ri \}.
\]
Therefore, \bel &  &  \li [ \Phi^{-1}
 \li ( 1 - \f{\de_\ell}{2} - \f{ \mscr{B} }{n_\ell^2 \vep^3 } \ri ) \ri ]^2 > \li [ \Phi^{-1}
 \li ( 1 - \f{\de_\ell}{2} - \f{\de_\ell}{2} K \vep \ri ) \ri ]^2, \la{inviewofa}\\
 &  & \li [ \Phi^{-1}
 \li ( 1 - \f{\de_\ell}{2} + \f{ \mscr{B} }{n_\ell^2 \vep^3 } \ri ) \ri ]^2 < \li [ \Phi^{-1}
 \li ( 1 - \f{\de_\ell}{2} + \f{\de_\ell}{2} K \vep \ri ) \ri ]^2 \la{inviewofb}
 \eel
for $\vep \in (0, \vep_2)$ and $\ell \geq \tau_\vep$.
 As a consequence of the fact that
 \be
 \la{fact89}
\lim_{\de  \downarrow 0} \f{ \li [ \Phi^{-1}
 \li ( 1 - \de  \ri ) \ri ]^2 } { 2 \ln \f{1}{\de}} = 1,
 \ee
there exists $\ba^* \in (0, 1)$ such that {\small  \bee &  & \f{ \li [ \Phi^{-1}
 \li ( 1 - \ba - \ba  K \vep  \ri ) \ri ]^2 } { \li [ \Phi^{-1}
 \li ( 1 - \ba  \ri ) \ri ]^2  }  >  \li ( \f{1 + \eta}{1 - \eta} \ri )^{\f{3}{2}}  \li ( \f{1}{1 + \vsi} \ri )^{\f{1}{2}} \; \f{ \ln \f{1}{ \ba (1 + K \vep) }   }
 {  \ln \f{1}{\ba}  }  > \li ( \f{1 + \eta}{1 - \eta} \ri )^{\f{3}{2}} \li ( \f{1}{1 + \vsi} \ri )^{\f{1}{2}} \; \li [ 1 -  \f{ \ln ( 1 + K \vep ) } { \ln \f{1}{\ba^*}  }  \ri
 ],\\
&  &  \f{ \li [ \Phi^{-1}
 \li ( 1 - \ba + \ba  K \vep  \ri ) \ri ]^2 } { \li [ \Phi^{-1}
 \li ( 1 - \ba  \ri ) \ri ]^2  } < \li ( \f{1 - \eta}{1 + \eta} \ri )^{\f{3}{2}} \li ( \f{1}{1 - \vsi} \ri )^{\f{1}{2}} \; \f{ \ln \f{1}{ \ba (1 - K \vep) }   } {  \ln \f{1}{\ba}  }
  <  \li ( \f{1 - \eta}{1 + \eta} \ri )^{\f{3}{2}}  \li ( \f{1}{1 - \vsi}  \ri )^{\f{1}{2}} \;  \li [ 1 -  \f{ \ln ( 1 - K \vep ) } { \ln \f{1}{\ba^*}  }  \ri ]
 \eee} for all $\ba \in (0, \ba^*)$ and $\vep \in (0, \vep_2)$. Hence, there exists $\vep_3 \in (0, \vep_2)$ such that \bee &  &  \f{ \li [ \Phi^{-1}
 \li ( 1 - \ba - \ba  K \vep  \ri ) \ri ]^2 } { \li [ \Phi^{-1}
 \li ( 1 - \ba  \ri ) \ri ]^2  } > \li ( \f{1 + \eta}{1 - \eta} \ri )^3 \f{1}{1 + \vsi},\\
 &  & \f{ \li [ \Phi^{-1}
 \li ( 1 - \ba + \ba  K \vep  \ri ) \ri ]^2 } { \li [ \Phi^{-1}
 \li ( 1 - \ba  \ri ) \ri ]^2  } < \li ( \f{1 - \eta}{1 + \eta} \ri )^3 \f{1}{1 - \vsi}
 \eee
for all $\ba \in (0, \ba^*)$ and $\vep \in (0, \vep_3)$.  Making use of (\ref{inviewofa}), (\ref{inviewofb}) and the above inequalities, for
$\ell \geq \tau_\vep$ with $\f{\de_\ell}{2} < \ba^*$, we have that (\ref{note8}) and (\ref{note8b}) hold for all $\vep \in (0, \vep_3)$.   On
the other hand, in view of (\ref{inviewofa}) and (\ref{inviewofb}), for $\ell \geq \tau_\vep$ with $\f{\de_\ell}{2} \geq \ba^*$, we can find
$\vep_4 \in (0, \vep_2)$ such that (\ref{note8}) and (\ref{note8b}) hold for all $\vep \in (0, \vep_4)$.  Thus, we have shown that there exists
$\vep^* \in (0, \vep_2)$ such that (\ref{note8}) and (\ref{note8b}) hold for all $\vep \in (0, \vep^*)$ and $\ell \geq \tau_\vep$. This
completes the proof of the lemma.

\epf

\beL

\la{normboundary}

For $\mu \in \Se$, if $\Lm_j (\mu)  > 1$ for some index $j \geq 1$, then there exist $\ga > 0$ and $\vep^* > 0$ such that $\{ \bs{D}_\ell = 0 \}
\subseteq \{ | Z_\ell - \mu  | \geq \ga  \}  \qu \tx{for all} \; \; \ell \geq j$ and $\vep \in (0, \vep^*)$.

Similarly, for $\mu \in \Se$, if $\Lm_j (\mu)  < 1$ for some index $j \geq 1$, then there exist $\ga > 0$ and $\vep^* > 0$ such that $\{
\bs{D}_\ell = 1 \} \subseteq \{ | Z_\ell - \mu  | \geq \ga  \}  \qu \tx{for all} \; \; \ell \leq j$ and $\vep \in (0, \vep^*)$.

 \eeL

\bpf

Define {\small \bee g_\ell(z, \vep) = n_\ell \min \li \{  \f{ [z -  \mscr{U} (z, \vep) ]^2 } { \mcal{V} ( \mscr{U} (z, \vep) )  \li [ \Phi^{-1}
\li ( 1 - \f{\de_\ell}{2} + \f{C \mscr{W} ( \mscr{U} (z, \vep)  ) }{ n_\ell^2 | \mscr{U} (z, \vep) - z|^3 } \ri ) \ri ]^2  }, \; \f{ [z -
\mscr{L} (z, \vep) ]^2 } { \mcal{V} ( \mscr{L} (z, \vep) ) \li [ \Phi^{-1} \li ( 1 - \f{\de_\ell}{2} +
\f{C \mscr{W} ( \mscr{L} (z, \vep) ) }{n_\ell^2 | \mscr{L} (z, \vep) - z  |^3 } \ri ) \ri ]^2 }  \ri \}, &  & \\
h_\ell(z, \vep) =  n_\ell \max \li \{ \f{  [z - \mscr{U} (z, \vep) ]^2 }{ \mcal{V} ( \mscr{U} (z, \vep) ) \li [ \Phi^{-1} \li ( 1 -
\f{\de_\ell}{2} - \f{C \mscr{W} ( \mscr{U} (z, \vep) ) }{n_\ell^2 | \mscr{U} (z, \vep) - z  |^3} \ri ) \ri ]^2 }, \; \f{  [z - \mscr{L} (z,
\vep) ]^2 }{ \mcal{V} ( \mscr{L} (z, \vep) ) \li [ \Phi^{-1} \li ( 1 - \f{\de_\ell}{2} - \f{C \mscr{W} ( \mscr{L} (z, \vep) ) }{n_\ell^2 |
\mscr{L} (z, \vep) - z  |^3} \ri ) \ri ]^2 } \ri \} &   &
 \eee} for $\ell \geq \tau_\vep$.  By virtue of nonuniform Berry-Essen inequality,
 it can be readily shown that {\small \bel &  &  \li \{ F_{Z_\ell} ( Z_\ell,  \mscr{U} (Z_\ell,
\vep) ) \leq \f{\de_\ell}{2},  \; \mscr{U}
(Z_\ell, \vep) \in \Se \ri \} \nonumber\\
 &  &  \supseteq
\li \{ \Phi \li ( \f{ \sq{{n_\ell}} [ Z_\ell - \mscr{U} (Z_\ell, \vep) ] }{ \sq{ \mcal{V} \li ( \mscr{U} (Z_\ell, \vep) \ri ) } } \ri ) \leq
\f{\de_\ell}{2} - \f{C \mscr{W}  ( \mscr{U} (Z_\ell, \vep)  ) }{n_\ell^2 |  \mscr{U} (Z_\ell, \vep) - Z_\ell  |^3},
\; \mscr{U} (Z_\ell, \vep) \in \Se \ri \}, \la{chaina} \\
&  &  \li \{ F_{Z_\ell} ( Z_\ell,  \mscr{U} (Z_\ell, \vep) ) > \f{\de_\ell}{2}, \; \mscr{U} (Z_\ell, \vep) \in \Se
\ri \} \nonumber\\
&  & \supseteq \li \{  \Phi \li ( \f{ \sq{{n_\ell}} [ Z_\ell - \mscr{U} (Z_\ell, \vep) ] }{ \sq{ \mcal{V} \li ( \mscr{U} (Z_\ell, \vep) \ri ) }
} \ri ) > \f{\de_\ell}{2} + \f{C \mscr{W}  ( \mscr{U} (Z_\ell, \vep)  ) }
{n_\ell^2 |  \mscr{U} (Z_\ell, \vep) - Z_\ell  |^3}, \; \mscr{U} (Z_\ell, \vep) \in \Se  \ri \}, \la{chainb} \\
&  & \li \{ G_{Z_\ell} ( Z_\ell,  \mscr{L} (Z_\ell, \vep) ) \leq \f{\de_\ell}{2}, \; \mscr{L} (Z_\ell, \vep) \in \Se
\ri \} \nonumber\\
&   & \supseteq \li \{  \Phi \li ( \f{ \sq{{n_\ell}} [ \mscr{L} (Z_\ell, \vep) - Z_\ell ] }{ \sq{ \mcal{V} \li ( \mscr{L} (Z_\ell, \vep)   \ri )
} } \ri ) \leq \f{\de_\ell}{2} - \f{C \mscr{W}  ( \mscr{L} (Z_\ell, \vep)  ) }{n_\ell^2 |  \mscr{L} (Z_\ell, \vep) - Z_\ell  |^3}, \;
\mscr{L} (Z_\ell, \vep) \in \Se  \ri \}, \la{chainc} \\
&  &  \li \{ G_{Z_\ell} ( Z_\ell,  \mscr{L} (Z_\ell, \vep) ) > \f{\de_\ell}{2}, \; \mscr{L} (Z_\ell, \vep) \in \Se
\ri \} \nonumber\\
&  & \supseteq \li \{  \Phi \li ( \f{ \sq{{n_\ell}} [ \mscr{L} (Z_\ell, \vep) - Z_\ell ] }{ \sq{ \mcal{V} \li ( \mscr{L} (Z_\ell, \vep)  \ri ) }
} \ri ) > \f{\de_\ell}{2} + \f{C \mscr{W}  ( \mscr{L} (Z_\ell, \vep)  ) }{n_\ell^2 |  \mscr{L} (Z_\ell, \vep) - Z_\ell  |^3}, \; \mscr{L}
(Z_\ell, \vep) \in \Se \ri \} \la{chaind} \eel} for $\ell \in \bb{Z}$ no less than $\tau_\vep$, where $C$ is an absolute constant.  Hence,
applying (\ref{chaina}), (\ref{chainb}), (\ref{chainc}), (\ref{chaind}), the definition of $g_\ell(z, \vep), \; h_\ell(z, \vep)$,  and the
definition of decision variable $\bs{D}_\ell$ for $\ell \geq \tau_\vep$, we have
 \be
 \la{homea}
 \{ g_\ell(
Z_\ell, \vep) \geq 1, \; \mscr{L} (Z_\ell, \vep) \in \Se,  \; \mscr{U} (Z_\ell, \vep) \in \Se \} \subseteq \{ \bs{D}_\ell = 1 \} \ee and \be
\la{homeb}
 \{ h_\ell(
Z_\ell, \vep) < 1, \; \mscr{L} (Z_\ell, \vep) \in \Se,  \; \mscr{U} (Z_\ell, \vep) \in \Se \} \subseteq \{ \bs{D}_\ell = 0 \} \ee
 for $\ell \geq \tau_\vep$.

Next, we claim that, for any $\se \in \Se$ and $\vsi \in (0, 1)$, there exist $\ga > 0$ and $\vep^* > 0$ such that
\[
\mscr{L} (z, \vep) \in \Se, \qqu \mscr{U} (z, \vep) \in \Se \] and that $(1 - \vsi) \Lm_\ell (z) < g_\ell(z, \vep) < h_\ell(z, \vep) < (1 +
\vsi) \Lm_\ell (z)$ for all $z \in (\se - \ga, \se + \ga) \subseteq \Se, \; \vep \in (0, \vep^*)$ and $\ell \geq \tau_\vep$.  To show this
claim, note that for $\vsi \in (0, 1)$, there exists $\eta \in (0, 1)$ such that  $\li ( \f{1 + \eta}{1 - \eta} \ri )^3 \f{1}{1 + \vsi} < 1 <
\li ( \f{1 - \eta}{1 + \eta} \ri )^3  \f{1}{1 - \vsi}$. As a consequence of $\ka(\se) > 0$, the assumption associated with (\ref{uniformasp})
and the continuity of $\ka(z)$ and $\mcal{V}(z)$, we have that there exist $\ga > 0$ and $\vep_1 > 0$ such that $\se + 2 \ga, \; \se - 2 \ga \in
\Se$,
\[
\mscr{L} (z, \vep) \in \Se, \qqu \mscr{U} (z, \vep) \in \Se
\]
and that
 \bel &   & (1 - \eta ) \ka(\se) \vep < | z - \mscr{U} (z, \vep) | < (1 + \eta ) \ka(\se) \vep, \la{6691}\\
 &  &  (1 - \eta ) \ka(\se) \vep < | z - \mscr{L} (z, \vep) |
< (1 + \eta ) \ka(\se) \vep, \la{6692} \\
&  & (1 - \eta)\ka(\se) < \ka(z) < (1 + \eta)\ka(\se), \la{6693} \\
&  & \se - 2 \ga < \se - \ga - (1 + \eta ) \ka(\se) \vep \leq \mscr{L} (z, \vep) < \mscr{U} (z, \vep) \leq \se  + \ga + (1 + \eta ) \ka(\se) \vep  < \se + 2 \ga \la{6694} \\
&  & (1 - \eta) \mcal{V} (z) < \mcal{V} ( \mscr{U} (z, \vep) ) < (1 + \eta) \mcal{V} (z), \la{addnewa} \\
&  & (1 - \eta) \mcal{V} (z) < \mcal{V} ( \mscr{L} (z, \vep) ) < (1 + \eta) \mcal{V} (z), \la{addnewb}
 \eel
 for all $z \in (\se - \ga, \se  + \ga)$ and $\vep \in (0, \vep_1)$.
 Define
 \[
 \mscr{B} = \f{ C }{  [ (1 - \eta) \ka (\se) ]^3 } \times \sup_{z \in (\se - 2 \ga, \se + 2 \ga)} \mscr{W} (z).
 \]
In view of (\ref{6691}), (\ref{6692}) and (\ref{6694}), we have \be \la{see88a}
 \f{C \mscr{W} ( \mscr{L} (z, \vep) ) }{n_\ell^2 | \mscr{L} (z,
\vep) - z |^3} \leq \f{ \mscr{B} }{n_\ell^2 \vep^3} \ee and \be \la{see88b}
 \f{C \mscr{W} ( \mscr{U} (z, \vep) ) }{n_\ell^2 | \mscr{U} (z, \vep)
- z |^3} \leq \f{ \mscr{B} }{n_\ell^2 \vep^3} \ee for all $z \in (\se - \ga, \se  + \ga)$ and $\vep \in (0, \vep_1)$.  According to Lemma
\ref{lemvip98},  there exists $\vep^* \in (0, \vep_1)$ such that (\ref{first8996}), (\ref{note8}) and (\ref{note8b}) hold for all $\vep \in (0,
\vep^*)$ and $\ell \geq \tau_\vep$. Since (\ref{note8}) and (\ref{note8b}) can be, respectively, written as \bee &   & \f{\f{ (1 + \eta)
}{\vep^2} C_\ell \li [ \Phi^{-1} \li (  1 - \f{\de_\ell}{2} \ri ) \ri ]^2} {(1 - \eta)  [ (1 - \eta)\ka(\se) ]^2 C_\ell} \times \f{ [ (1 + \eta
) \ka(\se) \vep ]^2 }{ \li [ \Phi^{-1} \li ( 1 - \f{\de_\ell}{2} - \f{\mscr{B} }{n_\ell^2 \vep^3} \ri ) \ri ]^2 }  < 1 + \vsi,\\
&  &  \f{\f{ (1 - \eta) }{\vep^2} C_\ell \li [ \Phi^{-1} \li (  1 - \f{\de_\ell}{2}   \ri ) \ri ]^2} {(1 + \eta) [ (1 + \eta)\ka(\se) ]^2
C_\ell} \times \f{ [ (1 - \eta ) \ka(\se) \vep ]^2 } { \li [ \Phi^{-1}  \li ( 1 - \f{\de_\ell}{2} + \f{\mscr{B}}{n_\ell^2 \vep^3} \ri ) \ri ]^2
} > 1 - \vsi, \eee it follows from (\ref{first8996}) that \bel &  & \f{n_\ell} {(1 - \eta)  [ (1 - \eta)\ka(\se) ]^2 C_\ell} \times \f{ [ (1 +
\eta ) \ka(\se) \vep ]^2 }{ \li [ \Phi^{-1} \li ( 1 - \f{\de_\ell}{2} -
\f{\mscr{B}}{n_\ell^2 \vep^3} \ri ) \ri ]^2 } < 1 + \vsi, \la{6695} \\
&  &  \f{ n_\ell} {(1 + \eta)  [ (1 + \eta)\ka(\se) ]^2 C_\ell} \times  \f{ [ (1 - \eta ) \ka(\se) \vep ]^2 } { \li [ \Phi^{-1}  \li ( 1 -
\f{\de_\ell}{2} + \f{\mscr{B}}{n_\ell^2 \vep^3} \ri ) \ri ]^2 }  > 1 - \vsi \la{6696} \eel for all $\vep \in (0, \vep^*)$ and $\ell \geq
\tau_\vep$. Therefore, \bel &  & \f{n_\ell}{(1 - \eta)  [ \ka(z) ]^2 C_\ell} \times  \f{ [z - \mscr{U} (z, \vep) ]^2 }{ \li [ \Phi^{-1} \li ( 1
- \f{\de_\ell}{2}
- \f{C \mscr{W} ( \mscr{U} (z, \vep) ) }{n_\ell^2 | \mscr{U} (z, \vep) - z  |^3} \ri ) \ri ]^2 } < 1 + \vsi, \la{6697}\\
&  & \f{ n_\ell} { (1 - \eta)  [ \ka(z) ]^2 C_\ell} \times  \f{ [z - \mscr{L} (z, \vep) ]^2 } { \li [ \Phi^{-1}  \li ( 1 - \f{\de_\ell}{2} -
\f{C \mscr{W} (
\mscr{L} (z, \vep) ) }{n_\ell^2 | \mscr{L} (z, \vep) - z  |^3} \ri ) \ri ]^2 } < 1 + \vsi, \la{6698}\\
&  & \f{n_\ell} {(1 + \eta)  [ \ka(z) ]^2 C_\ell} \times  \f{ [z - \mscr{U} (z, \vep) ]^2 }{ \li [ \Phi^{-1}  \li ( 1 - \f{\de_\ell}{2} + \f{C
\mscr{W} (
\mscr{U} (z, \vep) ) }{n_\ell^2 | \mscr{U} (z, \vep) - z  |^3} \ri ) \ri ]^2 } > 1 - \vsi, \la{6699}\\
&  & \f{ n_\ell}{(1 + \eta)  [ \ka(z) ]^2 C_\ell} \times  \f{ [z - \mscr{L} (z, \vep) ]^2 } { \li [ \Phi^{-1}  \li ( 1 - \f{\de_\ell}{2} + \f{C
\mscr{W} ( \mscr{L} (z, \vep) ) }{n_\ell^2 | \mscr{L} (z, \vep) - z  |^3} \ri ) \ri ]^2 } > 1 - \vsi \la{66910} \eel for all $z \in (\se - \ga,
\se + \ga), \; \vep \in (0, \vep^*)$ and $\ell \geq \tau_\vep$.  Here, equation (\ref{6697}) follows from (\ref{see88b}) and (\ref{6695}).
Equation (\ref{6698}) follows from (\ref{see88a}) and (\ref{6695}). Equation (\ref{6699}) follows from (\ref{see88b}) and (\ref{6696}). Equation
(\ref{66910}) follows from (\ref{see88a}) and (\ref{6696}).  It follows from (\ref{addnewa}), (\ref{addnewb}), (\ref{6697}), (\ref{6698}),
(\ref{6699}), (\ref{66910}) and the definition of $\Lm_\ell (z)$ that \bel &  & \f{n_\ell}{\mcal{V} ( \mscr{U} (z, \vep) )} \times \f{ [z -
\mscr{U} (z, \vep) ]^2  } { \li [ \Phi^{-1}
\li ( 1 - \f{\de_\ell}{2} - \f{C \mscr{W} ( \mscr{U} (z, \vep)  ) }{n_\ell^2 | \mscr{U} (z, \vep) - z  |^3} \ri ) \ri ]^2 } < (1 + \vsi) \Lm_\ell (z), \la{ineq99a} \\
&  &  \f{n_\ell}{\mcal{V} ( \mscr{L} (z, \vep) )} \times  \f { [z - \mscr{L} (z, \vep) ]^2  }{ \li [ \Phi^{-1}  \li ( 1 - \f{\de_\ell}{2} -
\f{C \mscr{W} ( \mscr{L} (z, \vep)  ) }{n_\ell^2 | \mscr{L} (z, \vep) - z  |^3} \ri ) \ri ]^2 } < (1 + \vsi) \Lm_\ell (z), \la{ineq99b} \\
&  &  \f{n_\ell}{\mcal{V} ( \mscr{U} (z, \vep) )} \times  \f{ [z - \mscr{U} (z, \vep) ]^2  }{ \li [ \Phi^{-1}  \li ( 1 - \f{\de_\ell}{2} +
\f{C \mscr{W} ( \mscr{U} (z, \vep)  ) }{n_\ell^2 | \mscr{U} (z, \vep) - z  |^3} \ri ) \ri ]^2 } > (1 - \vsi) \Lm_\ell (z), \la{ineq99c} \\
&  &   \f{n_\ell}{\mcal{V} ( \mscr{L} (z, \vep) )} \times  \f{ [z - \mscr{L} (z, \vep) ]^2  }{ \li [ \Phi^{-1}  \li ( 1 - \f{\de_\ell}{2} + \f{C
\mscr{W} ( \mscr{L} (z, \vep)  ) }{n_\ell^2 | \mscr{L} (z, \vep) - z  |^3} \ri ) \ri ]^2 } > (1 - \vsi) \Lm_\ell (z) \la{ineq99d} \eel for all
$z \in (\se - \ga, \se + \ga), \; \vep \in (0, \vep^*)$ and $\ell \geq \tau_\vep$.  Making use of (\ref{ineq99a}), (\ref{ineq99b}),
(\ref{ineq99c}), (\ref{ineq99d}) and the definitions of $g_\ell (z, \vep)$ and $h_\ell (z, \vep)$, we have that for any $\se \in \Se$ and $\vsi
\in (0, 1)$, there exist $\ga > 0$ and $\vep^* > 0$ such that $\mscr{L} (z, \vep) \in \Se, \; \mscr{U} (z, \vep) \in \Se$ and that $(1 - \vsi)
\Lm_\ell (z) < g_\ell(z, \vep) < h_\ell(z, \vep) < (1 + \vsi) \Lm_\ell (z)$ for all $z \in (\se - \ga, \se + \ga), \; \vep \in (0, \vep^*)$ and
$\ell \geq \tau_\vep$.  Our claim is thus proved.

For $\mu \in \Se$ with $\Lm_j (\mu) > 1$, let {\small $\vsi \in \li ( 0, 1 - \f{1}{ \sq{ \Lm_j (\mu) } } \ri )$}. By the continuity of $\Lm_\ell
(z)$ on $z \in \Se$, there exists a number $\ga_0 > 0$ such that
\[
\Lm_j (z) \geq (1 - \vsi) \Lm_j (\mu) \qu \tx{for all} \; z \in ( \mu - \ga_0, \mu + \ga_0 ).
\]
Since $\Lm_\ell(\mu)$ is increasing with respect to $\ell$, it must be true that
\[
\Lm_\ell (z) \geq (1 - \vsi) \Lm_j (\mu) \qu \tx{for all} \; \; \ell \geq j  \; \tx{and} \; z \in ( \mu - \ga_0, \mu + \ga_0 ).
\]
By our established claim,  there exist $\ga \in (0, \ga_0)$ and $\vep^* > 0$ such that
\[
g_\ell(z, \vep) > (1 - \vsi) \Lm_\ell (z) \geq (1 - \vsi)^2 \Lm_j (\mu) > 1, \qqu \mscr{L} (z, \vep) \in \Se, \qqu \mscr{U} (z, \vep) \in \Se
\]
for all $z \in (\mu - \ga, \mu + \ga), \; \vep \in (0, \vep^*)$ and $\ell \geq j$.  Making use of this result and (\ref{homea}), we have that
there exist $\ga > 0$ and $\vep^* > 0$ such that $\{ | Z_\ell - \mu  | < \ga \} \subseteq  \{ \bs{D}_\ell = 1  \}$ for all $\vep \in (0,
\vep^*)$ and $\ell \geq j$.  This implies that, there exist $\ga > 0$ and $\vep^* > 0$ such that $\{ \bs{D}_\ell = 0  \} \subseteq \{ | Z_\ell -
\mu  | \geq \ga \}$ for all $\vep \in (0, \vep^*)$ and $\ell \geq j$. Hence, we have shown the first statement of the lemma.

To show the second statement of the lemma, choose {\small $\vsi \in \li (0, \f{1}{ \sq{ \Lm_j (\mu) } } - 1 \ri )$} for $\mu \in \Se$ with
$\Lm_j (\mu) < 1$.  By the continuity of $\Lm_\ell (z)$ on $z \in \Se$, there exists a number $\ga_0 > 0$ such that
\[
\Lm_j (z) \leq (1 + \vsi) \Lm_j (\mu) \qu \tx{for all} \; z \in ( \mu - \ga_0, \mu + \ga_0 ).
\]
Since $\Lm_\ell(\mu)$ is increasing with respect to $\ell$, it must be true that
\[
\Lm_\ell (z) \leq (1 + \vsi) \Lm_j (\mu) \qu \tx{for all} \; \; \ell \leq j  \; \tx{and} \; z \in ( \mu - \ga_0, \mu + \ga_0 ).
\]
By our established claim,  there exist $\ga \in (0, \ga_0)$ and $\vep^* > 0$ such that
\[
h_\ell(z, \vep) < (1 + \vsi) \Lm_\ell (z) \leq (1 + \vsi)^2 \Lm_j (\mu) < 1, \qqu \mscr{L} (z, \vep) \in \Se, \qqu \mscr{U} (z, \vep) \in \Se
\]
for all $z \in (\mu - \ga, \mu + \ga), \; \vep \in (0, \vep^*)$ and $\ell \leq j$. Making use of this result and (\ref{homeb}), we have that
there exist $\ga > 0$ and $\vep^* > 0$ such that $\{ | Z_\ell - \mu  | < \ga \} \subseteq  \{ \bs{D}_\ell = 0  \}$ for all $\vep \in (0,
\vep^*)$ and $\ell \leq j$.  This implies that, there exist $\ga > 0$ and $\vep^* > 0$ such that $\{ \bs{D}_\ell = 1  \} \subseteq \{ | Z_\ell -
\mu  | \geq \ga \}$ for all $\vep \in (0, \vep^*)$ and $\ell \leq j$. So, we have shown the second statement of the lemma. This completes the
proof of the lemma.

\epf

\sect{Proof of Theorem \ref{CHST}} \la{CHST_app}

The proof of (\ref{VIP886d}) is similar to that of (\ref{VIP886b}).  The other properties are proved in the sequel.

\subsection{Proof of Property (I)}

To establish (\ref{perty891}), we can make use of a similar method as that of the counterpart of Theorem \ref{CHCDF} and Lemma
\ref{asn889222group} as follows.

\beL \la{asn889222group} There exist a real number $\ga > 0$ and an integer $m > \tau_\vep$ such that \be \la{makeuse222}
 \{ \bs{l} > \ell \} \subseteq  \{ | Z_\ell - \mu | \geq \ga \} \qqu \tx{for all $\ell > m$}. \ee \eeL

\bpf From Lemma \ref{remove asp}, we have that there exists a number $\ga > 0$ such that $(\mu - \ga, \mu + \ga) \subseteq \Se$ and that \[
\sup_{\se \in (\mu - \ga, \mu + \ga) } \mscr{M} (\se, \mscr{U} ( \se, \vep) ) < 0.
\]
Since $C_\ell \to \iy$ as $\ell \to \iy$, it follows that  $\f{ \ln \f{\de_\ell}{2} }{n_\ell} \to 0$ as $\ell \to \iy$.  Hence, there exists $m
> \tau_\vep$ such that
\[
\f{ \ln \f{\de_\ell}{2} }{n_\ell} > \sup_{\se \in (\mu - \ga, \mu + \ga) } \mscr{M} (\se, \mscr{U} ( \se, \vep) ) \qu \tx{for $\ell > m$}.
\]
This implies that
\[
\li \{ \mscr{M} ( Z_\ell, \mscr{U} ( Z_\ell, \vep) )  > \f{ \ln \f{\de_\ell}{2} }{n_\ell}, \; | Z_\ell - \mu | < \ga \ri \} = \emptyset \qu
\tx{for all $\ell > m$}.
\]
It follows that {\small \bee  \li \{ \mscr{M} ( Z_\ell, \mscr{U} ( Z_\ell, \vep) )  > \f{ \ln \f{\de_\ell}{2} }{n_\ell} \ri \} & \subseteq & \li
\{ \mscr{M} ( Z_\ell, \mscr{U} ( Z_\ell, \vep) )  > \f{ \ln \f{\de_\ell}{2} }{n_\ell}, \; | Z_\ell - \mu | < \ga \ri
\} \cup \{ | Z_\ell - \mu | \geq \ga \}\\
&  = & \{ | Z_\ell - \mu | \geq \ga \} \eee} for all $\ell > m$.  So, we have shown that there exist a real number $\ga
> 0$ and an integer $m > \tau_\vep$ such that $\li \{ \mscr{M} ( Z_\ell, \mscr{U} ( Z_\ell, \vep) )  > \f{ \ln \f{\de_\ell}{2} }{n_\ell} \ri \}  \subseteq \{ | Z_\ell - \mu |
\geq \ga \}$ for all $\ell > m$.  In a similar manner, we can show that there exist a real number $\ga
> 0$ and an integer $m > \tau_\vep$ such that
\[
\li \{ \mscr{M} ( Z_\ell, \mscr{L} ( Z_\ell, \vep) )  > \f{ \ln \f{\de_\ell}{2} }{n_\ell} \ri \} \subseteq  \{ | Z_\ell - \mu | \geq \ga \} \qu
\tx{for all $\ell > m$}.
\]
Therefore, there exist a real number $\ga > 0$ and an integer $m > \tau_\vep$ such that {\small \[ \li \{ \mscr{M} ( Z_\ell, \mscr{U} (Z_\ell,
\vep) ) \leq \f{\ln \f{\de_\ell}{2} }{n_\ell}, \; \mscr{M} ( Z_\ell, \mscr{L} (Z_\ell, \vep) ) \leq \f{ \ln \f{\de_\ell}{2} }{n_\ell} \ri \}
\supseteq   \{ | Z_\ell - \mu | < \ga \}
\]}
for all $\ell > m$.  Making use of this result and the definition of the stopping rule, we have that {\small $\{ \bs{l} > \ell \} \subseteq \{ |
Z_\ell - \mu | \geq \ga \}$} for all $\ell > m$.  This completes the proof of the lemma. \epf

\subsection{Proof of Property (II)}

To prove (\ref{propertyii1}) and (\ref{propertyii2}), we can apply Lemma \ref{pro689a} and Lemma \ref{pro689bch} as follows.

\beL

\la{pro689bch}

\bel &  & \Pr \li \{ \liminf_{\vep \downarrow 0} \Lm_{\bs{l}} (\mu) \geq 1 \ri \} = 1, \la{cite369ach}\\
&  & \Pr \li \{ \tx{Both} \; \limsup_{\vep \downarrow 0} \Lm_{\bs{l} - 1} (\mu) \leq 1 \; \tx{and} \;  \liminf_{\vep \downarrow 0} \bs{l} > 1 \;
\tx{hold or} \; \lim_{\vep \downarrow 0} \bs{l}  = 1 \ri \} = 1.  \la{cite369bch} \eel

\eeL

\bpf

To simplify notations, define
\[
\bs{L}_{\bs{l}} = \mscr{L} ( Z_{\bs{l}}, \vep ), \qqu \bs{L}_{\bs{l} - 1} = \mscr{L} ( Z_{\bs{l} - 1}, \vep ), \qqu \bs{U}_{\bs{l}} = \mscr{U} (
Z_{\bs{l}}, \vep ), \qqu \bs{U}_{\bs{l} - 1} = \mscr{U} ( Z_{\bs{l} - 1}, \vep ).
\]
Define {\small \bee &  & \mscr{E}_1 = \li \{ \mscr{M} (Z_{\bs{l}}, \bs{U}_{\bs{l}} ) \leq \f{\ln \f{\de_{\bs{l}}}{2} }{n_{\bs{l}} }, \;
\mscr{M} (Z_{\bs{l}}, \bs{L}_{\bs{l}} ) \leq \f{\ln \f{\de_{\bs{l}}}{2} }{n_{\bs{l}}} \ri \}, \\
&  & \mscr{E}_2 = \{  \bs{l} = \tau_\vep  \} \bigcup \li \{  \mscr{M} (Z_{\bs{l}-1}, \bs{U}_{\bs{l}-1} ) > \f{\ln
\f{\de_{\bs{l}-1}}{2}}{n_{\bs{l}-1}}, \; \bs{l} > \tau_\vep \ri \}
 \bigcup \li \{ \mscr{M} (Z_{\bs{l}-1}, \bs{L}_{\bs{l}-1} ) > \f{\ln \f{\de_{\bs{l}-1}}{2}}{n_{\bs{l}-1}}, \;  \bs{l} > \tau_\vep \ri \}, \\
&  & \mscr{E}_3 = \li \{ \lim_{\vep \downarrow 0} \mbf{n} = \iy, \; \lim_{\vep \downarrow 0} Z_{\bs{l}} = \mu, \; \lim_{\vep \downarrow 0}
Z_{\bs{l} - 1} = \mu \ri
\}, \\
&  & \mscr{E}_4 = \mscr{E}_1 \cap \mscr{E}_2 \cap  \mscr{E}_3, \\
&  &  \mscr{E}_5 =  \li \{ \liminf_{\vep \downarrow 0} \Lm_{\bs{l}} (\mu) \geq 1 \ri \},\\
&  &  \mscr{E}_6 =  \li \{ \tx{Both} \; \limsup_{\vep \downarrow 0} \Lm_{\bs{l} - 1} (\mu) \leq 1 \; \tx{and} \; \liminf_{\vep \downarrow 0}
\bs{l} > 1 \; \tx{hold or} \; \lim_{\vep \downarrow 0} \bs{l}  = 1 \ri \}. \eee} By the assumption that $n_{\tau_\vep} \to \iy$ as $\vep
\downarrow 0$, we have that $\{ \lim_{\vep \downarrow 0} \mbf{n} = \iy \}$ is a sure event. It follows from the strong law of large numbers that
$\mscr{E}_3$ is an almost sure event. By the definition of the stopping rule, we have that  $\mscr{E}_1 \cap \mscr{E}_2$ is an almost sure
event. Hence, $ \mscr{E}_4$ is an almost sure event.  To show (\ref{cite369ach}), it suffices to show that $\mscr{E}_4 \subseteq \mscr{E}_5$.
For this purpose, we let $\om \in \mscr{E}_4$ and attempt to show that $\om \in \mscr{E}_5$.  For simplicity of notations, let $\ell = \bs{l}
(\om), \;  z_\ell = Z_{\bs{l}} (\om), \; z_{\ell-1} = Z_{\bs{l} - 1} (\om)$,
\[
L_\ell = \mscr{L} ( z_\ell, \vep), \qu U_\ell = \mscr{U} ( z_\ell, \vep), \qu L_{\ell-1} = \mscr{L} ( z_{\ell-1}, \vep), \qu U_{\ell-1} =
\mscr{U} ( z_{\ell-1}, \vep).
\]
For small $\vep
> 0$, $z_\ell$ and $z_{\ell - 1}$ will be bounded in a neighborhood of $\mu$. We restrict $\vep > 0$ to be sufficiently small so that \[
L_\ell \in \Se, \qu U_\ell \in \Se, \qu L_{\ell-1} \in \Se, \qu U_{\ell-1} \in \Se.
\]
As a consequence of $\om \in \mscr{E}_1$,  \bel & & n_\ell \geq \f{ \ln (\f{\de_\ell}{2}) }{\mscr{M} (z_\ell, U_\ell
 )}, \la{imp18new} \\
&  & n_\ell \geq \f{ \ln (\f{\de_\ell}{2}) }{\mscr{M} (z_\ell, L_\ell
 )}.  \la{imp2new}
\eel  From (\ref{imp18new}), we have \bee \Lm_\ell (\mu) = \f{ \ka^2  (\mu) C_\ell }{ \mcal{V} (\mu) }
 \geq \f{C_\ell}{\vep^2 n_\ell} \f{ 2 \mcal{V} ( U_\ell ) \ln (\f{\de_\ell}{2}) } { [ U_\ell - z_\ell ]^2 } \f{  [ U_\ell - z_\ell ]^2 }
{ 2 \mcal{V} ( U_\ell ) \mscr{M} (z_\ell, U_\ell  ) }  \f{ [\ka  (\mu) \vep]^2 }{ \mcal{V} (\mu) }.  \eee Making use of this result and Lemma
\ref{uniformcon},  we have  \be \la{imp3new} \liminf_{\vep \downarrow 0} \Lm_\ell (\mu) \geq 1. \ee Similarly, (\ref{imp3new}) can be deduced
from (\ref{imp2new}) by virtue of Lemma \ref{uniformcon}.  Hence, $\om \in \mscr{E}_5$ and thus (\ref{cite369ach}) holds.

To show (\ref{cite369bch}), it suffices to show that $\om \in \mscr{E}_6$ for $\om \in \mscr{E}_4$.  Since $\om \in \mscr{E}_2$, we have that
either $\ell = \tau_\vep$ or \bel & & n_{\ell - 1}  < \f{ \ln (\f{\de_{\ell - 1}}{2}) }{\mscr{M}
(z_{\ell - 1}, U_{\ell-1} )}  \qqu \tx{or} \la{imp18newb} \\
&  & n_{\ell - 1}  < \f{ \ln (\f{\de_{\ell - 1}}{2}) }{\mscr{M} (z_{\ell - 1}, L_{\ell-1}
 )}  \la{imp2newb}
\eel must be true for small enough $\vep > 0$.  From (\ref{imp18newb}), we have \bee \Lm_{\ell - 1} (\mu)  =  \f{ \ka^2  (\mu) C_{\ell - 1} }{
\mcal{V} (\mu) }
 <  \f{C_{\ell - 1}}{\vep^2 n_{\ell - 1}} \f{ 2 \mcal{V} ( U_{\ell-1} ) \ln (\f{\de_{\ell - 1}}{2}) }
 { [ U_{\ell-1} - z_{\ell - 1} ]^2 } \f{  [ U_{\ell-1} - z_{\ell - 1} ]^2 }
{ 2 \mcal{V} ( U_{\ell-1} ) \mscr{M} (z_{\ell - 1}, U_{\ell-1}  ) }  \f{ [\ka  (\mu) \vep]^2 }{ \mcal{V} (\mu) }.  \eee Using this result and
Lemma \ref{uniformcon}, we have that \be \la{imp3newb} \limsup_{\vep \downarrow 0} \Lm_{\ell - 1} (\mu) \leq 1. \ee Similarly, (\ref{imp3newb})
can be deduced from (\ref{imp2newb}) by virtue of Lemma \ref{uniformcon}.  This implies that either $\lim_{\vep \downarrow 0} \ell = 1$ or
\[
\limsup_{\vep \downarrow 0} \Lm_{\ell - 1} (\mu) \leq 1, \qqu \liminf_{\vep \downarrow 0} \ell > 1 \]
 must be true.  Hence, $\om \in \mscr{E}_6$ and thus (\ref{cite369bch}) holds.

\epf

\subsection{Proof of Properties (III) -- (VII)}

Property (III) follows from Lemma \ref{lem22888} and the established Property (II).  A key result to establish Properties (IV) -- (VII) is Lemma
\ref{CHboundary}, which will be developed in the sequel.  Once Lemma \ref{CHboundary} is established,

Property (IV) follows from Lemmas \ref{lem20888}, \ref{lem21888} and \ref{CHboundary}.

Property (V) follows from Lemmas \ref{lem20888}, \ref{lem23888} and \ref{CHboundary}.

Property (VI) follows from Lemmas \ref{lem20888}, \ref{lem24888} and \ref{CHboundary}.

Property (VII) follows from Lemmas \ref{lem20888}, \ref{lem25888} and \ref{CHboundary}.

\beL

\la{CHboundary}

For $\mu \in \Se$, if $\Lm_j (\mu)  > 1$ for some index $j \geq 1$, then there exist a number $\ga > 0$ and $\vep^* > 0$ such that $\{
\bs{D}_\ell = 0 \} \subseteq \{ | Z_\ell - \mu  | \geq \ga  \}  \qu \tx{for all} \; \; \ell \geq j$ and $\vep \in (0, \vep^*)$.

Similarly, for $\mu \in \Se$, if $\Lm_j (\mu)  < 1$ for some index $j \geq 1$, then there exist a number $\ga > 0$ and $\vep^* > 0$ such that
$\{ \bs{D}_\ell = 1 \} \subseteq \{ | Z_\ell - \mu  | \geq \ga  \}  \qu \tx{for all} \; \; \ell \leq j$ and $\vep \in (0, \vep^*)$.

\eeL

\bpf

Define \be \la{defHell} H_\ell(z, \vep) = \f{n_\ell}{ \ln \f{\de_\ell}{2}  } \max \{ \mscr{M} ( z, \mscr{L} (z, \vep) ), \; \mscr{M} ( z,
\mscr{U} (z, \vep) ) \} \ee  for $\ell \geq \tau_\vep$. By the definitions of $H_\ell(z, \vep)$, the decision variable $\bs{D}_\ell$, and the
stopping rule, we have that \be
 \la{homed}
 \{ H_\ell(
Z_\ell, \vep) \geq 1, \; \mscr{L} (Z_\ell, \vep) \in \Se, \; \mscr{U} (Z_\ell, \vep) \in \Se \} \subseteq \{ \bs{D}_\ell = 1 \} \ee  and \be
\la{homec}
 \{
H_\ell( Z_\ell, \vep) < 1, \; \mscr{L} (Z_\ell, \vep) \in \Se, \; \mscr{U} (Z_\ell, \vep) \in \Se \} \subseteq \{ \bs{D}_\ell = 0 \} \ee for
$\ell \geq \tau_\vep$.

Next, we claim that for any $\se \in \Se$ and $\vsi \in (0, 1)$, there exist $\ga > 0$ and $\vep^* > 0$ such that $\mscr{L} (z, \vep) \in \Se,
\; \mscr{U} (z, \vep) \in \Se$ and that $(1 - \vsi) \Lm_\ell (z) < H_\ell(z, \vep) < (1 + \vsi) \Lm_\ell (z)$ for all $z \in (\se - \ga, \se +
\ga) \subseteq \Se, \; \vep \in (0, \vep^*)$ and $\ell \geq \tau_\vep$.  To prove this claim, note that for $\vsi \in (0, 1)$, there exists
$\eta \in (0, 1)$ such that  $\f{(1 + \eta)^4}{ 1 + \vsi } < 1 < \f{(1 - \eta)^4}{ 1 - \vsi }$. Note that there exist $\ga > 0$ and $\vep_1 > 0$
such that $\mscr{L} (z, \vep) \in \Se, \; \mscr{U} (z, \vep) \in \Se$ and that \bel &  &  1 - \eta < \f{ - 2 \mcal{V} (\mscr{L} (z, \vep))
\mscr{M} ( z, \mscr{L} (z, \vep) )}{[ \mscr{L} (z, \vep) - z ]^2 }  < 1 +
\eta, \la{inq33881} \\
&  &  1 - \eta < \f{- 2 \mcal{V} (\mscr{U} (z, \vep)) \mscr{M} ( z, \mscr{U} (z, \vep) )}{ [  \mscr{U} (z, \vep) - z ]^2}  < 1 + \eta,
\la{inq33882} \\
&  &  1 - \eta < \li [ \f{z - \mscr{L} (z, \vep)  }{  \ka (z) \vep } \ri ]^2 < 1 + \eta, \la{inq33883} \\
&  &  1 - \eta < \li [ \f{ \mscr{U} (z, \vep) - z }{  \ka (z) \vep } \ri ]^2 < 1 + \eta, \la{inq33884} \\
&  & 1 - \eta < \f{ \mcal{V} (z) } { \mcal{V} (\mscr{L} (z, \vep))  } < 1 + \eta,  \la{inq33885} \\
&  &  1 - \eta < \f{ \mcal{V} (z) } { \mcal{V} (\mscr{U} (z, \vep)) } < 1 + \eta \la{inq33886} \eel for all $z \in (\se - \ga, \se + \ga)$ and
$\vep \in (0, \vep_1)$. Here, inequalities (\ref{inq33881}) and (\ref{inq33882}) follow from Lemma \ref{uniformcon}.  Inequalities
(\ref{inq33883}) and (\ref{inq33884})  are due to the continuity of $\ka(z)$ and the assumption associated with (\ref{uniformasp}).
Inequalities (\ref{inq33885}) and (\ref{inq33886}) follows from the fact that $\f{ \mcal{V} (z) } { \mcal{V} (\mscr{L} (z, \vep))  }$ and $\f{
\mcal{V} (z) } { \mcal{V} (\mscr{U} (z, \vep))  }$ are continuous functions of $z$ and $\vep$ which tend to $1$ as $(z, \vep)$ tends to $(\se,
0)$.   Multiplying inequalities (\ref{inq33881}), (\ref{inq33883}) and (\ref{inq33885}) yields \be \la{hou1}
 (1 - \eta)^3 < \f{ - 2 \mcal{V} (z)  \mscr{M} ( z, \mscr{L} (z, \vep)
)}{ [\ka (z) \vep]^2 }  < (1 + \eta)^3 \ee for all $z \in (\se - \ga, \se + \ga)$ and $\vep \in (0, \vep_1)$. Multiplying inequalities
(\ref{inq33882}), (\ref{inq33884}) and (\ref{inq33886}) yields \be \la{hou2} (1 - \eta)^3 < \f{ - 2 \mcal{V} (z) \mscr{M} ( z, \mscr{U} (z,
\vep) )}{ [\ka (z) \vep]^2 }  < (1 + \eta)^3 \ee for all $z \in (\se - \ga, \se + \ga)$ and $\vep \in (0, \vep_1)$.  By the definition of the
sample sizes, there exists $\vep^* \in (0, \vep_1)$ such that \be \la{hou3}
 \f{1 - \eta}{\vep^2} C_\ell \; 2 \ln \f{2}{\de_\ell} < n_\ell <
\f{1 + \eta}{\vep^2} C_\ell \; 2 \ln \f{2}{\de_\ell} \ee for all $\vep \in (0, \vep^*)$ and $\ell \geq \tau_\vep$.  It follows from
(\ref{hou1}), (\ref{hou2}) and (\ref{hou3}) that \bee H_\ell(z, \vep) & = & \f{n_\ell}{ \ln \f{2}{\de_\ell} } \min \{ - \mscr{M} (  z, \mscr{L}
(z, \vep) ),
\; - \mscr{M} ( z, \mscr{U} (z, \vep) ) \}\\
&  < & \f{ \f{1 + \eta}{\vep^2} C_\ell \; 2 \ln \f{2}{\de_\ell} }{ \ln \f{2}{\de_\ell} } \times  (1 + \eta)^3  \f{ \ka^2 (z) \vep^2 } { 2
\mcal{V} (z) } = (1 + \eta)^{4} \Lm_\ell (z)  < (1 + \vsi) \Lm_\ell (z) \eee and \bee  &  & H_\ell(z, \vep) > \f{ \f{1 - \eta}{\vep^2} C_\ell \;
2 \ln \f{2}{\de_\ell} } { \ln \f{2}{\de_\ell} }  \times  (1 - \eta)^3  \f{ \ka^2 (z) \vep^2 }{ 2 \mcal{V} (z) }  = (1 - \eta)^{4} \Lm_\ell (z) >
(1 - \vsi) \Lm_\ell (z) \eee for all $z \in (\se - \ga, \se + \ga), \; \vep \in (0, \vep^*)$ and $\ell \geq \tau_\vep$.  The claim is thus
established.

For $\mu \in \Se$ with $\Lm_j (\mu) > 1$, let {\small $\vsi \in \li ( 0, 1 - \f{1}{ \sq{ \Lm_j (\mu) } } \ri )$}. By the continuity of $\Lm_\ell
(z)$ on $z \in \Se$, there exists a number $\ga_0 > 0$ such that
\[
\Lm_j (z) \geq (1 - \vsi) \Lm_j (\mu) \qu \tx{for all} \; z \in ( \mu - \ga_0, \mu + \ga_0 ).
\]
Since $\Lm_\ell(\mu)$ is increasing with respect to $\ell$, it must be true that
\[
\Lm_\ell (z) \geq (1 - \vsi) \Lm_j (\mu) \qu \tx{for all} \; \; \ell \geq j  \; \tx{and} \; z \in ( \mu - \ga_0, \mu + \ga_0 ).
\]
By our established claim,  there exist $\ga \in (0, \ga_0)$ and $\vep^* > 0$ such that
\[
H_\ell(z, \vep) > (1 - \vsi) \Lm_\ell (z) \geq (1 - \vsi)^2 \Lm_j (\mu) > 1, \qqu \mscr{L} (z, \vep) \in \Se, \qqu \mscr{U} (z, \vep) \in \Se
\]
for all $z \in (\mu - \ga, \mu + \ga), \; \vep \in (0, \vep^*)$ and $\ell \geq j$.  Making use of this result and (\ref{homed}), we have that
there exist $\ga > 0$ and $\vep^* > 0$ such that $\{ | Z_\ell - \mu  | < \ga \} \subseteq  \{ \bs{D}_\ell = 1  \}$ for all $\vep \in (0,
\vep^*)$ and $\ell \geq j$.  This implies that, there exist $\ga > 0$ and $\vep^* > 0$ such that $\{ \bs{D}_\ell = 0  \} \subseteq \{ | Z_\ell -
\mu  | \geq \ga \}$ for all $\vep \in (0, \vep^*)$ and $\ell \geq j$. Hence, we have shown the first statement of the lemma.

To show the second statement of the lemma, choose {\small $\vsi \in \li (0, \f{1}{ \sq{ \Lm_j (\mu) } } - 1 \ri )$} for $\mu \in \Se$ with
$\Lm_j (\mu) < 1$.  By the continuity of $\Lm_\ell (z)$ on $z \in \Se$, there exists a number $\ga_0 > 0$ such that
\[
\Lm_j (z) \leq (1 + \vsi) \Lm_j (\mu) \qu \tx{for all} \; z \in ( \mu - \ga_0, \mu + \ga_0 ).
\]
Since $\Lm_\ell(\mu)$ is increasing with respect to $\ell$, it must be true that
\[
\Lm_\ell (z) \leq (1 + \vsi) \Lm_j (\mu) \qu \tx{for all} \; \; \ell \leq j  \; \tx{and} \; z \in ( \mu - \ga_0, \mu + \ga_0 ).
\]
By our established claim,  there exist $\ga \in (0, \ga_0)$ and $\vep^* > 0$ such that
\[
H_\ell(z, \vep) < (1 + \vsi) \Lm_\ell (z) \leq (1 + \vsi)^2 \Lm_j (\mu) < 1, \qqu \mscr{L} (z, \vep) \in \Se, \qqu \mscr{U} (z, \vep) \in \Se
\]
for all $z \in (\mu - \ga, \mu + \ga), \; \vep \in (0, \vep^*)$ and $\ell \leq j$. Making use of this result and (\ref{homec}), we have that
there exist $\ga > 0$ and $\vep^* > 0$ such that $\{ | Z_\ell - \mu  | < \ga \} \subseteq  \{ \bs{D}_\ell = 0  \}$ for all $\vep \in (0,
\vep^*)$ and $\ell \leq j$.  This implies that, there exist $\ga > 0$ and $\vep^* > 0$ such that $\{ \bs{D}_\ell = 1  \} \subseteq \{ | Z_\ell -
\mu  | \geq \ga \}$ for all $\vep \in (0, \vep^*)$ and $\ell \leq j$. So, we have shown the second statement of the lemma. This completes the
proof of the lemma.

\epf

\sect{Proof of Theorem \ref{NLST}}  \la{NLSTapp}

\subsection{Proof of Property (I)}

To establish (\ref{perty891}), we can make use of a similar method as that of the counterpart of Theorem \ref{CHCDF} and Lemma
\ref{finitenl899group}  as follows.

\beL \la{finitenl899group} There exist a number $\ga > 0$ and an integer $m
> \tau_\vep$ such that $\{ \bs{l} > \ell \} \subseteq \{ | Z_\ell - \mu | \geq \ga \}$ for all $\ell > m$. \eeL

\bpf  To show the lemma, we need to prove two claims as follows.

(i): For $\mu \in \Se$, there exist a number $\ga > 0$ and  a positive integer $m > \tau_\vep$  such that
\[
\li \{  \mcal{V} \li ( Z_\ell + \ro [\mscr{U} (Z_\ell, \vep) - Z_\ell ] \ri )  > \f{ n_\ell [ Z_\ell - \mscr{U} (Z_\ell, \vep) ]^2 } { \ln
\f{1}{\de_\ell}  }  \ri \} \subseteq \{  | Z_\ell - \mu | \geq \ga   \} \qu \tx{for all $\ell > m$}.
\]

(ii): For $\mu \in \Se$, there exist a number $\ga > 0$ and  a positive integer $m > \tau_\vep$  such that
\[
\li \{  \mcal{V} \li ( Z_\ell + \ro [\mscr{L} (Z_\ell, \vep) - Z_\ell ] \ri )  > \f{ n_\ell [ Z_\ell - \mscr{L} (Z_\ell, \vep) ]^2 } { \ln
\f{1}{\de_\ell}  }  \ri \} \subseteq \{  | Z_\ell - \mu | \geq \ga   \} \qu \tx{for all $\ell > m$}.
\]

To show the first claim, note that there exists $\ga > 0$ such that  $\mu - \ga, \; \mu + \ga \in \Se$ and that $\mscr{U} (\se, \vep) - \se \geq
\f{ \mscr{U} (\mu, \vep) - \mu }{2} > 0$ for all $\se \in (\mu - \ga, \mu + \ga) \subseteq \Se$.  Define $S = \{ \se \in (\mu - \ga, \mu + \ga):
\se + \ro [ \mscr{U} (\se, \vep) - \se ] \in \Se \}$. If $S = \emptyset$, then \[ \li \{ \mcal{V} \li ( Z_\ell + \ro [\mscr{U} (Z_\ell, \vep) -
Z_\ell ] \ri ) > \f{n_\ell}{\ln \f{1}{\de_\ell}} [ Z_\ell - \mscr{U} (Z_\ell, \vep) ]^2, \; | Z_\ell - \mu | < \ga \ri \} = \emptyset \]
 for $\ell \geq \tau_\vep$.  In the other case that $S \neq \emptyset$, there exists a constant $D < \iy$ such that
\[
\sup_{ z \in S } \mcal{V}( Z_\ell + \ro [\mscr{U} (Z_\ell, \vep) - Z_\ell ] ) < D.
\]
Making use of this fact and the observation that $\lim_{\ell \to \iy} \f{\ln \f{1}{\de_\ell}} {n_\ell} = 0$, we have that there exists an
integer $m > \tau_\vep$ such that {\small \bee & & \li \{ \mcal{V} \li ( Z_\ell + \ro [\mscr{U} (Z_\ell,
\vep) - Z_\ell ] \ri ) > \f{n_\ell}{\ln \f{1}{\de_\ell}} [ Z_\ell - \mscr{U} (Z_\ell, \vep) ]^2, \; | Z_\ell - \mu | < \ga  \ri \}\\
&  & = \li \{ n_\ell [ Z_\ell - \mscr{U} (Z_\ell, \vep) ]^2 < \mcal{V} \li ( Z_\ell + \ro [\mscr{U} (Z_\ell, \vep) - Z_\ell ] \ri ) \ln
\f{1}{\de_\ell}, \; Z_\ell + \ro [\mscr{U} (Z_\ell, \vep) - Z_\ell ]  \in \Se, \; |
Z_\ell - \mu | < \ga  \ri \}\\
&  & \subseteq \li \{ n_\ell [ Z_\ell - \mscr{U} (Z_\ell, \vep) ]^2 < D \ln \f{1}{\de_\ell}, \;  |
Z_\ell - \mu | < \ga  \ri \}\\
&  & \subseteq \li \{  \mscr{U} (Z_\ell, \vep) - Z_\ell < \f{ \mscr{U} (\mu, \vep) - \mu  }{2}, \; | Z_\ell - \mu | < \ga \ri \} = \emptyset
\eee} for all $\ell > m$.   It follows that {\small \bee &  & \li \{ \mcal{V} \li ( Z_\ell + \ro [\mscr{U} (Z_\ell, \vep) -
Z_\ell ] \ri ) > \f{n_\ell}{\ln \f{1}{\de_\ell}} [ Z_\ell - \mscr{U} (Z_\ell, \vep) ]^2  \ri \}\\
 & \subseteq & \li \{ n_\ell [
Z_\ell - \mscr{U} (Z_\ell, \vep) ]^2 < \mcal{V} \li ( Z_\ell + \ro [\mscr{U} (Z_\ell, \vep) - Z_\ell ] \ri ) \ln \f{1}{\de_\ell}, \; | Z_\ell -
\mu | < \ga  \ri \}  \cup \{ |
Z_\ell - \mu | \geq \ga \}\\
& = & \{ | Z_\ell - \mu | \geq \ga \}  \eee} for all $\ell > m$. This establishes the first claim. In a similar manner, we can show the second
claim.  Finally, making use of the two established claims and the definition of the stopping rule completes the proof of the lemma.  \epf

\subsection{Proof of Property (II)}

To prove (\ref{propertyii1}) and (\ref{propertyii2}), we can apply Lemma \ref{pro689a} and Lemma \ref{pro689bnl} as follows.

\beL

\la{pro689bnl}

\bel &  & \Pr \li \{ \liminf_{\vep \downarrow 0} \Lm_{\bs{l}} (\mu) \geq 1 \ri \} = 1, \la{cite369anl}\\
&  & \Pr \li \{ \tx{Both} \; \limsup_{\vep \downarrow 0} \Lm_{\bs{l} - 1} (\mu) \leq 1 \; \tx{and} \;   \liminf_{\vep \downarrow 0} \bs{l} > 1
\; \tx{hold or} \; \lim_{\vep \downarrow 0} \bs{l}  = 1 \ri \} = 1.  \la{cite369bnl} \eel

\eeL

\bpf To simplify notations, define
\[
\bs{L}_{\bs{l}} = \mscr{L} ( Z_{\bs{l}}, \vep ), \qqu \bs{L}_{\bs{l} - 1} = \mscr{L} ( Z_{\bs{l} - 1}, \vep ), \qqu \bs{U}_{\bs{l}} = \mscr{U} (
Z_{\bs{l}}, \vep ), \qqu \bs{U}_{\bs{l} - 1} = \mscr{U} ( Z_{\bs{l} - 1}, \vep ).
\]
Define {\small \bee &  & \mscr{E}_1 = \li \{  \mcal{V} \li ( Z_{\bs{l}} + \ro [ \bs{U}_{\bs{l}} - Z_{\bs{l}} ] \ri )  \leq \f{ n_{\bs{l}} [
Z_{\bs{l}} - \bs{U}_{\bs{l}} ]^2  }{ \ln \f{1}{ \de_{{\bs{l}}} }  }, \;  \mcal{V} \li ( Z_{\bs{l}} + \ro [ \bs{L}_{\bs{l}} - Z_{\bs{l}} ] \ri )
\leq \f{ n_{\bs{l}} [ Z_{\bs{l}} - \bs{L}_{\bs{l}} ]^2  }{ \ln \f{1}{ \de_{{\bs{l}}} }  } \ri \}, \\
&  & \mscr{E}_2 =  \li \{  \mcal{V} \li ( Z_{{\bs{l}} - 1} + \ro [ \bs{U}_{{\bs{l}}-1} - Z_{{\bs{l}}-1} ] \ri ) > \f{ n_{\bs{l} - 1} [
Z_{{\bs{l}}-1} - \bs{U}_{{\bs{l}}-1} ]^2  }{ \ln \f{1}{ \de_{{\bs{l}}-1} }  }  \ri \}\\
&  & \qqu \; \bigcup  \li \{  \mcal{V} \li ( Z_{{\bs{l}} - 1} + \ro [ \bs{L}_{{\bs{l}}-1} - Z_{{\bs{l}}-1} ] \ri ) > \f{ n_{\bs{l} - 1} [
Z_{{\bs{l}}-1} - \bs{L}_{{\bs{l}}-1} ]^2  }{ \ln \f{1}{ \de_{{\bs{l}}-1} }  }  \ri \},\\
&  & \mscr{E}_3 =   [  \mscr{E}_2 \cap \{ {\bs{l}} > \tau_\vep  \} ] \cup \{ {\bs{l}} = \tau_\vep  \},\\
&  & \mscr{E}_4 = \li \{ \lim_{\vep \downarrow 0} {\mbf{n}} = \iy, \;  \lim_{\vep \downarrow 0} Z_{{\bs{l}}}
= \mu, \; \lim_{\vep \downarrow 0} Z_{{\bs{l}} - 1} = \mu \ri \},\\
&  &  \mscr{E}_5 = \mscr{E}_1 \cap \mscr{E}_3 \cap \mscr{E}_4,\\
&  &  \mscr{E}_6 =  \li \{ \liminf_{\vep \downarrow 0} \Lm_{\bs{l}} (\mu) \geq 1 \ri \},\\
&  &  \mscr{E}_7 =  \li \{ \limsup_{\vep \downarrow 0} \Lm_{\bs{l} - 1} (\mu) \leq 1, \;  \liminf_{\vep \downarrow 0} \bs{l} > 1 \; \tx{or} \;
\lim_{\vep \downarrow 0} \bs{l}  = 1 \ri \}.  \eee} By the assumption that $n_{\tau_\vep} \to \iy$ as $\vep \downarrow 0$, we have that $\{
\lim_{\vep \downarrow 0} \mbf{n} = \iy \}$ is a sure event. It follows from the strong law of large numbers that $\mscr{E}_4$ is an almost sure
event. By the definition of the stopping rule, we have that  $\mscr{E}_1  \cap \mscr{E}_3$ is an almost sure event. Hence, $ \mscr{E}_5$ is an
almost sure event.  To show (\ref{cite369anl}), it suffices to show that $\mscr{E}_5 \subseteq \mscr{E}_6$. For this purpose, we let $\om \in
\mscr{E}_5$ and attempt to show that $\om \in \mscr{E}_6$.  For simplicity of notations, let $\ell = \bs{l} (\om), \;  z_\ell = Z_{\bs{l}}
(\om), \; z_{\ell-1} = Z_{\bs{l} - 1} (\om)$,
\[
L_\ell = \mscr{L} ( z_\ell, \vep), \qu U_\ell = \mscr{U} ( z_\ell, \vep), \qu L_{\ell-1} = \mscr{L} ( z_{\ell-1}, \vep), \qu U_{\ell-1} =
\mscr{U} ( z_{\ell-1}, \vep).
\]
We restrict $\vep > 0$ to be sufficiently small such that
\[ z_\ell + \ro [U_\ell - z_\ell ] \in \Se, \qqu z_\ell + \ro [L_\ell - z_\ell ] \in \Se.  \]  As a consequence of $\om \in \mscr{E}_1$,
\[
\f{ n_\ell [ \ka (\mu) \vep ]^2 }{ \mcal{V} (\mu) }  \geq  \f{ \mcal{V} \li ( z_\ell + \ro [U_\ell - z_\ell ] \ri ) } { \mcal{V} (\mu) } \li [
\f{ \ka (\mu) \vep  }{ U_\ell - z_\ell } \ri ]^2 \ln \f{1}{\de_\ell} \qu \tx{and}
\]
\[
\f{ n_\ell [ \ka (\mu) \vep ]^2 }{ \mcal{V} (\mu) }  \geq  \f{ \mcal{V} \li ( z_\ell + \ro [L_\ell - z_\ell ] \ri ) } { \mcal{V} (\mu) }  \li [
\f{ \ka (\mu) \vep  }{ z_\ell - L_\ell } \ri ]^2 \ln \f{1}{\de_\ell}.
\]
These inequalities can be written as
\[
\Lm_\ell (\mu) = \f{ \ka^2 (\mu) C_\ell }{ \mcal{V} (\mu) }  \geq \f{ C_\ell } {\vep^2 n_\ell} \f{ \mcal{V} \li ( z_\ell + \ro [U_\ell - z_\ell
] \ri ) } { \mcal{V} (\mu) } \li [ \f{ \ka (\mu) \vep  }{ U_\ell - z_\ell } \ri ]^2 \ln \f{1}{\de_\ell} \qu \tx{and}
\]
\[
\Lm_\ell (\mu) = \f{ \ka^2 (\mu) C_\ell }{ \mcal{V} (\mu) }  \geq \f{ C_\ell } {\vep^2 n_\ell} \f{ \mcal{V} \li ( z_\ell + \ro [L_\ell - z_\ell
] \ri ) } { \mcal{V} (\mu) }  \li [ \f{  \ka (\mu) \vep  }{ z_\ell - L_\ell } \ri ]^2 \ln \f{1}{\de_\ell}.
\]
It follows that $\liminf_{\vep \downarrow 0} \Lm_\ell (\mu) \geq 1$. Hence, $\om \in \mscr{E}_6$ and thus (\ref{cite369anl}) holds.

To show (\ref{cite369bnl}), it suffices to show that $\om \in \mscr{E}_7$ for $\om \in \mscr{E}_5$.  Since $\om \in \mscr{E}_3$, we have that
either $\ell = \tau_\vep$ or
\[
\f{ n_{\ell - 1} [ \ka (\mu) \vep ]^2 }{ \mcal{V} (\mu) }  < \f{ \mcal{V} \li ( z_{\ell - 1} + \ro [U_{\ell-1} - z_{\ell - 1} ] \ri ) } {
\mcal{V} (\mu) }  \li [ \f{  \ka (\mu) \vep  }{ U_{\ell-1} - z_{\ell - 1} } \ri ]^2  \ln \f{1}{\de_{\ell-1}} \qu \tx{or}
\]
\[
\f{ n_{\ell - 1} [ \ka (\mu) \vep ]^2 }{ \mcal{V} (\mu) }  < \f{ \mcal{V} \li ( z_{\ell - 1} + \ro [L_{\ell-1} - z_{\ell - 1} ] \ri ) } {
\mcal{V} (\mu) }  \li [ \f{  \ka (\mu) \vep  }{ z_{\ell - 1} - L_{\ell-1} } \ri ]^2 \ln \f{1}{\de_{\ell-1}}
\]
must be true for small enough $\vep > 0$.  It follows that either $\ell = \tau_\vep$ or
\[
\Lm_{\ell - 1} (\mu) = \f{ \ka^2 (\mu) C_{\ell - 1} }{ \mcal{V} (\mu) }  < \f{ C_{\ell - 1} } {\vep^2 n_{\ell - 1}}  \f{ \mcal{V} \li ( z_{\ell
- 1} + \ro [U_{\ell-1} - z_{\ell - 1} ] \ri ) } { \mcal{V} (\mu) } \li [ \f{ \ka (\mu) \vep }{ U_{\ell-1} - z_{\ell - 1} } \ri ]^2 \ln
\f{1}{\de_{\ell-1}} \qu \tx{or}
\]
\[
\Lm_{\ell - 1} (\mu) = \f{ \ka^2 (\mu) C_{\ell - 1} }{ \mcal{V} (\mu) }  < \f{ C_{\ell - 1} } {\vep^2 n_{\ell - 1}}  \f{ \mcal{V} \li ( z_{\ell
- 1} + \ro [L_{\ell-1} - z_{\ell - 1} ] \ri ) } { \mcal{V} (\mu) } \li [ \f{ \ka (\mu) \vep }{ z_{\ell - 1} - L_{\ell-1} } \ri ]^2 \ln
\f{1}{\de_{\ell-1}}.
\]
must be true for small enough $\vep > 0$.  This implies that either $\lim_{\vep \downarrow 0} \ell = 1$ or
\[
\limsup_{\vep \downarrow 0} \Lm_{\ell - 1} (\mu) \leq 1, \qqu \liminf_{\vep \downarrow 0} \ell > 1 \]
 must be true.  Hence, $\om \in \mscr{E}_7$ and thus (\ref{cite369bnl}) holds.

\epf

\subsection{Proof of Properties (III) -- (VII)}

Property (III) follows from Lemma \ref{lem22888} and the established Property (II). A key result to establish Properties (IV) -- (VII) is Lemma
\ref{tubeNAP}, which will be developed in the sequel.  Once Lemma \ref{tubeNAP} is established,

Property (IV) follows from Lemmas \ref{lem20888}, \ref{lem21888} and \ref{tubeNAP};

Property (V) follows from Lemmas \ref{lem20888}, \ref{lem23888} and \ref{tubeNAP};

Property (VI) follows from Lemmas \ref{lem20888}, \ref{lem24888} and \ref{tubeNAP};

Property (VII) follows from Lemmas \ref{lem20888}, \ref{lem25888} and \ref{tubeNAP}.

\beL

\la{tubeNAP}

For $\mu \in \Se$, if $\Lm_j (\mu)  > 1$ for some index $j \geq 1$, then there exist $\ga > 0$ and $\vep^* > 0$ such that $\{ \bs{D}_\ell = 0 \}
\subseteq \{ | Z_\ell - \mu  | \geq \ga  \}  \qu \tx{for all} \; \; \ell \geq j$ and $\vep \in (0, \vep^*)$.

Similarly, for $\mu \in \Se$, if $\Lm_j (\mu)  < 1$ for some index $j \geq 1$, then there exist $\ga > 0$ and $\vep^* > 0$ such that $\{
\bs{D}_\ell = 1 \} \subseteq \{ | Z_\ell - \mu  | \geq \ga  \}  \qu \tx{for all} \; \; \ell \leq j$ and $\vep \in (0, \vep^*)$. \eeL

\bpf

Define \be \la{reuse} H_\ell(z, \vep) = \f{n_\ell} { \ln \f{1}{\de_\ell}  } \min \li \{  \f{ [ \mscr{U} (z, \vep) - z ]^2 } { \mcal{V} \li ( \ro
\mscr{U} (z, \vep) + (1 - \ro) z \ri ) }, \; \f{ [ \mscr{L} (z, \vep)  - z ]^2 } { \mcal{V} \li ( \ro \mscr{L} (z, \vep) + (1 - \ro) z \ri ) }
\ri \} \ee for $\ell \geq \tau_\vep$.  By the definitions of $H_\ell(z, \vep)$, the decision variable $\bs{D}_\ell$, and the stopping rule, we
have that \be \la{homef}
 \{ H_\ell(
Z_\ell, \vep) \geq 1, \; \ro \mscr{L} (Z_\ell, \vep) + (1 - \ro) Z_\ell \in \Se, \; \ro \mscr{U} (Z_\ell, \vep) + (1 - \ro) Z_\ell \in \Se \}
\subseteq \{ \bs{D}_\ell = 1 \} \ee and \be \la{homee}
 \{ H_\ell(
Z_\ell, \vep) < 1, \; \ro \mscr{L} (Z_\ell, \vep) + (1 - \ro) Z_\ell \in \Se, \; \ro \mscr{U} (Z_\ell, \vep) + (1 - \ro) Z_\ell \in \Se\}
\subseteq \{ \bs{D}_\ell = 0 \} \ee for $\ell \geq \tau_\vep$.

Next, we claim that for any $\se \in \Se$ and $\vsi \in (0, 1)$, there exist $\ga > 0$ and $\vep^* > 0$ such that $\ro \mscr{L} (z, \vep) + (1 -
\ro) z \in \Se, \; \ro \mscr{U} (z, \vep) + (1 - \ro) z \in \Se$ and that $(1 - \vsi) \Lm_\ell (z) < H_\ell(z, \vep) < (1 + \vsi) \Lm_\ell (z)$
for all $z \in (\se - \ga, \se + \ga), \; \vep \in (0, \vep^*)$ and $\ell \geq \tau_\vep$.

To prove this claim, note that for $\vsi \in (0, 1)$, there exists $\eta \in (0, 1)$ such that $\f{1 + \eta }{(1 - \eta)^2 (1 + \vsi)}   < 1 <
\f{1 - \eta}{(1 + \eta)^2 (1 - \vsi)}$. Note that  there exist $\ga > 0$ and $\vep_1
> 0$ such that
\[
\ro \mscr{L} (z, \vep) + (1 - \ro) z \in \Se, \qqu \ro \mscr{U} (z, \vep) + (1 - \ro) z \in \Se \]
 and that
\bel &   &  1 - \eta < \li [ \f{  \ka (z) \vep } {z - \mscr{L} (z, \vep)  } \ri ]^2 < 1 + \eta, \la{inq33883a}\\
&  &  1 - \eta < \li [ \f{  \ka (z) \vep } { \mscr{U} (z, \vep) - z } \ri ]^2 < 1 + \eta, \la{inq33884b}\\
&  &  1 - \eta < \f{\mcal{V} \li ( \ro \mscr{L} (z, \vep) + (1 - \ro) z
 \ri )} { \mcal{V} (z) } < 1 + \eta, \la{inq33885c}\\
 &  & 1 - \eta < \f{\mcal{V} \li ( \ro \mscr{U} (z, \vep) + (1 - \ro) z
 \ri )} { \mcal{V} (z) }  < 1 + \eta \la{inq33886d}
\eel for all $z \in (\se - \ga, \se + \ga)$ and $\vep \in (0, \vep_1)$.  Here, inequalities (\ref{inq33883a}) and (\ref{inq33884b})  are due to
the continuity of $\ka(z)$ and the assumption associated with (\ref{uniformasp}).  Inequalities (\ref{inq33885c}) and (\ref{inq33886d}) follow
from the fact that $\f{\mcal{V} \li ( \ro \mscr{L} (z, \vep) + (1 - \ro) z
 \ri )} { \mcal{V} (z) }$ and $\f{\mcal{V} \li ( \ro \mscr{U} (z, \vep) + (1 - \ro) z
 \ri )} { \mcal{V} (z) }$ are
continuous functions of $z$ and $\vep$ which tend to $1$ as $(z, \vep)$ tends to $(\se, 0)$.  Therefore, \bel &  & \f{\mcal{V} (z)}{ [\ka(z)
\vep]^2 } (1 - \eta)^2 < \f{ \mcal{V} \li ( \ro \mscr{L} (z, \vep) + (1 - \ro) z
 \ri ) } { [ \mscr{L} (z, \vep)  - z
]^2 }   < (1 + \eta)^2 \f{\mcal{V} (z)}{ [\ka(z) \vep]^2 }, \la{here8a}\\
 &  &  \f{\mcal{V} (z)}{ [\ka(z) \vep]^2 } (1 - \eta)^2 < \f{ \mcal{V}
\li ( \ro \mscr{U} (z, \vep) + (1 - \ro) z \ri ) } { [ \mscr{U} (z, \vep) - z ]^2 }  < (1 + \eta)^2 \f{\mcal{V} (z)}{ [\ka(z) \vep]^2 }
\la{here8b} \eel for all $z \in (\se - \ga, \se + \ga)$ and $\vep \in (0, \vep_1)$.  Here,  (\ref{here8a}) follows from (\ref{inq33883a}) and
(\ref{inq33885c}).  On the other hand, (\ref{here8b}) is due to (\ref{inq33884b}) and (\ref{inq33886d}). By the definition of the sample sizes,
there exists $\vep^* \in (0, \vep_1)$ such that \be \la{here69}
 \f{1 - \eta}{\vep^2} C_\ell \ln \f{1}{\de_\ell} < n_\ell < \f{1 + \eta}{\vep^2}
C_\ell \ln \f{1}{\de_\ell}. \ee It follows from (\ref{here8a}), (\ref{here8b}) and (\ref{here69}) that \bee  &  & H_\ell(z, \vep)  < \f{ \f{1 +
\eta}{\vep^2} C_\ell \ln \f{1}{\de_\ell} } { \ln \f{1}{\de_\ell} }  \f{ [\ka(z) \vep]^2 } {\mcal{V} (z) (1 - \eta)^2 }   = \f{1 + \eta} {(1 -
\eta)^2} \Lm_\ell (z) < (1 + \vsi) \Lm_\ell (z) \eee and \bee  & & H_\ell(z, \vep)
> \f{ \f{1 - \eta}{\vep^2} C_\ell \ln \f{1}{\de_\ell} }{\ln \f{1}{\de_\ell} }  \f{ [\ka(z) \vep]^2 } {\mcal{V} (z) (1 + \eta)^2}
 =  \f{1 - \eta} {(1 + \eta)^2} \Lm_\ell (z) > (1 - \vsi) \Lm_\ell (z) \eee for all $z \in (\se - \ga, \se + \ga)$ and $\vep \in (0,
 \vep^*)$.  Our claim is thus established.

 For $\mu \in \Se$ with $\Lm_j (\mu) > 1$, let {\small $\vsi \in \li ( 0, 1 - \f{1}{ \sq{ \Lm_j (\mu) } } \ri )$}. By the continuity of $\Lm_\ell
(z)$ on $z \in \Se$, there exists a number $\ga_0 > 0$ such that
\[
\Lm_j (z) \geq (1 - \vsi) \Lm_j (\mu) \qu \tx{for all} \; z \in ( \mu - \ga_0, \mu + \ga_0 ).
\]
Since $\Lm_\ell(\mu)$ is increasing with respect to $\ell$, it must be true that
\[
\Lm_\ell (z) \geq (1 - \vsi) \Lm_j (\mu) \qu \tx{for all} \; \; \ell \geq j  \; \tx{and} \; z \in ( \mu - \ga_0, \mu + \ga_0 ).
\]
By our established claim,  there exist $\ga \in (0, \ga_0)$ and $\vep^* > 0$ such that
\[
H_\ell(z, \vep) > (1 - \vsi) \Lm_\ell (z) \geq (1 - \vsi)^2 \Lm_j (\mu) > 1, \qqu \ro \mscr{L} (z, \vep) + (1 - \ro) z \in \Se, \qu \ro \mscr{U}
(z, \vep) + (1 - \ro) z \in \Se
\]
for all $z \in (\mu - \ga, \mu + \ga), \; \vep \in (0, \vep^*)$ and $\ell \geq j$.  Making use of this result and (\ref{homef}), we have that
there exist $\ga > 0$ and $\vep^* > 0$ such that $\{ | Z_\ell - \mu  | < \ga \} \subseteq  \{ \bs{D}_\ell = 1  \}$ for all $\vep \in (0,
\vep^*)$ and $\ell \geq j$.  This implies that, there exist $\ga > 0$ and $\vep^* > 0$ such that $\{ \bs{D}_\ell = 0  \} \subseteq \{ | Z_\ell -
\mu  | \geq \ga \}$ for all $\vep \in (0, \vep^*)$ and $\ell \geq j$. Hence, we have shown the first statement of the lemma.

To show the second statement of the lemma, choose {\small $\vsi \in \li (0, \f{1}{ \sq{ \Lm_j (\mu) } } - 1 \ri )$} for $\mu \in \Se$ with
$\Lm_j (\mu) < 1$.  By the continuity of $\Lm_\ell (z)$ on $z \in \Se$, there exists a number $\ga_0 > 0$ such that
\[
\Lm_j (z) \leq (1 + \vsi) \Lm_j (\mu) \qu \tx{for all} \; z \in ( \mu - \ga_0, \mu + \ga_0 ).
\]
Since $\Lm_\ell(\mu)$ is increasing with respect to $\ell$, it must be true that
\[
\Lm_\ell (z) \leq (1 + \vsi) \Lm_j (\mu) \qu \tx{for all} \; \; \ell \leq j  \; \tx{and} \; z \in ( \mu - \ga_0, \mu + \ga_0 ).
\]
By our established claim,  there exist $\ga \in (0, \ga_0)$ and $\vep^* > 0$ such that
\[
H_\ell(z, \vep) < (1 + \vsi) \Lm_\ell (z) \leq (1 + \vsi)^2 \Lm_j (\mu) < 1, \qqu \ro \mscr{L} (z, \vep) + (1 - \ro) z \in \Se, \qu \ro \mscr{U}
(z, \vep) + (1 - \ro) z \in \Se
\]
for all $z \in (\mu - \ga, \mu + \ga), \; \vep \in (0, \vep^*)$ and $\ell \leq j$. Making use of this result and (\ref{homee}), we have that
there exist $\ga > 0$ and $\vep^* > 0$ such that $\{ | Z_\ell - \mu  | < \ga \} \subseteq  \{ \bs{D}_\ell = 0  \}$ for all $\vep \in (0,
\vep^*)$ and $\ell \leq j$.  This implies that, there exist $\ga > 0$ and $\vep^* > 0$ such that $\{ \bs{D}_\ell = 1  \} \subseteq \{ | Z_\ell -
\mu  | \geq \ga \}$ for all $\vep \in (0, \vep^*)$ and $\ell \leq j$. So, we have shown the second statement of the lemma. This completes the
proof of the lemma.

\epf

 \sect{Proof of Theorem \ref{Seqgen} }
\la{Seqgen_app}

\subsection{Proof of Property (I) }

To establish (\ref{perty891}), we can make use of the following Lemma \ref{finitegroup}.

\beL \la{finitegroup} If the random variable $X$ has mean $\mu$ and variance $\nu$, then there exist a number $\ga > 0$ and an integer $m >
\tau_\vep$ such that $\{ \bs{l} > \ell \} \subseteq \{ | Z_\ell - \mu | \geq \ga \; \tx{or} \; | V_{n_\ell} - \nu | \geq \ga \}$ for all $\ell
> m$.

\eeL

\bpf

Note that there exists $\ga > 0$ such that
 \[
\mscr{U}( \se, \vep  ) - \se > \f{  \mscr{U}( \mu, \vep  ) - \mu }{2} > 0
 \]
 for all $\se \in (\mu - \ga, \mu + \ga)$.  Since $C_\ell \to \iy$ as $\ell \to \iy$, it follows that
 $\f{ \ln \f{1}{\de_\ell} }{n_\ell} \to 0$ as $\ell \to \iy$.  Therefore, there exists an integer $m > \tau_\vep$ such that
\bee &  &  \li \{ V_{n_\ell} + \f{\ro}{n_\ell} > \f{ n_\ell [ \mscr{U}( Z_\ell, \vep  ) - Z_\ell ]^2}{  \ln \f{1}{\de_\ell} }, \;  | Z_\ell -
\mu |
< \ga, \; | V_{n_\ell} - \nu | < \ga \ri \}\\
&  & \subseteq  \li \{ V_{n_\ell} + \f{\ro}{n_\ell} > \f{ n_\ell [ \f{  \mscr{U}( \mu, \vep  ) - \mu }{2} ]^2}{  \ln \f{1}{\de_\ell} }, \;  |
Z_\ell - \mu | < \ga, \; | V_{n_\ell} - \nu | < \ga \ri \} = \emptyset \eee for all $\ell > m$.  It follows that \be \la{cite9886g}  \li \{
V_{n_\ell} + \f{\ro}{n_\ell} > \f{ n_\ell [ \mscr{U}( Z_\ell, \vep  ) - Z_\ell ]^2}{  \ln \f{1}{\de_\ell} } \ri \} \subseteq  \li \{ | Z_\ell -
\mu | \geq \ga \; \tx{or} \; | V_{n_\ell} - \nu | \geq \ga \ri \} \qu \tx{for all $\ell > m$}. \ee Thus, we have shown that there exist a number
$\ga > 0$ and an integer $m > \tau_\vep$ such that (\ref{cite9886g}) holds.  In a similar manner, we can show that there exist a number $\ga >
0$ and an integer $m > \tau_\vep$ such that
\[ \li \{ V_{n_\ell} + \f{\ro}{n_\ell} > \f{ n_\ell [ \mscr{L}( Z_\ell, \vep  ) - Z_\ell ]^2}{  \ln \f{1}{\de_\ell} } \ri \} \subseteq  \li \{ | Z_\ell - \mu
| \geq \ga \; \tx{or} \; | V_{n_\ell} - \nu | \geq \ga \ri \} \qu \tx{for all $\ell > m$}. \] Finally, making use of these results and the
definition of the stopping rule completes the proof of the lemma.

\epf

We are now in a position to show (\ref{perty891}). According to Lemma \ref{finitegroup}, there exist an integer $m > \tau_\vep$ and a positive
number $\ga
> 0$ such that \be \la{makeusegroup}
 \Pr \{ \bs{l} > \ell \}
\leq \Pr \{ | Z_\ell - \mu | \geq \ga \} + \Pr \{ | V_{n_\ell} - \nu | \geq \ga  \} \qqu \tx{for all $\ell > m$}. \ee Making use of
(\ref{makeusegroup}), the weak law of large numbers, and Lemma \ref{vipnow},  we have $\lim_{\ell \to \iy} \Pr \{ \bs{l} > \ell \} = 0$, which
implies $\Pr \{ \mbf{n} < \iy \} = 1$.  By virtue of (\ref{mostvip}) of Lemma \ref{nonuniformBE}, Lemma \ref{vipnow},  and (\ref{makeusegroup}),
we have \bee  \bb{E} [
\mbf{n} ] & \leq &  n_{m+1} + \sum_{\ell > m} [ \Pr \{ | Z_\ell - \mu | \geq \ga \} + \Pr \{ | V_{n_\ell} - \nu | \geq \ga  \} ]\\
&  \leq &  n_{m+1} + \sum_{\ell > m}  (n_{\ell + 1} - n_\ell) \li [ \exp \li ( - \f{n_\ell}{2} \f{ \ga^2 }{ \nu } \ri ) + \f{2 C }{n_\ell^2} \f{ \mscr{W} }{ \ga^3 } \ri ]\\
&   & + \sum_{\ell > m} (n_{\ell + 1} - n_\ell) \li [ \exp \li ( - \f{n_\ell}{4} \f{ \ga}{ \nu } \ri ) +  \exp \li ( - \f{n_\ell}{8} \f{ \ga^2
}{ \varpi } \ri ) + \f{4 C }{n_\ell^2 \ga^3} \li (  \sq{2} \mscr{W} \ga^{3\sh 2} +  4 \mscr{V} \ri ) \ri ]\\
& < & \iy. \eee

\subsection{Proof of Property (II) }

To prove (\ref{propertyii1}) and (\ref{propertyii2}), we can apply Lemma \ref{pro689a} and Lemma \ref{pro689bfull} as follows.

\beL

\la{pro689bfull}

\bel &  & \Pr \li \{ \liminf_{\vep \downarrow 0} \Lm_{\bs{l}} (\mu, \nu) \geq 1 \ri \} = 1, \la{cite369afull}\\
&  & \Pr \li \{ \tx{Both} \; \limsup_{\vep \downarrow 0} \Lm_{\bs{l} - 1} (\mu, \nu) \leq 1 \; \tx{and} \;  \liminf_{\vep \downarrow 0} \bs{l} >
1 \; \tx{hold or} \; \lim_{\vep \downarrow 0} \bs{l}  = 1 \ri \} = 1.  \la{cite369bfull} \eel

\eeL

\bpf

To simplify notations, define $W_\ell = V_{n_\ell}, \; \ell \geq \tau_\vep$ and
\[
\bs{L}_{\bs{l}} = \mscr{L} ( Z_{\bs{l}}, \vep ), \qqu \bs{L}_{\bs{l} - 1} = \mscr{L} ( Z_{\bs{l} - 1}, \vep ), \qqu \bs{U}_{\bs{l}} = \mscr{U} (
Z_{\bs{l}}, \vep ), \qqu \bs{U}_{\bs{l} - 1} = \mscr{U} ( Z_{\bs{l} - 1}, \vep ).
\]
Define {\small \bee &  & \mscr{E}_1 = \li \{ W_{\bs{l}} + \f{\ro}{n_{\bs{l}}} \leq \f{ n_{\bs{l}} [ Z_{\bs{l}} -  \bs{U}_{\bs{l}}  ]^2 }{ \ln
\f{1}{\de_{{\bs{l}}}} }, \;  W_{\bs{l}} + \f{\ro}{n_{\bs{l}}} \leq \f{ n_{\bs{l}} [  \bs{L}_{\bs{l}}  - Z_{\bs{l}} ]^2 }{ \ln
\f{1}{\de_{{\bs{l}}}} } \ri \}, \\
&  & \mscr{E}_2 = \li \{ W_{{\bs{l}} - 1} + \f{\ro}{n_{\bs{l}-1} } > \f{ n_{\bs{l} - 1}  [ Z_{{\bs{l}} - 1} -  \bs{U}_{{\bs{l}} - 1} ]^2
}{ \ln \f{1}{\de_{{\bs{l}}-1}} }, \; \bs{l} > \tau_\vep \ri \} \\
&  & \qqu \bigcup \li \{ W_{{\bs{l}}-1} + \f{\ro}{n_{\bs{l} - 1}} > \f{ n_{\bs{l} - 1} [  \bs{L}_{{\bs{l}}-1}  - Z_{{\bs{l}}-1} ]^2 }{ \ln
\f{1}{\de_{{\bs{l}}-1}} }, \; \bs{l} > \tau_\vep \ri \} \bigcup \{ \bs{l} = \tau_\vep \},\\
&   & \mscr{E}_3 = \li \{ \lim_{\vep \downarrow 0} {\mbf{n}} = \iy, \; \lim_{\vep \downarrow 0} Z_{{\bs{l}}} = \mu, \; \lim_{\vep \downarrow 0}
W_{{\bs{l}}} = \nu, \; \lim_{\vep \downarrow 0} Z_{{\bs{l}} - 1} = \mu, \; \lim_{\vep \downarrow
0} W_{{\bs{l}} - 1} = \nu \ri \}, \\
&   & \mscr{E}_4 = \mscr{E}_1 \cap \mscr{E}_2 \cap \mscr{E}_3, \\
&  &  \mscr{E}_5 =  \li \{ \liminf_{\vep \downarrow 0} \Lm_{\bs{l}} (\mu, \nu) \geq 1 \ri \},\\
&  &  \mscr{E}_6 =  \li \{ \tx{Both} \; \limsup_{\vep \downarrow 0} \Lm_{\bs{l} - 1} (\mu, \nu) \leq 1 \; \tx{and} \;  \liminf_{\vep \downarrow
0} \bs{l} > 1 \; \tx{hold or} \; \lim_{\vep \downarrow 0} \bs{l}  = 1 \ri \}. \eee} By the assumption that $n_{\tau_\vep} \to \iy$ as $\vep
\downarrow 0$, we have that $\{ \lim_{\vep \downarrow 0} \mbf{n} = \iy \}$ is a sure event. It follows from the strong law of large numbers that
$\mscr{E}_3$ is an almost sure event. By the definition of the stopping rule, we have that  $\mscr{E}_1 \cap \mscr{E}_2$ is an almost sure
event. Hence, $ \mscr{E}_4$ is an almost sure event.  To show (\ref{cite369afull}), it suffices to show that $\mscr{E}_4 \subseteq \mscr{E}_5$.
For this purpose, we let $\om \in \mscr{E}_4$ and attempt to show that $\om \in \mscr{E}_5$.  For simplicity of notations, let $\ell = \bs{l}
(\om), \;  z_\ell = Z_{\bs{l}} (\om), \; z_{\ell-1} = Z_{\bs{l} - 1} (\om), \; w_\ell = W_{\bs{l}} (\om), \; w_{\ell-1} = W_{\bs{l} - 1} (\om)$,
\[
L_\ell = \mscr{L} ( z_\ell, \vep), \qu U_\ell = \mscr{U} ( z_\ell, \vep), \qu L_{\ell-1} = \mscr{L} ( z_{\ell-1}, \vep), \qu U_{\ell-1} =
\mscr{U} ( z_{\ell-1}, \vep).
\]
For small $\vep > 0$, $z_\ell$ and $z_{\ell - 1}$ will be bounded in a neighborhood of $\mu$. For small $\vep > 0$, $w_\ell$ and $w_{\ell - 1}$
will be bounded in a neighborhood of $\nu$. Note that, as a consequence of $\om \in \mscr{E}_1$,
\[
\Lm_\ell (\mu, \nu) = \f{ \ka^2 (\mu) C_\ell }{  \nu  }  \geq \f{C_\ell \ln \f{1}{\de_\ell} }{\vep^2 n_\ell}  \f{ w_\ell + \f{\ro}{n_\ell}  } {
\nu } \li [ \f{ \ka (\mu) \vep }{ U_\ell - z_\ell } \ri ]^2, \] \[ \Lm_\ell (\mu, \nu) = \f{ \ka^2 (\mu) C_\ell }{  \nu  } \geq \f{C_\ell \ln
\f{1}{\de_\ell} }{\vep^2 n_\ell} \f{ w_\ell + \f{\ro}{n_\ell}  } { \nu } \li [ \f{ \ka (\mu) \vep }{ z_\ell - L_\ell } \ri ]^2.
\]
These two inequalities imply that $\liminf_{\vep \downarrow 0} \Lm_\ell (\mu, \nu) \geq 1$.  Hence, $\om \in \mscr{E}_5$ and thus
(\ref{cite369afull}) holds.

To show (\ref{cite369bfull}), it suffices to show that $\om \in \mscr{E}_6$ for $\om \in \mscr{E}_4$.  As a consequence of $\om \in \mscr{E}_2$,
we have that either $\ell = \tau_\vep$ or
\[
\Lm_{\ell - 1} (\mu, \nu) = \f{ \ka^2 (\mu) C_{\ell - 1} }{  \nu  }  < \f{C_{\ell - 1} \ln \f{1}{\de_{\ell-1}} }{\vep^2 n_{\ell - 1}}  \f{
w_{\ell - 1} + \f{\ro}{n_{\ell - 1}}  } { \nu } \li [ \f{ \ka (\mu) \vep }{ U_{\ell-1} - z_{\ell - 1} } \ri ]^2 \] or
\[ \Lm_{\ell -
1} (\mu, \nu) = \f{ \ka^2 (\mu) C_{\ell - 1} }{ \nu  } < \f{C_{\ell - 1} \ln \f{1}{\de_{\ell-1}} }{\vep^2 n_{\ell - 1}} \f{ w_{\ell - 1} +
\f{\ro}{n_{\ell - 1}} } { \nu } \li [ \f{ \ka (\mu) \vep }{ z_{\ell - 1} - L_{\ell-1} } \ri ]^2
\]
must be true.  This implies that either $\lim_{\vep \downarrow 0} \ell = 1$ or
\[
\limsup_{\vep \downarrow 0} \Lm_{\ell - 1} (\mu, \nu) \leq 1, \qqu \liminf_{\vep \downarrow 0} \ell > 1 \]
 must be true.  Hence, $\om \in \mscr{E}_6$ and thus (\ref{cite369bfull}) holds.

 \epf

 \subsection{Proof of Properties (III) -- (VII)}

Property (III) follows from Lemma \ref{lem22888} and the established Property (II).  To prove Properties (IV) -- (VII), we need to establish
Lemmas \ref{doublep} and \ref{viplast} in the sequel.  Once these two lemmas are established,

Property (IV) follows from Lemmas \ref{lem21888},  \ref{doublep} and \ref{viplast};

Property (V) follows from Lemmas \ref{lem23888},  \ref{doublep} and \ref{viplast};

Property (VI) follows from Lemmas  \ref{lem24888},  \ref{doublep} and \ref{viplast};

Property (VII) follows from Lemmas \ref{lem25888},  \ref{doublep} and \ref{viplast}.

\beL

\la{doublep}

If the random variable $X$ has mean $\mu$ and variance $\nu$ such that $\Lm_j (\mu, \nu)  > 1$ for some index $j \geq 1$, then there exist $\ga
> 0$ and $\vep^* > 0$ such that $\{ \bs{D}_\ell = 0 \} \subseteq
\{ | Z_\ell - \mu  | \geq \ga \; \tx{or} \; | V_{n_\ell} - \nu  | \geq \ga \}$ for all $\ell \geq j$ and $\vep \in (0, \vep^*)$.

Similarly, if the random variable $X$ has mean $\mu$ and variance $\nu$ such that  $\Lm_j (\mu, \nu)  < 1$ for some index $j \geq 1$, then there
exist $\ga > 0$ and $\vep^* > 0$ such that $\{ \bs{D}_\ell = 1 \} \subseteq \{ | Z_\ell - \mu  | \geq \ga  \; \tx{or} \; | V_{n_\ell} - \nu |
\geq \ga  \}$ for all $\ell \leq j$ and $\vep \in (0, \vep^*)$.

 \eeL

\bpf

Define \[ H_\ell(z, \se, \vep) = \f{n_\ell}{(\se + \f{\ro}{n_\ell} ) \ln \f{1}{\de_\ell} } \min \{ [z - \mscr{U} (z, \vep) ]^2, \; [z - \mscr{L}
(z, \vep) ]^2 \}
\]
 for $\ell \geq \tau_\vep$.  By the definitions of $H_\ell(z, \se, \vep)$, the
decision variable $\bs{D}_\ell$, and the stopping rule, we have that $\{ H_\ell( Z_\ell, V_{n_\ell}, \vep) < 1 \} \subseteq \{ \bs{D}_\ell = 0
\}$ and $\{ H_\ell( Z_\ell, V_{n_\ell}, \vep) \geq 1 \} \subseteq \{ \bs{D}_\ell = 1 \}$ for $\ell \geq \tau_\vep$.

Next, we claim that for any $\vsi \in (0, 1)$, there exist $\ga > 0$ and $\vep^*
> 0$ such that $(1 - \vsi) \Lm_\ell (z, \se) <  H_\ell(z, \se, \vep) < (1 + \vsi) \Lm_\ell (z, \se)$ for all $z \in (\mu - \ga, \mu + \ga), \;
\se \in (\nu - \ga, \nu + \ga), \;  \vep \in (0, \vep^*)$ and $\ell \geq \tau_\vep$.

To show this claim, note that for $\vsi \in (0, 1)$, there exists $\eta \in (0, 1)$ such that $\f{1 + \eta}{(1 - \eta) (1 + \vsi)}   < 1 < \f{1
- \eta}{(1 + \eta) (1 - \vsi)}$. Note that there exist $\ga > 0$ and $\vep_1 > 0$ such that
\bel &   &  1 - \eta < \li [ \f{z - \mscr{L} (z, \vep)  } {  \ka (z) \vep }  \ri ]^2 < 1 + \eta, \la{inq33883e}\\
&  &  1 - \eta < \li [ \f{ \mscr{U} (z, \vep) - z } {  \ka (z) \vep }  \ri ]^2 < 1 + \eta \la{inq33884f} \eel for all $z \in (\mu - \ga, \mu +
\ga)$ and $\vep \in (0, \vep_1)$.  Here, inequalities (\ref{inq33883e}) and (\ref{inq33884f})  are due to the continuity of $\ka(z)$ and the
assumption associated with (\ref{uniformasp}).  By the definition of the sample sizes, there exists $\vep^* \in (0, \vep_1)$ such that \be
\la{where8}
 \f{(1 - \eta) \se \vep^2}{C_\ell} <  \li ( \se +
\f{\ro}{n_\ell} \ri ) \f{ \ln \f{1}{\de_{\ell}} }{n_\ell}  < \f{(1 + \eta) \se \vep^2}{C_\ell} \ee for all $\se \in (\nu - \ga, \nu + \ga), \;
\vep \in (0, \vep^*)$ and $\ell \geq \tau_\vep$.  It follows from (\ref{inq33883e}), (\ref{inq33884f}) and  (\ref{where8}) that \bee  &  &
H_\ell(z, \se, \vep) < \f{ (1 + \eta) C_\ell [\ka(z) \vep]^2  }{ (1 - \eta) \se \vep^2}  = \f{ 1 + \eta }{ 1 - \eta } \Lm_\ell (z, \se) < (1 +
\vsi) \Lm_\ell (z, \se) \eee and \bee & & H_\ell(z, \se, \vep)
> \f{ (1 - \eta) C_\ell [\ka(z) \vep]^2  } { (1 + \eta) \se \vep^2} = \f{ 1 - \eta } { 1 + \eta } \Lm_\ell (z, \se) > (1 - \vsi)
\Lm_\ell (z, \se)  \eee for all $z \in (\mu - \ga, \mu + \ga), \; \se \in (\nu - \ga, \nu + \ga), \; \vep \in (0, \vep^*)$ and $\ell \geq
\tau_\vep$.  The claim is thus proved.

In the case that $\Lm_j (\mu, \nu) > 1$, let {\small $\vsi \in \li ( 0, 1 - \f{1}{ \sq{ \Lm_j (\mu, \nu) } } \ri )$}. By the continuity of
$\Lm_\ell (z, \se)$, there exists a number $\ga_0 > 0$ such that
\[
\Lm_j (z, \se) \geq (1 - \vsi) \Lm_j (\mu, \nu) \qu \tx{for all} \; z \in ( \mu - \ga_0, \mu + \ga_0 ), \qu \se \in ( \nu - \ga_0, \nu + \ga_0
).
\]
Since $\Lm_\ell(\mu, \nu)$ is increasing with respect to $\ell$, it must be true that
\[
\Lm_\ell (z, \se) \geq (1 - \vsi) \Lm_j (\mu, \nu) \qu \tx{for all} \; \; \ell \geq j  \; \tx{and} \; z \in ( \mu - \ga_0, \mu + \ga_0 ), \qu
\se \in ( \nu - \ga_0, \nu + \ga_0 ).
\]
By our established claim,  there exist $\ga \in (0, \ga_0)$ and $\vep^* > 0$ such that
\[
H_\ell(z, \se, \vep) > (1 - \vsi) \Lm_\ell (z, \se) \geq (1 - \vsi)^2 \Lm_j (\mu, \nu) > 1, \qqu \mscr{L} (z, \vep) \in \Se, \qqu \mscr{U} (z,
\vep) \in \Se
\]
for all $z \in (\mu - \ga, \mu + \ga), \; \se \in (\nu - \ga, \nu + \ga), \; \vep \in (0, \vep^*)$ and $\ell \geq j$.  Making use of this
result, we have that there exist $\ga > 0$ and $\vep^* > 0$ such that $\{ | Z_\ell - \mu  | < \ga, \; | V_{n_\ell} - \nu  | < \ga \} \subseteq
\{ \bs{D}_\ell = 1 \}$ for all $\vep \in (0, \vep^*)$ and $\ell \geq j$.  This implies that, there exist $\ga > 0$ and $\vep^* > 0$ such that
$\{ \bs{D}_\ell = 0 \} \subseteq \{ | Z_\ell - \mu  | \geq \ga \; \tx{or} \; | V_{n_\ell} - \nu  | \geq \ga \}$ for all $\vep \in (0, \vep^*)$
and $\ell \geq j$. Hence, we have shown the first statement of the lemma.

To show the second statement of the lemma, choose {\small $\vsi \in \li (0, \f{1}{ \sq{ \Lm_j (\mu, \nu) } } - 1 \ri )$} for $\mu, \nu$ with
$\Lm_j (\mu, \nu) < 1$.  By the continuity of $\Lm_\ell (z, \se)$, there exists a number $\ga_0 > 0$ such that
\[
\Lm_j (z, \se) \leq (1 + \vsi) \Lm_j (\mu, \nu) \qu \tx{for all} \; z \in ( \mu - \ga_0, \mu + \ga_0 ), \qu \se \in ( \nu - \ga_0, \nu + \ga_0
).
\]
Since $\Lm_\ell(\mu, \nu)$ is increasing with respect to $\ell$, it must be true that
\[
\Lm_\ell (z, \se) \leq (1 + \vsi) \Lm_j (\mu, \nu) \qu \tx{for all} \; \; \ell \leq j  \; \tx{and} \; z \in ( \mu - \ga_0, \mu + \ga_0 ), \qu
\se \in ( \nu - \ga_0, \nu + \ga_0 ).
\]
By our established claim,  there exist $\ga \in (0, \ga_0)$ and $\vep^* > 0$ such that
\[
H_\ell(z, \se, \vep) < (1 + \vsi) \Lm_\ell (z, \se) \leq (1 + \vsi)^2 \Lm_j (\mu, \nu) < 1, \qqu \mscr{L} (z, \vep) \in \Se, \qqu \mscr{U} (z,
\vep) \in \Se
\]
for all $z \in (\mu - \ga, \mu + \ga), \; \se \in (\nu - \ga, \nu + \ga), \; \vep \in (0, \vep^*)$ and $\ell \leq j$. Making use of this result,
we have that there exist $\ga
> 0$ and $\vep^* > 0$ such that $\{ | Z_\ell - \mu  | < \ga, \; | V_{n_\ell} - \nu  | < \ga \} \subseteq  \{ \bs{D}_\ell = 0  \}$ for all $\vep \in (0, \vep^*)$ and
$\ell \leq j$.  This implies that, there exist $\ga > 0$ and $\vep^* > 0$ such that $\{ \bs{D}_\ell = 1  \} \subseteq \{ | Z_\ell - \mu  | \geq
\ga \; \tx{or} \; | V_{n_\ell} - \nu  | \geq \ga \}$ for all $\vep \in (0, \vep^*)$ and $\ell \leq j$. So, we have shown the second statement of
the lemma. This completes the proof of the lemma.

\epf

\beL

\la{viplast}

If the random variable $X$ has mean $\mu$ and variance $\nu$ such that $\Lm_{i-1} (\mu, \nu) \leq 1 < \Lm_i (\mu, \nu)$ for some index $i \geq
1$, then \bee & & \lim_{\vep \downarrow 0} \sum_{\ell = \tau_\vep}^{i - 2} n_\ell \Pr \{ \bs{D}_\ell = 1 \} = 0,  \qqu \lim_{\vep \downarrow 0}
\sum_{\ell = i}^\iy n_\ell \Pr \{ \bs{D}_{\ell} = 0 \} = 0.  \eee Moreover, \be \la{last888}
 \lim_{\vep \downarrow 0} \sum_{\ell = \tau_\vep}^{i
- 1} n_\ell \Pr \{ \bs{D}_\ell = 1 \} = 0 \ee if the random variable $X$ has mean $\mu$ and variance $\nu$ such that $\Lm_{i-1} (\mu, \nu) < 1 <
\Lm_i (\mu, \nu)$ for some index $i \geq 1$.

\eeL

\bpf

From Lemma \ref{doublep},   we have that there exists $\ga > 0$ such that \bee  \Pr \{ \bs{D}_\ell = 1 \} & \leq & \Pr \{ |
Z_\ell - \mu | \geq \ga \; \tx{or} \; | V_{n_\ell} - \nu | \geq \ga \}\\
& \leq & \Pr \{  | Z_\ell - \mu | \geq \ga \} + \Pr \{  | V_{n_\ell} - \nu | \geq \ga \} \eee for $\tau_\vep \leq \ell \leq i - 2$, and that
\bee  \Pr \{ \bs{D}_\ell = 0 \} & \leq & \Pr \{  | Z_\ell - \mu | \geq \ga \; \tx{or}
\; | V_{n_\ell} - \nu | \geq \ga \}\\
& \leq & \Pr \{  | Z_\ell - \mu | \geq \ga \} + \Pr \{  | V_{n_\ell} - \nu | \geq \ga \} \eee for $\ell \geq i$. From Lemma \ref{vipnow} and
(\ref{mostvip}) of Lemma \ref{nonuniformBE}, we have that there exists an absolute  constant $C > 0$ such that \bee &  & \Pr \{  | Z_\ell - \mu
| \geq \ga \} \leq \exp \li ( - \f{n_\ell}{2} \f{ \ga^2 }{  \nu } \ri ) + \f{2 C }{n_\ell^2} \f{ \mscr{W} }{ \ga^3 } \eee and  \bee & & \Pr \{ |
V_{n_\ell} - \nu | \geq \ga  \}  \leq  \exp \li ( - \f{{n_\ell}}{4} \f{ \ga}{ \nu } \ri ) +  \exp \li ( - \f{{n_\ell}}{8} \f{ \ga^2 }{ \varpi }
\ri ) + \f{4 C }{{n_\ell}^2 \ga^3}  \li (  \sq{2} \mscr{W} \ga^{3\sh 2} +  4 \mscr{V} \ri ), \eee where \[ \nu = \bb{E} [ | X - \mu |^2 ] < \iy,
\qqu \mscr{W} = \bb{E} [ | X - \mu |^3 ] < \iy, \]
\[ \varpi = \bb{E} [ | ( X - \mu
)^2 - \nu |^2 ] < \iy, \qqu \mscr{V} = \bb{E} [ | ( X - \mu )^2 - \nu |^3 ] < \iy,
\]
as a consequence of $\bb{E} [ | X |^6 ] < \iy$.  Therefore, \bee &   & \sum_{\ell = \tau_\vep}^{i-2} n_\ell \Pr \{ \bs{D}_{\ell} = 1 \}\\
&  & \leq \sum_{\ell = \tau_\vep}^{i-2} n_\ell \li [ \exp \li ( - \f{{n_\ell}}{2} \f{ \ga^2 }{ \nu } \ri ) +  \exp \li ( - \f{{n_\ell}}{4} \f{
\ga}{ \nu } \ri ) +  \exp \li ( - \f{{n_\ell}}{8} \f{ \ga^2 }{ \varpi } \ri ) + \f{2 C }{{n_\ell}^2 \ga^3}  \li ( \mscr{W} + 2 \sq{2} \mscr{W}
\ga^{3\sh 2} + 8 \mscr{V} \ri )  \ri ] \to 0 \eee and \bee & & \sum_{\ell = i}^\iy n_\ell
\Pr \{ \bs{D}_{\ell} = 0 \}\\
&  & \leq \sum_{\ell = i}^\iy n_\ell \li [ \exp \li ( - \f{{n_\ell}}{2} \f{ \ga^2 }{ \nu } \ri ) +  \exp \li ( - \f{{n_\ell}}{4} \f{ \ga}{ \nu }
\ri ) +  \exp \li ( - \f{{n_\ell}}{8} \f{ \ga^2 }{ \varpi } \ri ) + \f{2 C }{{n_\ell}^2 \ga^3}  \li ( \mscr{W} + 2 \sq{2} \mscr{W} \ga^{3\sh 2}
+ 8 \mscr{V} \ri ) \ri ] \to 0 \eee  as $\vep \downarrow 0$. In a similar manner, we can establish (\ref{last888}). This completes the proof of
the lemma.

 \epf

\end{document}